\magnification 1200

\newif\ifuseauxfile \useauxfilefalse
\newif\ifshowrefs \showrefsfalse
\newif\ifshowlabels \showlabelsfalse

\newif\ifexistsauxfile \existsauxfilefalse
\newif\ifexistsauxdata \existsauxdatafalse
\newread\inauxfile \newwrite\outauxfile
\ifuseauxfile
\immediate\openin\inauxfile=\jobname.aux
\ifeof\inauxfile \existsauxfilefalse \showrefsfalse \else \existsauxfiletrue \fi
\immediate\closein\inauxfile
\fi

\catcode`\@=11 \newcount\el@ven \el@ven=11 \newcount\tw@lve \tw@lve=12
\catcode`\-=\el@ven \catcode`\_=\el@ven
\catcode`\0=\el@ven \catcode`\1=\el@ven \catcode`\2=\el@ven \catcode`\3=\el@ven
\catcode`\4=\el@ven \catcode`\5=\el@ven \catcode`\6=\el@ven \catcode`\7=\el@ven
\catcode`\8=\el@ven \catcode`\9=\el@ven

\existsauxdatatrue

\def\REFanal-sing{0.4}
\def\REFsing-ordering{0.5}
\def\REFlc-threshold{0.6}
\def\REFlct-semicontinuity{0.7}
\def\REFthm-strong-openness{0.8}
\def\REFusc-ideal-regul{0.9}
\def\REFcor-upper-regularization{0.10}
\def\REFcao-nd-current{0.12}
\def\REFhard-lefschetz{0.13}
\def\REFlefschetz-semipositive{0.14}
\def\REFkvn-cao-theorem{0.15}
\def\REFusc-ideal-regul-metric{0.16}
\def\REFsection-approx-psh{1}
\def\REFsection-approx-psh{1}
\def\REFlelong-num{1.1}
\def\REFpsh-bergman-approx{1.2}
\def\REFbuild-potential{1.5}
\def\REFcurrent-bergman-approx{1.6}
\def\REFequisingular-approx{1.7}
\def\REFint-decomp{1.8}
\def\REFright-hand-int-bergman{1.9}
\def\REFbergman-convergence-estimate{1.10}
\def\REFmultiplier-ideals-comparison{1.12}
\def\REFsection-strong-openness{2}
\def\REFsection-strong-openness{2}
\def\REFopenness-conj{2.1}
\def\REFthm-hiep{2.2}
\def\REFcor-strong-openness{2.3}
\def\REFcor-convergence-below{2.4}
\def\REFlemma-hiep{2.9}
\def\REFinequality-Fj{2.10}
\def\REFintegral-bound{2.12)}
\def\REFsection-hard-lefschetz{3}
\def\REFsection-hard-lefschetz{3}
\def\REFbochner-variant{3.1.1}
\def\REFbochner-variant-complete{3.1.2}
\def\REFharmonic-rep{3.2.3}
\def\REFnef-weird-example{3.2.4}
\def\REFsection-num-dim{4}
\def\REFsection-num-dim{4}
\def\REFsubsection-equising-approx{4.1}
\def\REFsing-equiv-fibration{4.1.1}
\def\REFmin-operation{4.1.2}
\def\REFasympt-equisingular-def{4.1.3}
\def\REFtheorem-bergman-approx-fund-observ{4.1.5}
\def\REFsing-comparison{4.1.6}
\def\REFprop-bergman-max-relation{4.1.8}
\def\REFtheorem-bergman-additivity{4.1.9}

\def\REFadd2{4.1.11}
\def\REFadd1{4.1.12}
\def\REFintegral-condition{4.1.13}
\def\REFadd3{4.1.14}
\def\REFsing-compare{4.1.15}
\def\REFweak-sing-relation{4.1.16}
\def\REFbergman-approx-map{4.1.18}
\def\REFweak-inequality{4.1.19}
\def\REFweak-equivalence{4.1.20}
\def\REFbergman-singularity-isomorphism{4.1.21}
\def\REFsubsection-intersection-theory{4.2}
\def\REFlog-res-current{4.2.1}
\def\REFdirect-image-equation{4.2.2}
\def\REFmonotonicity-lemma{4.2.3}
\def\REFintersection-product{4.2.4}
\def\REFnumerical-dimension{4.3.1}
\def\REFnumerical-dimension-line-bundle{4.3.2}
\def\REFsection-cao-proof{5}
\def\REFsection-cao-proof{5}
\def\REFcurvature-estimates{5.1}
\def\REFLtwo-estimate{5.3}
\def\REFerror-estimate{5.4}
\def\REFratio-lemma{5.5}
\def\REFproduct-eigenvalues-bound{5.6}
\def\REFcohomology-bound{6.5}
\def\REFtheorem-kummer{6.8}


\ifexistsauxfile \input \jobname.aux \existsauxdatatrue \fi
\ifuseauxfile \immediate\openout\outauxfile=\jobname.aux \fi
\catcode`\0=\tw@lve \catcode`\1=\tw@lve \catcode`\2=\tw@lve \catcode`\3=\tw@lve
\catcode`\4=\tw@lve \catcode`\5=\tw@lve \catcode`\6=\tw@lve \catcode`\7=\tw@lve
\catcode`\8=\tw@lve \catcode`\9=\tw@lve
\catcode`\-=\tw@lve \catcode`\_=8

\parindent=3.5mm
\vsize=21.4cm \voffset=-3.5mm
\hsize=14.5cm \hoffset=-5.5mm
\parskip=4pt plus 2pt minus 1pt
\parindent=5.2mm

\font\twentybf = cmbx10 at 20.736pt

\font\twelvebf = cmbx10 at 12pt

\font\twelvebsy=cmbsy10 at 12pt
\font\tenbsy=cmbsy10
\font\eightbsy=cmbsy8
\font\sevenbsy=cmbsy7
\font\sixbsy=cmbsy6
\font\fivebsy=cmbsy5


\font\tenCal=eusm10 at 10pt
\font\eightCal=eusm10 at 8pt
\font\sevenCal=eusm10 at 7pt
\font\sixCal=eusm10 at 6pt
\font\fiveCal=eusm10 at 5pt

\newfam\Calfam
  \textfont\Calfam=\tenCal
  \scriptfont\Calfam=\sevenCal
  \scriptscriptfont\Calfam=\fiveCal
\def\Cal{\fam\Calfam\tenCal}


\font\twelvebf=cmbx10 at 12pt

\font\teneuf=eufm10
\font\seveneuf=eufm7
\font\fiveeuf=eufm5
\newfam\euffam
  \textfont\euffam=\teneuf
  \scriptfont\euffam=\seveneuf
  \scriptscriptfont\euffam=\fiveeuf

\font\tenmsa=msam10 \font\sevenmsa=msam7 \font\fivemsa=msam5
\newfam\msafam
  \textfont\msafam=\tenmsa
  \scriptfont\msafam=\sevenmsa
  \scriptscriptfont\msafam=\fivemsa

\font\tenmsb=msbm10 \font\sevenmsb=msbm7 \font\fivemsb=msbm5
\newfam\msbfam
  \textfont\msbfam=\tenmsb
  \scriptfont\msbfam=\sevenmsb
  \scriptscriptfont\msbfam=\fivemsb
\def\Bbb{\fam\msbfam\tenmsb}

\font\twelvemib = cmmib10 at 12pt
\font\tenmib = cmmib10
\font\sevenmib = cmmib7

\font\eightrm = cmr8
\font\eightbf = cmbx8
\font\eighti = cmmi8
\font\eightit = cmti8
\font\eightsy = cmsy8
\font\eighttt = cmtt8
\font\eightex = cmex8
\font\eightmib = cmmib8
\font\eightmsa = msam8
\font\eightmsb = msbm8
\font\eighteuf = eufm8

\font\sixrm = cmr6
\font\sixbf = cmbx6
\font\sixi = cmmi6

\font\sixsy = cmsy6

\font\sixmib = cmmib6
\font\sixmsa = msam6
\font\sixmsb = msbm6
\font\sixeuf = eufm6

\font\fiverm = cmr5
\font\fivebf = cmbx5
\font\fivei = cmmi5

\font\fivesy = cmsy5

\font\fivemib = cmmib5
\font\fivemsa = msam5
\font\fivemsb = msbm5

\def\eightpoint{\def\rm{\fam0\eightrm}
   \textfont0=\eightrm \scriptfont0=\sixrm \scriptscriptfont0=\fiverm
   \textfont1=\eighti  \scriptfont1=\sixi  \scriptscriptfont1=\fivei
   \textfont2=\eightsy \scriptfont2=\sixsy \scriptscriptfont2=\fivesy
   \textfont3=\eightex \scriptfont3=\eightex\scriptscriptfont3=\eightex
   \def\it{\fam\itfam\eightit}%
   \textfont\itfam=\eightit
   \def\bf{\fam\bffam\eightbf}%
   \textfont\bffam=\eightbf \scriptfont\bffam=\sixbf
   \scriptscriptfont\bffam=\fivebf
   \def\tt{\fam\ttfam\eighttt}%
   \textfont\ttfam=\eighttt
   \textfont\msafam=\eightmsa \scriptfont\msafam=\sixmsa
   \textfont\msbfam=\eightmsb \scriptfont\msbfam=\sixmsb
   \textfont\euffam=\eighteuf \scriptfont\euffam=\sixeuf
   \textfont\Calfam=\eightCal
   \scriptfont\Calfam=\sixCal
   \normalbaselineskip=9.6pt
   \normalbaselines\rm}

\def\tenpointbf{%
 \textfont0=\tenbf   \scriptfont0=\sevenbf   \scriptscriptfont0=\fivebf
 \textfont1=\tenmib  \scriptfont1=\sevenmib  \scriptscriptfont1=\fivemib
 \textfont2=\tenbsy  \scriptfont2=\sevenbsy  \scriptscriptfont2=\fivebsy
 \tenbf
 \baselineskip=14.4pt}

\def\twelvepointbf{%
 \textfont0=\twelvebf   \scriptfont0=\eightbf   \scriptscriptfont0=\sixbf
 \textfont1=\twelvemib  \scriptfont1=\eightmib \scriptscriptfont1=\sixmib
 \textfont2=\twelvebsy  \scriptfont2=\eightbsy  \scriptscriptfont2=\sixbsy
 \twelvebf
 \baselineskip=14.4pt}

\def\proof#1{\removelastskip\vskip4pt\noindent{\it #1.}}
\def\endproof{\vskip5pt plus 1pt minus 2pt}

\catcode`\@=11
\def\hexnumber@#1{\ifnum#1<10 \number#1\else
 \ifnum#1=10 A\else\ifnum#1=11 B\else\ifnum#1=12 C\else
 \ifnum#1=13 D\else\ifnum#1=14 E\else\ifnum#1=15 F\fi\fi\fi\fi\fi\fi\fi}
\def\msa@{\hexnumber@\msafam}
\def\msb@{\hexnumber@\msbfam}
\mathchardef\compact="3\msa@62
\mathchardef\complement="0\msa@7B
\mathchardef\preccurlyeq="3\msa@34
\mathchardef\succcurlyeq="3\msa@3C
\mathchardef\square="0\msa@03
\mathchardef\smallsetminus="2\msb@72
\mathchardef\subsetneq="2\msb@28
\catcode`\@=12


\def\bC{{\Bbb C}}
\def\bN{{\Bbb N}}
\def\bP{{\Bbb P}}
\def\bQ{{\Bbb Q}}
\def\bR{{\Bbb R}}
\def\bZ{{\Bbb Z}}

\def\cC{{\Cal C}}
\def\cE{{\Cal E}}
\def\cF{{\Cal F}}
\def\cH{{\Cal H}}
\def\cI{{\Cal I}}
\def\cJ{{\Cal J}}
\def\cK{{\Cal K}}
\def\cO{{\Cal O}}
\def\cS{{\Cal S}}
\def\cU{{\Cal U}}

\def\bfB{{\bf B}}
\def\hatcS{\smash{\widehat{\Cal S}}}

\def\fraka{\hbox{\teneuf a}}
\def\frakm{\hbox{\teneuf m}}
\let\ssm=\smallsetminus
\let\em=\it

\def\nd{\mathop{\rm nd}}

\def\num{\mathop{\rm num}}
\def\ord{\mathop{\rm ord}\nolimits}
\def\ev{\mathop{\rm ev}}
\def\Alb{\mathop{\rm Alb}}
\def\Div{\mathop{\rm Div}\nolimits}
\def\PSH{\mathop{\rm PSH}\nolimits}
\def\Pic{\mathop{\rm Pic}}
\def\ii{\,{\rm{i}}\,}
\def\loc{{\rm loc}}
\def\bu{{\scriptstyle\bullet}}
\def\bfone{{\bf 1}}
\def\intn{\int\nolimits}
\def\QED{{\hfill$\square$}}
\def\dbar{{\overline\partial}}
\def\ddbar{{\partial\overline\partial}}
\def\ul#1{{\smash{\rm#1}}}
\def\SUPP{\mathop{\ul{SUPP}}}
\def\IM{\mathop{\ul{IM}}}
\def\IC{\mathop{\ul{IC}}}
\def\IT{\mathop{\ul{IT}}}

\def\cfudot#1{\ifmmode\setbox7\hbox{$\accent"5E#1$}\else
  \setbox7\hbox{\accent"5E#1}\penalty 10000\relax\fi\raise 1\ht7
  \hbox{\raise.1ex\hbox to 1\wd7{\hss.\hss}}\penalty 10000 \hskip-1\wd7\penalty
  10000\box7}

\def\RGBColor#1#2{#2}

\long\def\ignore#1{}
\long\def\cite#1{[#1]}

\ifshowlabels
\def\label#1#2{\catcode`\_=11\hbox{\RGBColor{0 .7 0}{#2 [[#1]]}}\catcode`\_=8{}
\ifuseauxfile
\immediate\write\outauxfile{\string\def\string\REF\string#1\string{%
#2\string}}#2\fi}
\else
\ifuseauxfile\def\label#1#2{%
\immediate\write\outauxfile{\string\def\string\REF\string#1\string{%
#2\string}}#2}\else \def\label#1#2{#2}\fi\fi

\ifexistsauxdata
\ifshowrefs
\def\ref#1{\catcode`\_=11\hbox{\RGBColor{1 0 0}{[[#1]]}}\catcode`\_=8{}}
\else
\def\ref#1{\catcode`\-=11\csname REF#1\endcsname{}\catcode`\-=12{}}
\fi
\else
\def\ref#1{\catcode`\_=11\hbox{\RGBColor{1 0 0}{[[#1]]}}\catcode`\_=8{}}
\fi

\def\frac#1#2{{#1\over #2}}

\long\def\claim#1 #2\endclaim
{\removelastskip\vskip7pt plus 1pt minus 1pt\noindent{\bf#1.}
{\it\ignorespaces#2}\vskip\baselineskip}

\newdimen\plainitemindent \plainitemindent=6mm
\def\item#1{\par\parindent=\plainitemindent\noindent
\hangindent\parindent\hbox to\parindent{#1\hss}\ignorespaces}


\newcount\void \void=-1
\newcount\newtitle \newtitle=\void
\newbox\titlebox   \setbox\titlebox\hbox{\hfil}
\newbox\sectionbox \setbox\sectionbox\hbox{\hfil}
\def\folio{\ifnum\pageno=\newtitle {\hfil}
           \else
           \ifodd\pageno
           {\eightpoint\hfil\copy\sectionbox\hfil\hfil\bf\number\pageno}\else
           {\eightpoint\bf\number\pageno\hfil\hfil\copy\titlebox\hfil}
           \fi\fi}
\footline={\hfil}
\headline={\folio}

\def\titlerunning#1{\setbox\titlebox\hbox{\eightpoint\hskip-25mm Jean-Pierre Demailly, #1}}
\def\title#1{%
  \newtitle=\pageno
  \vskip5mm
  \centerline{\twentybf\strut\kern4.25cm \vbox{\baselineskip=28pt #1}}\vskip1mm\noindent
  \titlerunning{#1}}

\def\sectionrunning#1#2{\setbox\sectionbox\hbox{\eightpoint #1. #2}}
\def\section#1#2{%
  \par\vskip0.75cm\penalty -100
  \vbox{\baselineskip=14.4pt\noindent{{\twelvepointbf #1.\kern6pt #2}}}
  \sectionrunning{#1}{#2}
  \vskip5pt
  \penalty 500}

\def\subsection#1#2{%
  \par\vskip0.2cm\penalty -100
  \vbox{\noindent{{\tenpointbf $\,$#1. #2}}}
  \vskip1pt
  \penalty 500}

\def\subsubsection#1#2{%
  \par\vskip0.1cm\penalty -100
  \noindent{{\tenpointbf #1.} {\it #2}}}
\def\\{\hfil\break}

\def\bibitem[#1]#2 #3\newblock#4\newblock#5\endbib
{\hangindent=1.6cm\hangafter=1
\noindent\rlap{\hbox{\eightbf[#1]}}\kern1.6cm{\rm #3}{\it #4}{\rm #5}}

\def\today{\ifcase\month\or
January\or February\or March\or April\or May\or June\or July\or August\or
September\or October\or November\or December\fi \space\number\day,
\number\year}


\title{On the cohomology of\\ pseudoeffective line bundles}

\centerline{\twelvebf Jean-Pierre Demailly}
\medskip

\centerline{Universit\'e de Grenoble~I, Institut Fourier, UMR 5582 du CNRS}
\centerline{BP$\,$74, 100 rue des Maths, 38402 Saint-Martin d'H\`eres, France}
\vskip10mm

\line{\it\hfill Dedicated to Professor Yum-Tong Siu on the occasion of his 70th birthday}
\vskip5mm

\noindent
{\eightpoint {\eightbf Abstract.} The goal of this survey is to present various results concerning the cohomology of pseudoeffective line bundles on compact K\"ahler manifolds, and related properties of their multiplier ideal sheaves. In case the curvature is strictly positive, the prototype is the well known Nadel vanishing theorem, which is itself a generalized analytic version of the fundamental Kawamata-Viehweg vanishing theorem of algebraic geometry. We are interested here in the case where the curvature is merely semipositive in the sense of currents, and the base manifold is not neces\-sarily projective. In this situation, one can still obtain interesting information on cohomology, e.g.\ a Hard Lefschetz theorem with pseudoeffective coefficients, in the form of a surjectivity statement for the Lefschetz map. More recently, Junyan Cao, in his PhD thesis defended in Grenoble, obtained a general K\"ahler vanishing theorem that depends on the concept of numerical dimension of a given pseudoeffective line bundle. The proof of these results depends in a crucial way on a general approximation result for closed $(1,1)$-currents, based on the use of Bergman kernels, and the related intersection theory of currents. Another important ingredient is the recent proof by Guan and Zhou of the strong openness conjecture. As an application, we discuss a structure theorem for compact K\"ahler threefolds without nontrivial subvarieties, following a joint work with F.~Campana and M.~Verbitsky. We hope that these notes will serve as a useful guide to the more detailed and more technical papers in the literature; in some cases, we provide here substantially simplified proofs and unifying viewpoints.
\medskip

\noindent
{\eightbf Key-words.}
Closed positive current, plurisubharmonic function, Ohsawa-Takegoshi extension theorem, curvature current, pseudoeffective line bundle, Bergman approximation, multiplier ideal sheaf, Nadel vanishing theorem, hard Lefschetz theorem, intersection theory, numerical dimension, openness conjecture, simple K\"ahler manifold, complex torus
\medskip

\noindent
{\eightbf MSC Classification 2010.} 14B05, 14F18, 14J30, 32C30, 32J25, 32L20
}

\section{0}{Introduction and statement of the main results}

Let $X$ be a compact K\"ahler $n$-dimensional manifold, equipped with a K\"ahler metric, i.e.\ a positive definite Hermitian $(1,1)$-form $\omega=\ii\sum_{1\le j,k\le n}
\omega_{jk}(z)\,dz_j\wedge d\overline z_k$ such that $d\omega=0$. By definition 
a holomorphic line bundle $L$ on $X$ is said to be {\it pseudo\-effective} if there exists a singular hermitian metric $h$ on $L$, given by $h(z)=e^{-\varphi(z)}$ with respect to a local trivialization $L_{|U}\simeq U\times\bC$,
such that the curvature form
$$
\ii\Theta_{L,h}:=\ii\ddbar\varphi
\leqno(0.1)$$
is (semi)positive in the sense of currents, i.e.\ $\varphi$ is locally
integrable and $i\Theta_{L,h}\ge 0\,$: in other words, the weight
function $\varphi$ is plurisubharmonic (psh) on the corresponding trivializing
open set~$U$. A basic concept is the notion of {\it multiplier ideal
sheaf}, introduced in \cite{Nad90}.

\claim{0.2.\ Definition} To any psh function $\varphi$ on an open subset $U$ of
a complex manifold~$X$, one associates the ``multiplier ideal sheaf'' $\cI(\varphi)\subset\cO_{X|U}$ of germs of
holomorphic functions $f\in\cO_{X,x}$, $x\in U$, such that $|f|^2e^{-\varphi}$
is integrable with respect to the Lebesgue measure in some local
coordinates near~$x$. We also define the global multiplier ideal sheaf $\cI(h)\subset\cO_X$ of a hermitian metric $h$ on $L\in\Pic(X)$ to be equal to $\cI(\varphi)$ on any open subset $U$ where $L_{|U}$ is trivial and $h=e^{-\varphi}$. In such a definition, we may in fact assume $\ii\Theta_{L,h}\geq -C\omega$,
i.e.\ locally $\varphi=\hbox{psh}+C^\infty$, we say in that case that
$\varphi$ is {\rm quasi-psh}.
\endclaim

Let us observe that a multiplier ideal sheaf $\cI(\varphi)$ is left unmodified
by adding a smooth function to $\varphi\,$; for such purposes, the 
additional $C^\infty$ terms are irrelevant in quasi-psh functions.
A crucial and well-known fact is that the ideal sheaves $\cI(\varphi)\subset
\cO_{X|U}$ and $\cI(h)\subset\cO_X$ are always {\it coherent analytic 
sheaves}; when $U\subset X$ is a coordinate open ball, this can be shown
by observing that $\cI(\varphi)$ coincides with the locally stationary limit 
$\cJ=\lim{\uparrow}_{N\to+\infty}\cJ_N$
of the increasing sequence of coherent ideals $\cJ_N=(g_j)_{0\le j<N}$
associated with a Hilbert basis
$(g_j)_{j\in\bN}$ of the Hilbert space of holomorphic functions
$f\in\cO_X(U)$ such that $\int_U|f|^2e^{-\varphi}dV_\omega<+\infty$.
The proof is a consequence of H\"ormander's $L^2$ estimates applied to 
weights of the form 
$$\psi(z)=\varphi(z)+(n+k)\log|z-x|^2.$$
This easily shows that $\cI(\varphi)_x+\frakm_x^k=\cJ_x+\frakm_x^k$, and
one then concludes that $\cI(\varphi)_x=\cJ_x$ by the Krull lemma.
When $X$ is {\it projective algebraic}, Serre's GAGA theorem implies that
$\cI(h)$ is in fact a {\it coherent algebraic sheaf}, in spite of the fact
that $\varphi$ may have very ``wild'' analytic singularities -- e.g.\
they might be everywhere dense in $X$ in the Euclidean
topology. Therefore, in some sense, the multiplier ideal sheaf is a
powerful tool to extract algebraic (or at least analytic) data from
arbitrary singularities of psh functions. In this context, assuming
strict positivity of the curvature, one has the following well-known
fundamental vanishing theorem.

\claim{\label{nadel}{0.3}.\ Theorem} {\rm (Nadel Vanishing Theorem, \cite{Nad90}, \cite{Dem93b})} Let
$(X,\omega)$ be a compact K\"ahler $n$-dimensional manifold, and let $L$ be a
holomorphic line bundle over $X$ equipped with a singular Hermitian 
metric~$h$. Assume that $\ii\Theta_{L,h}\ge\varepsilon\omega$
for some $\varepsilon>0$ on~$X$. Then
$$H^q\big(X,\cO(K_X\otimes L)\otimes\cI(h)\big)=0\qquad
\hbox{for all $q\ge 1$},$$
where $K_X=\Omega^n_X=\Lambda^nT^*_X$ denotes the canonical line bundle.
\endclaim

The proof follows from an application of H\"ormander's $L^2$ estimates
with singular weights, themselves derived from the Bochner-Kodaira identity
(see \cite{H\"or66}, \cite{Dem82}, \cite{Dem92}).
One should observe that the strict positivity assumption implies $L$ to be
big, hence $X$ must be projective, since every compact manifold that 
is K\"ahler and Moishezon is also projective (cf.~\cite{Moi66}, \cite{Pet86}, \cite{Pet98a}). However, when
relaxing the strict positivity assumption, one can enter the world of
general compact K\"ahler manifolds, and their study is one of our main goals.

In many cases, one has to assume that the psh functions involved have milder singularities. We say that a psh or quasi-psh function $\varphi$ has {\it analytic singularities} if locally on the domain of definition $U$ of $\varphi$ one can write

$$\varphi(z)=c\log\sum_{j=1}^N|g_j|^2+O(1)
\leqno(\label{anal-sing}{0.4})$$
where the $g_j$'s are holomorphic functions, $c\in\bR_+$ and $O(1)$ means
a locally bounded remainder term. Assumption (\ref{anal-sing}) implies that
the set of poles $Z=\varphi^{-1}(-\infty)$ is an analytic set, locally defined
as $Z=\bigcap g_j^{-1}(0)$, and that $\varphi$ is locally bounded on $U\ssm Z$.
We also refer to this situation by saying that $\varphi$ has 
{\it logarithmic poles}.
In general, one introduces the following comparison relations for psh
or quasi-psh functions $\varphi$ and hermitian metrics $h=e^{-\varphi}\,$;
a more flexible comparison relation will be introduced in 
Section~\ref{section-num-dim}.

\claim{\label{sing-ordering}{0.5}.\ Definition}
{\it Let $\varphi_1,$ $\varphi_2$ be psh functions on an
open subset $U$ of a complex manifold~$X$. We say that 
\item{\rm(a)} $\varphi_1$ has less singularities than
$\varphi_2$, and write $\varphi_1\preccurlyeq\varphi_2$, if for every point 
$x\in U$, there exists a neighborhood $V$ of $x$ and a constant $C\geq 0$ 
such that $\varphi_1\ge \varphi_2-C$ on~$V$. 
\vskip2pt
\item{\rm(b)} $\varphi_1$ and $\varphi_2$ have
equivalent singularities, and write $\varphi_1\sim\varphi_2$, if 
locally near any point of $U$ we have $\varphi_1-C\le \varphi_2\le \varphi_1+C$.
\vskip2pt\noindent
Similarly, given a pair of hermitian metrics $h_1$, $h_2$ on a line bundle 
$L\to X$,
\item{\rm(a')} we say that $h_1$ is less singular than $h_2$, and
write $h_1\preccurlyeq h_2$, if locally there exists a constant $C>0$ such that
$h_1\le C h_2$.
\vskip2pt
\item{\rm(b')} we say that $h_1$, $h_2$ have equivalent singularities, and
write $h_1\sim h_2$, if locally there exists a constant
$C>0$ such that $C^{-1}h_2\le h_1\le C h_2$.
\vskip2pt\noindent
$($of course when $h_1$ and $h_2$ are defined on a compact manifold $X$, 
the constant $C$ can be taken global on $X$ in {\rm (a')} and  {\rm(b')}$)$.}
\endclaim

Important features of psh singularities are the semi-continuity theorem (see \cite{DK01}) and the strong openness property recently proved by Guan and Zhou \cite{GZ13}, \cite{GZ14a}, \cite{GZ14b}. Let $U$ be an open set in a complex manifold $X$ and $\varphi$ a psh function on~$U$. Following \cite{DK01}, we define the log canonical threshold of $\varphi$ at a point $z_0\in U$ by
$$
c_{z_0}(\varphi) = \sup\big\{ c>0:\ e^{-2c\,\varphi} \hbox{ is } L^1 \hbox{ on a neighborhood of }z_0\big\}\in {}]0,+\infty]
\leqno(\label{lc-threshold}{0.6})
$$
(Here $L^1$ integrability refers to the Lebesgue measure with respect to local coordinates). It is an important invariant of the singularity of $\varphi$ at $z_0$. We refer to \cite{FEM03}, \cite{DH14}, \cite{DK01}, \cite{FEM10}, \cite{Kis94}, \cite{Nad90}, \cite{PS00}, \cite{Sko72b} for further information about properties of the log canonical threshold. In this setting, the semi-continuity theorem can be stated as follows.

\claim{\label{lct-semicontinuity}{0.7}.\ Theorem} {\rm(cf.\ \cite{DK01})} For any given $z_0\in U$, the map $\PSH(U)\to{}]0,+\infty]$, \hbox{$\varphi\mapsto c_{z_0}(\varphi)$} is upper semi-continuous with respect to the topology of weak con\-vergence on the space of psh functions $($the latter topology being actually the same as the topology of $L^1_{\rm loc}$ convergence$)$.
\endclaim

The original proof of \cite{DK01} was rather involved and depended on uniform polynomial approximation, combined with a reduction to a semi-continuity theorem for algebraic singularities; the Ohsawa-Takegoshi $L^2$ extension theorem \cite{OT87} was used in a crucial way. We will give here a simpler and more powerful derivation due to Hiep \cite{Hiep14}, still depending on the Ohsawa-Takegoshi theorem, that simultaneously yields effective versions of Berndtsson's result \cite{Bern13} on the openness conjecture, as well as Guan and Zhou's proof of the strong openness conjecture for multiplier ideal sheaves.

\claim{\label{thm-strong-openness}{0.8}.\ Theorem} {\rm(\cite{GZ13}, \cite{GZ14a}, \cite{GZ14b})} Let $\varphi,\,\psi_j$, $j\in\bN$, be psh functions on an open set $U$ in a complex manifold~$X$. Assume that $\psi_j\le\varphi$ and that $\psi_j$ converges to $\varphi$ in $L^1_{\rm loc}$ topology as $j\to+\infty$. Then for every relatively compact subset $U'\compact U$, the multiplier ideal sheaves $\cI(\psi_j)$ coincide with $\cI(\varphi)$ on $U'$ for $j\ge j_0(U')\gg 1$.
\endclaim

Before going further, notice that the family of multiplier ideals $\lambda\mapsto\cI(\lambda\varphi)$ associated with a psh function $\varphi$ is nonincreasing in $\lambda\in\bR_+$. By the Noetherian property of ideal sheaves, they can jump only for a locally finite set of values $\lambda$ in $[0,+\infty[$, and in particular, there exists a real value $\lambda_0>1$ such that
$$
\cI_+(\varphi):=\lim_{\varepsilon\to 0+}\cI((1+\varepsilon)\varphi)=\cI(\lambda\varphi),~~~\forall\lambda\in{}]1,\lambda_0].
\leqno(\label{usc-ideal-regul}{0.9})
$$
We will say that $\cI_+(\varphi)$ is the upper semicontinuous regularization of the multiplier ideal sheaf. Berndtsson's result \cite{Bern13} states that the equality $\cI(\varphi)=\cO_X$ implies $\cI_+(\varphi)=\cO_X$. If we take $\psi_j=(1+1/j)\varphi$ and assume (without loss of generality) that~$\varphi\le 0$, Theorem~\ref{thm-strong-openness} implies in fact

\claim{\label{cor-upper-regularization}{0.10}.\ Corollary} For every psh function $\varphi$, the upper semicontinuous regularization coin\-cides with the multiplier ideal sheaf, i.e.\ $\cI_+(\varphi)=\cI(\varphi)$.
\endclaim

Now, if $L$ is a pseudoeffective line bundle, it was observed in \cite{Dem00}
that there always exist a unique
equivalence class $h_{\min}$ of singular hermitian metrics with minimal 
singularities, such that $\ii\Theta_{L,h_{\min}}\geq 0$ (by this we mean 
that $h_{\min}$ is unique up to equivalence of singularities). In fact, 
if $h_\infty$ is a smooth metric on $L$, one can define the corresponding
weight $\varphi_{\min}$ of $h_{\min}$ as an upper envelope
$$
\varphi_{\min}(z)=\sup\big\{\varphi(z)\,;\;\ii\Theta_{L,h_\infty}+\ii\ddbar\varphi\ge 0,~\varphi\leq 0~\hbox{on}~X\big\},
\leqno(0.11)$$
and put $h_{\min}=h_\infty e^{-\varphi_{\min}}$. In general, $h_{\min}$ need not have analytic singularities.

An important fact is that one can approximate arbitrary psh functions by psh functions with analytic singularities. The appropriate technique consists of using an asymptotic Bergman kernel procedure (cf.\ \cite{Dem92} and Section~\ref{section-approx-psh}). If $\varphi$ is a holomorphic function on a ball $B\subset\bC^n$, one puts
$$
\varphi_m(z)=\frac{1}{2m}\log\sum_{\ell\in\bN}|g_{m,\ell}(z)|^2
$$
where $(g_{m,\ell})_{\ell\in\bN}$ is a Hilbert basis of the space 
$\cH(B,m\varphi)$ of $L^2$ holomorphic functions on $B$ such that
$\int_B|f|^2e^{-2m\varphi}dV<+\infty$. When $T=\alpha+dd^c\varphi$ is a
closed $(1,1)$-current on $X$ in the same cohomology class
as a smooth $(1,1)$-form $\alpha$ and $\varphi$ is a quasi-psh potential
on~$X$, a sequence of global approximations $T_m$ can be produced by 
taking a finite covering of $X$
by coordinate balls $(B_j)$. A partition of unity argument allows
to glue the local approximations $\varphi_{m,j}$ of $\varphi$ on $B_j$ into
a global potential $\varphi_m$, and one sets $T_m=\alpha+dd^c\varphi_m$.
These currents $T_m$ converge weakly to~$T$, are smooth in the complement
$X\ssm Z_m$ of an increasing family of analytic subsets $Z_m\subset X$, and
their singularities approach those of~$T$. More precisely,
the Lelong numbers $\nu(T_m,z)$ converge uniformly to those of $T$, and
whenever $T\ge 0$, it is possible to produce a current $T_m$ that only 
suffers a small loss of positivity, namely $T_m\ge -\varepsilon_m\omega$ 
where $\lim_{m\to+\infty}\varepsilon_m=0$. These considerations lead in a 
natural way to the concept of {\it numerical dimension} of a 
closed positive $(1,1)$-current $T$. We define
$$\nd(T)=\max\big\{p=0,1,\ldots,n\,;\;\limsup_{m\to+\infty}\int_{X\ssm Z_m}(T_m+\varepsilon_m\omega)^p\wedge\omega^{n-p}>0\big\}.
\leqno(\label{cao-nd-current}{0.12})$$
One can easily show (see Section~\ref{section-num-dim}) that the right hand side of (\ref{cao-nd-current}) does not depend
on the sequence~$(T_m)$, provided that the singularities approach those
of $T$ (we call this an ``asymptotically equisingular approximation'').

These concepts are very useful to study cohomology groups
with values in pseudoeffective line bundles $(L,h)$.
Without assuming any strict posi\-tivity of the curvature, one can obtain 
at least a hard Lefschetz theorem with coefficients in $L$. The technique 
is based on a use of harmonic forms with respect to suitable 
``equi\-singular approximations'' 
$\varphi_m$ of the weight $\varphi$ of $h$ (in that 
case we demand that $\cI(\varphi_m)=\cI(\varphi)$ for all~$m$); the main
idea is to work with complete K\"ahler metrics in the open complements 
$X\ssm Z_m$ where $\varphi_m$ is smooth, and to apply a variant of
the Bochner formula on these sets. More details can be found in 
Section~\ref{section-hard-lefschetz} and in \cite{DPS01}.

\claim{\label{hard-lefschetz}{0.13}.\ Theorem} {\rm(\cite{DPS01})}
Let $(L,h)$ be a pseudo-effective line bundle on a
compact K\"ahler manifold $(X,\omega)$ of dimension $n$, let
$\Theta_{L,h}\ge 0$ be its curvature current and $\cI(h)$ the
associated multiplier ideal sheaf. Then, the wedge multiplication
operator $\omega^q\wedge\bu$ induces a surjective morphism
$$
\Phi^q_{\omega,h}:
H^0(X,\Omega_X^{n-q}\otimes L\otimes\cI(h))\longrightarrow
H^q(X,\Omega_X^n\otimes L\otimes\cI(h)).
$$
\endclaim

The special case when $L$ is nef is due to Takegoshi \cite{Tak97}. An even
more special case is when $L$ is semipositive, i.e.\ possesses a
smooth metric with semipositive curvature.
In that case the multiple ideal sheaf $\cI(h)$ coincides with $\cO_X$
and we get the following consequence already observed by Enoki \cite{Eno93}
and Mourougane \cite{Mou95}.

\claim{\label{lefschetz-semipositive}{0.14}.\ Corollary}
Let $(L,h)$ be a semipositive line bundle
on a compact K\"ahler mani\-fold $(X,\omega)$ of dimension~$n$. Then, the wedge
multiplication operator $\omega^q\wedge\bu$ induces a surjective morphism
$$
\Phi^q_\omega:H^0(X,\Omega_X^{n-q}\otimes L)\longrightarrow
H^q(X,\Omega_X^n\otimes L).
$$
\endclaim

It should be observed that although all objects involved in Th.~\ref{hard-lefschetz}
are algebraic when $X$ is a projective manifold, there is no known
algebraic proof of the statement; it is not even clear how to define
algebraically $\cI(h)$ for the case when $h=h_{min}$ is a metric
with minimal singularity. However, even in the special circumstance
when $L$ is nef, the multiplier ideal sheaf is crucially needed.

Our next statement is taken from the PhD thesis of Junyan Cao \cite{JC13}. The proof is a combination of our Bergman regularization techniques, together with an argument of Ch.~Mourou\-gane \cite{Mou95} relying on a use of the Calabi-Yau theorem for Monge-Amp\`ere equations.

\claim{\label{kvn-cao-theorem}{0.15}.\ Theorem} {\rm(\cite{JC13}, \cite{JC14})}
Let $(L,h)$ be a pseudoeffective
line bundle on a compact K\"ahler $n$-dimensional manifold~$X$. Then
$$ H^q(X, K_{X}\otimes L\otimes\cI(h))=0\qquad
\hbox{for every}~~q\geq n-\nd (L,h)+1 ,$$
where $\nd (L,h):=\nd(\ii\Theta_{L,h})$.
\endclaim

Cao's technique of proof actually yields the result for the upper semicontinuous regularization
$$
\cI_+(h)=\lim_{\varepsilon\to 0}\cI(h^{1+\varepsilon})
\leqno(\label{usc-ideal-regul-metric}{0.16})
$$
instead of $\cI(h)$, but we can apply Guan-Zhou's Theorem~\ref{thm-strong-openness} to see that the equality $\cI_+(h)=\cI(h)$ always holds. As~a final geometric application of this circle of ideas, we present the following result which was obtained in \cite{CDV13}.

\claim{0.17.\ Theorem}  {\rm (\cite{CDV13})} Let $X$ be a compact K\"ahler threefold
that is ``strongly simple'' in the sense that $X$ has no nontrivial 
analytic subvariety. Then the Albanese morphism \hbox{$\alpha:X\to\Alb(X)$}
is a biholomorphism, and therefore $X$ is biholomorphic to a $3$-dimensional
complex torus $\bC^3/\Lambda$.
\endclaim

I would like to thank the referee wholeheartedly for numerous suggestions that led to substantial improvements of the exposition.

\section{\label{section-approx-psh}{1}}{Approximation of psh functions and of closed (1,1)-currents}

We first recall here the basic result on the approximation of psh functions
by psh functions with analytic singularities. The main idea
is taken from [Dem92] and relies on the Ohsawa-Takegoshi extension theorem,
For other applications to algebraic geometry,
see [Dem93b] and Demailly-Koll\'ar [DK01]. Let $\varphi$ be a
psh function on an open set $\Omega\subset\bC^n$. Recall that
the Lelong number of $\varphi$ at a point $x_0\in\Omega$ is defined to be
$$\nu(\varphi,x_0)=\liminf_{z\to x_0}\frac{\varphi(z)}{\log|z-x_0|}=
\lim_{r\to 0_+}\frac{\sup_{B(x_0,r)}\varphi}{\log r}.
\leqno(\label{lelong-num}{1.1})$$
In particular, if $\varphi=\log|f|$ with $f\in\cO(\Omega)$, then
$\nu(\varphi,x_0)$ is equal to
the vanishing order
$$
\ord_{x_0}(f)=\sup\{k\in\bN\,;D^{\alpha}f(x_0)=0,~
\forall|\alpha|<k\}.
$$
\vskip0pt

\claim{\label{psh-bergman-approx}{1.2}.\ Theorem}
Let $\varphi$ be a plurisubharmonic function on a
bounded pseudoconvex open set $\Omega\subset\bC^n$. For every $m>0$, let
$\cH_\Omega(m\varphi)$ be the Hilbert space of holomorphic functions $f$
on $\Omega$ such that $\int_\Omega|f|^2e^{-2m\varphi}dV_{2n}<+\infty$ and
let $\varphi_m=\frac{1}{2m}\log\sum|g_{m,\ell}|^2$ where $(g_{m,\ell})$
is an orthonormal basis of~$\cH_\Omega(m\varphi)$. Then there are
constants $C_1,C_2>0$ independent of $m$ such that
\vskip2pt
\item{\rm (a)} $\displaystyle\varphi(z)-\frac{C_1}{m}\le
\varphi_m(z)\le\sup_{|\zeta-z|<r}\varphi(\zeta)+\frac{1}{m}
\log\frac{C_2}{r^n}$ for every $z\in\Omega$ and $r<d(z,\partial\Omega)$. In particular,
$\varphi_m$ converges to $\varphi$ pointwise and in $L^1_{\rm loc}$ topology
on~$\Omega$ when $m\to+\infty$ and
\vskip7pt
\item{\rm (b)} $\displaystyle\nu(\varphi,z)-\frac{n}{m}\le\nu(\varphi_m,z)\le
\nu(\varphi,z)$~ for every $z\in\Omega$.
\endclaim

\proof{Proof}  (a) Note that $\sum|g_{m,\ell}(z)|^2$ is the square of
the norm
of the evaluation linear form $\ev_z:f\mapsto f(z)$ on $\cH_\Omega(m\varphi)$,
since $g_{m,\ell}(z)=\ev_z(g_{m,\ell})$ is the $\ell$-th coordinate of $\ev_z$
in the orthonormal basis $(g_{m,\ell})$. In other words, we have
$$\sum|g_{m,\ell}(z)|^2=\sup_{f\in B(1)}|f(z)|^2$$
where $B(1)$ is the unit ball of $\cH_\Omega(m\varphi)$ (The sum
is called the {\it Bergman kernel} associated with $\cH_\Omega(m\varphi)$).
As $\varphi$ is locally bounded from above, the $L^2$ topology is
actually stronger
than the topology of uniform convergence on compact subsets of~$\Omega$.
It follows that the series $\sum|g_{m,\ell}|^2$ converges uniformly on
$\Omega$ and that its sum is real analytic. Moreover, by what we just
explained, we have
$$\varphi_m(z)=\sup_{f\in B(1)}\frac{1}{2m}\log|f(z)|^2
=\sup_{f\in B(1)}\frac{1}{m}\log|f(z)|.$$
For $z_0\in\Omega$ and $r<d(z_0,\partial\Omega)$, the mean value
inequality applied to the psh function $|f|^2$ implies
$$\eqalign{
|f(z_0)|^2&\le\frac{1}{\pi^nr^{2n}/n!}\int_{|z-z_0|<r}
|f(z)|^2dV_{2n}(z)\cr
&\le\frac{1}{\pi^nr^{2n}/n!}\exp\Big(2m\sup_{|z-z_0|<r}\varphi(z)\Big)
\int_\Omega|f|^2e^{-2m\varphi}dV_{2n}.\cr}
$$
If we take the supremum over all $f\in B(1)$ we get
$$
\varphi_m(z_0)\le\sup_{|z-z_0|<r}\varphi(z)+\frac{1}{2m}
\log\frac{1}{\pi^nr^{2n}/n!}
$$
and the second inequality in (a) is proved -- as we see, this is an easy
consequence of the mean value inequality. Conversely, the
Ohsawa-Takegoshi extension theorem (\cite{OT87}) applied
to the $0$-dimensional subvariety $\{z_0\}\subset\Omega$ shows that for any
$a\in\bC$ there is a holomorphic function $f$ on $\Omega$ such that
$f(z_0)=a$ and
$$\int_\Omega|f|^2e^{-2m\varphi}dV_{2n}\le C_3|a|^2e^{-2m\varphi(z_0)},$$
where $C_3$ only depends on $n$ and~${\rm diam}\,\Omega$. We fix $a$
such that the right hand side is~$1$. Then $\Vert f\Vert \le 1$ and so
we get
$$\varphi_m(z_0)\ge\frac{1}{m}\log|f(z_0)|=
\frac{1}{m}\log|a|=\varphi(z)-\log\frac{C_3}{m}.$$
The inequalities given in (a) are thus proved. Taking $r=1/m$, we find
that
$$\lim_{m\to+\infty}\sup_{|\zeta-z|<1/m}\varphi(\zeta)=\varphi(z)$$
by the upper semicontinuity of $\varphi$, and therefore
$\lim\varphi_m(z)=\varphi(z)$, since $\lim\frac{1}{m}\log (C_2m^n)=0$.
\medskip

\noindent
(b) The above estimates imply
$$
\sup_{|z-z_0|<r}\varphi(z)-\frac{C_1}{m}\le
\sup_{|z-z_0|<r}\varphi_m(z)\le\sup_{|z-z_0|<2r}\varphi(z)+
\frac{1}{m}\log\frac{C_2}{r^n}.$$
After dividing by $\log r<0$ when $r\to 0$, we infer
$$
\frac{\sup_{|z-z_0|<2r}\varphi(z)+\frac{1}{m}\log\frac{C_2}{r^n}}{\log r}\le
\frac{\sup_{|z-z_0|<r}\varphi_m(z)}{\log r}\le
\frac{\sup_{|z-z_0|<r}\varphi(z)-\frac{C_1}{m}}{\log r},
$$
and from this and definition (\ref{lelong-num}), it follows immediately that
$$
\nu(\varphi,x)-\frac{n}{m}\le\nu(\varphi_m,z)\le\nu(\varphi,z).
\eqno\square
$$
\endproof

Theorem~\ref{psh-bergman-approx} implies in a straightforward manner
the deep result of \cite{Siu74} on the analyticity of the Lelong number
upperlevel sets.

\claim{1.3.\ Corollary}  {\rm \cite{Siu74}} Let $\varphi$ be a plurisubharmonic
function on a complex manifold~$X$. Then, for every $c>0$, the Lelong number
upperlevel set
$$
E_c(\varphi)=\big\{z\in X\,;\;\nu(\varphi,z)\ge c\big\}
$$
is an analytic subset of~$X$.
\endclaim

\proof{Proof}  Since analyticity is a local property, it
is enough to consider the case of a psh function $\varphi$ on a
pseudoconvex open set $\Omega\subset\bC^n$. The
inequalities obtained in Theorem 13.2~(b) imply that
$$
E_c(\varphi)=\bigcap_{m\ge m_0}E_{c-n/m}(\varphi_m).
$$
Now, it is clear that $E_c(\varphi_m)$ is the analytic set defined by the
equations $g_{m,\ell}^{(\alpha)}(z)=0$ for all multi-indices $\alpha$
such that~$|\alpha|<mc$. Thus $E_c(\varphi)$ is analytic as a (countable)
intersection of analytic sets.\QED
\endproof

\claim{1.4.\ Remark} {\rm  It has been observed by Dano Kim \cite{Kim13} that the functions
$\varphi_m$ produced by Th.~\ref{psh-bergman-approx} do not in general
satisfy $\varphi_{m+1}\succcurlyeq\varphi_m$, in other words their singularities
may not always increase monotonically to those of $\varphi$. Thanks to 
the subbadditivity result of 
\cite{DEL00}, this is however the case for any subsequence
$\varphi_{m_k}$ such that $m_k$ divides $m_{k+1}$, e.g.\ $m_k=2^k$ or
$m_k=k!$ (we will refer to such a sequence below as being a 
``{\it multiplicative sequence}''). In that case, a use of the 
Ohsawa-Takegoshi theorem on the diagonal of $\Omega\times\Omega$ shows that one
can obtain $\varphi_{m_{k+1}}\le \varphi_{m_k}$ (after possibly replacing
$\varphi_{m_k}$ by $\varphi_{m_k}+C/m_k$ with $C$ large enough), see
\cite{DEL00} and \cite{DPS01}.}
\endclaim

Our next goal is to study the regularization process more globally,
i.e.\ on a compact complex manifold~$X$. For this, we have to take care
of cohomology class. It is convenient to introduce
$d^c=\frac{\ii}{4\pi}(\dbar-\partial)$, so that $dd^c=\frac{\ii}{2\pi}\ddbar$.
Let $T$ be a closed $(1,1)$-current on~$X$. We assume that $T$ is
{\it quasi-positive}, i.e.\ that there exists a $(1,1)$-form $\gamma$
with continuous coefficients such that~$T\ge\gamma\;$; observe that
a function $\varphi$ is quasi-psh iff its complex Hessian is bounded below
by a $(1,1)$-form with continuous or locally bounded coefficients, that is, 
if $dd^c\varphi$ is quasi-positive. The case of positive currents 
($\gamma=0$) is of course the most important.

\claim{\label{build-potential}{1.5}.\ Lemma} There exists a smooth closed 
$(1,1)$-form $\alpha$ representing the same
$\ddbar$-cohomology class as $T$ and an {\rm quasi-psh} function $\varphi$
on $X$ such that $T=\alpha+dd^c\varphi$.
\endclaim

\proof{Proof}  Select
an open covering $(B_j)$ of $X$ by coordinate balls such that
$T=dd^c\varphi_j$ over $B_j$, and construct a global function
$\varphi=\sum\theta_j \varphi_j$ by means of a partition of unity $\{\theta_j\}$
subordinate to~$B_j$. Now, we observe that $\varphi-\varphi_k$ is smooth
on $B_k$  because all differences $\varphi_j-\varphi_k$ are smooth in
the intersections $B_j\cap B_k$ and we can write 
$\varphi-\varphi_k=\sum\theta_j(\varphi_j-\varphi_k)$.
Therefore $\alpha:=T-dd^c\varphi$ is smooth.\QED
\endproof

Thanks to Lemma~\ref{build-potential}, the problem of approximating a quasi-positive closed $(1,1)$-current is reduced to approximating a quasi-psh function. In this way, we get

\claim{\label{current-bergman-approx}{1.6}.\ Theorem}
Let $T=\alpha+dd^c\varphi$ be a quasi-positive  closed
$(1,1)$-current on a compact Hermitian
manifold $(X,\omega)$ such that $T\ge\gamma$ for some
continuous $(1,1)$-form~$\gamma$. Then there exists a sequence of
quasi-positive currents $T_m=\alpha+dd^c\varphi_m$ whose local
potentials have the
same singularities as $1/2m$ times a logarithm of a sum of squares
of holomorphic functions and
a decreasing sequence $\varepsilon_m>0$ converging to~$0$, such that
\vskip3pt
\item{\rm (a)} $T_m$ converges weakly to $T$,
\vskip3pt
\item{\rm (b)} $\displaystyle\nu(T,x)-\frac{n}{m}\le\nu(T_m,x)\le
\nu(T,x)$~ for every $x\in X\,;$
\vskip3pt
\item{\rm (c)} $T_m\ge\gamma-\varepsilon_m\omega$.
\vskip3pt
We say that our currents $T_m$ are approximations of $T$ with analytic
singularities $($possessing logarithmic poles\/$)$. Moreover,
for any multiplicative
subsequence $m_k$, one can arrange that $T_{m_k}=\alpha+dd^c\varphi_{m_k}$
where $(\varphi_{m_k})$ is a non-increasing sequence of potentials.
\endclaim

\proof{Proof}  We just briefly sketch the idea -- essentially a partition of unity argument~-- and refer to \cite{Dem92} for the details. Let us write $T=\alpha+dd^c\varphi$ with $\alpha$ smooth, according to Lemma~\ref{build-potential}. After replacing $T$ with $T-\alpha$ and $\gamma$ with
$\gamma-\alpha$, we can assume without loss of
generality that $\{T\}=0$, i.e. that $T=dd^c\varphi$
with a quasi-psh function~$\varphi$ on~$X$ such that 
$dd^c\varphi\ge\gamma$. Now, for $\varepsilon>0$ small, we select a finite covering $(B_j)_{1\le j\le N(\varepsilon)}$ of $X$ by coordinate balls on which there exists an $\varepsilon$-approximation of
$\gamma$ as 
$$
\sum_{1\le\ell\le n}\lambda_{j,\ell}\ii dz^{j}_\ell\wedge d\overline z^{j}_\ell\le
\gamma_{|B_j}\le 
\sum_{1\le\ell\le n} (\lambda_{j,\ell}+\varepsilon)
\ii dz^{j}_\ell\wedge d\overline z^{j}_\ell
$$
in terms of holomorphic coordinates $(z^j_\ell)_{1\le\ell\le n}$ on $B_j$ (for this we just diagonalize $\gamma(a_j)$ at the center $a_j$ of $B_j$, and take the radius of $B_j$ small enough). By construction 
$\psi_{j,\varepsilon}(z)=\varphi(z)-\sum_{1\le\ell\le n}\lambda_{j,\ell}|z^j_\ell|^2$ is psh on $B_\ell$, and we can thus obtain approximations $\psi_{j,\varepsilon,m}$ of $\psi_j$ by the Bergman kernel process applied on each ball~$B_j$. The idea is to define a global approximation of $\varphi$ by putting
$$
\varphi_{\varepsilon,m}(x)=\frac{1}{m}\log\bigg(\sum_{1\le j\le N(\varepsilon)}\theta_{j,\varepsilon}(x)\,\exp\Big(m\Big(\psi_{j,\varepsilon,m}(x)+\sum_{1\le\ell\le n}(\lambda_{j,\ell}-\varepsilon)|z^j_\ell|^2\Big)\Big)\bigg)
$$
where $(\theta_{j,\varepsilon})_{1\le j\le N(\varepsilon)}$ is a partition of unity subordinate to the $B_j$'s. If we take $\varepsilon=\varepsilon_m$ and $\varphi_m=\varphi_{\varepsilon_m,m}$ where $\varepsilon_m$ decays very slowly, then it is not hard to check that $T_m=dd^c\varphi_m$ satisfies the required estimates; it is essentially enough to observe that the derivatives of $\theta_{j,\varepsilon}$ are ``killed'' by the factor $\frac{1}{m}$ when $m\gg\frac{1}{\varepsilon}$.\QED
\endproof

We need a variant of Th.~\ref{current-bergman-approx} providing more
``equisingularity'' in the sense that the multiplier ideal sheaves are
preserved. If one adds the requirement to obtain a non-increasing
sequence of approximations of the potential, one can do this only at
the expense of accepting ``transcendental'' singularities, which can
no longer be guaranteed to be logarithmic poles.

\claim{\label{equisingular-approx}{1.7}.\ Theorem}
Let $T=\alpha+dd^c\varphi$ be a closed
$(1,1)$-current on a compact Hermitian manifold $(X,\omega)$, where
$\alpha$ is a smooth closed $(1,1)$-form and $\varphi$ a quasi-psh
function. Let $\gamma$ be a continuous real $(1,1)$-form such that
$T\ge\gamma$. Then one can write
\hbox{$\varphi=\lim_{m\to+\infty}\widetilde\varphi_m$} where
\vskip3pt
\item{\rm (a)} $\widetilde\varphi_m$ is smooth in the complement $X\ssm Z_m$
of an analytic set $Z_m\subset X\,;$
\vskip3pt
\item{\rm (b)} $\{\widetilde\varphi_m\}$ is a non-increasing sequence, and $Z_m\subset Z_{m+1}$ for all~$m\,;$
\vskip3pt
\item{\rm (c)} $\int_X(e^{-\varphi}-e^{-\widetilde\varphi_m})dV_\omega$
is finite for every $m$ and converges to $0$ as $m\to+\infty\,;$
\vskip3pt
\item{\rm (d)} $($``equisingularity''$)$ $\cI(\widetilde\varphi_m)=\cI(\varphi)$ for all $m$
$\,;$
\vskip3pt
\item{\rm (e)} $T_m=\alpha+dd^c\widetilde\varphi_m$ satisfies $T_m\ge
\gamma-\varepsilon_m\omega$, where $\lim_{m\to+\infty}\varepsilon_m=0$.
\endclaim

\proof{Proof}  (A substantial simplication of the original proof in 
\cite{DPS01}.) As in the previous proof, we may assume that $\alpha=0$
and $T=dd^c\varphi$, and after subtracting a constant to
$\varphi$ we can also achieve that $\varphi\le -1$ everywhere on~$X$.
For every germ $f\in\cO_{X,x}$, (c) implies 
$\int_U|f|^2(e^{-\varphi}-e^{-\widetilde\varphi_m})dV_\omega<+\infty$ on
some neighborhood $U$ of~$x$, hence the integrals
$\int_U|f|^2e^{-\varphi}dV_\omega$ and $\int_U|f|^2e^{-\widetilde\varphi_m}dV_\omega$ are simultaneously convergent or divergent, and (d) follows trivially.
We define
$$
\widetilde\varphi_m(x)=\sup_{k\ge m}(1+2^{-k})\varphi_{p_k}
$$
where $(p_k)$ is a multiplicative sequence that grows 
fast enough, with $\varphi_{p_{k+1}}\le\varphi_{p_k}\le 0$ for all $k$.
Clearly $\widetilde\varphi_m$ is a non-increasing sequence, and
$$\lim_{m\to+\infty}\widetilde\varphi_m(x)=\lim_{k\to+\infty}
\varphi_{p_k}(x)=\varphi(x)
$$
at every point~$x\in X$. If $Z_m$ is the set of poles of $\varphi_{p_m}$,
it is easy to see that 
$$\widetilde\varphi_m(x)=\lim_{\ell\to+\infty}
\sup_{k\in[m,\ell]}(1+2^{-k})\varphi_{p_k}
$$
converges uniformly on every compact subset of $X\ssm Z_m$, since any new
term $(1+2^{-\ell})\varphi_{p_\ell}$ may contribute to the sup only in case 
$$\varphi_{p_\ell}\ge \frac{1+2^{-p_m}}{1+2^{-p_\ell}}\varphi_{p_m}\quad(\ge 2
\varphi_{p_m}),$$
and the difference of that term with respect to the previous term 
$(1+2^{-(\ell-1)})\varphi_{p_{\ell-1}}\ge (1+2^{-(\ell-1)})\varphi_{p_\ell}$
is less than 
$2^{-\ell}|\varphi_{p_\ell}|\le 2^{1-\ell}|\varphi_{p_m}|$. Therefore $\widetilde
\varphi_m$ is continuous on~$X\ssm Z_m$, and getting it to be smooth is only a matter of applying Richberg's approximation technique (\cite{Ric68}, \cite{Dem12}). The only serious thing to prove is property~(c). To achieve this, we observe that 
$\{\varphi<\widetilde\varphi_m\}$ is contained in the union $\bigcup_{k\ge m}
\{\varphi<(1+2^{-k})\varphi_{p_k}\}$, therefore
$$
\int_X\big(e^{-\varphi}-e^{-\widetilde\varphi_m}\big)dV_\omega\le
\sum_{k=m}^{+\infty}
\int_X{\bf 1}_{\varphi<(1+2^{-k})\varphi_{p_k}} e^{-\varphi}dV_\omega
\leqno(\label{int-decomp}{1.8})$$
and
$$\leqalignno{
\int_X{\bf 1}_{\varphi<(1+2^{-k})\varphi_{p_k}} e^{-\varphi}dV_\omega
&=\int_X{\bf 1}_{\varphi<(1+2^{-k})\varphi_{p_k}} \exp\big(2^k\varphi-(2^k+1)\varphi\big)dV_\omega\cr
&\le\int_X{\bf 1}_{\varphi<(1+2^{-k})\varphi_{p_k}} \exp\big((2^k+1)(\varphi_{p_k}-\varphi)\big)dV_\omega\cr
&\le\int_X{\bf 1}_{\varphi<(1+2^{-k})\varphi_{p_k}} \exp\big(2p_k(\varphi_{p_k}-\varphi)\big)dV_\omega&(\label{right-hand-int-bergman}{1.9})\cr}
$$
if we take $p_k>2^{k-1}$ (notice that $\varphi_{p_k}-\varphi\ge 0$). Now, by 
Lemma~\ref{bergman-convergence-estimate} below, our integral (\ref{right-hand-int-bergman}) is finite. By Lebesgue's monotone convergence theorem, we have for $k$ fixed
$$
\lim_{p\to+\infty}\int_X{\bf 1}_{\varphi<(1+2^{-k})\varphi_p} e^{-\varphi}dV_\omega=0
$$
as a decreasing limit, and we can take $p_k$ so large that
$\int_{\varphi<(1+2^{-k})\varphi_{p_k}} e^{-\varphi}dV_\omega\le 2^{-k}$.
This ensures that property (c) holds true by (\ref{int-decomp}).\QED
\endproof

\claim{\label{bergman-convergence-estimate}{1.10}.\ Lemma}
On a compact complex manifold, for any quasi-psh potential $\varphi$, the Bergman kernel procedure leads to quasi-psh potentials $\varphi_m$ with analytic singularities such that
$$\int_Xe^{2m(\varphi_m-\varphi)}dV_\omega<+\infty.$$
\endclaim

\proof{Proof} 
By definition of $\varphi_m$ in Th.~\ref{psh-bergman-approx}, $\exp(2m(\varphi_{m}))$ is (up to the irrelevant partition of unity procedure) equal to the Bergman kernel $\sum_{\ell\in\bN}|g_{m,\ell}|^2$. By local uniform convergence and the Noetherian property, it has the same local vanishing behavior as a finite sum $\sum_{\ell\le N(m)}|g_{m,\ell}|^2$ with $N(m)$ sufficiently large. Since all terms $g_{m,\ell}$ have $L^2$ norm equal to~$1$ with respect to the weight $e^{-2m\varphi}$, our contention follows.\QED
\endproof

\claim{1.11.\ Remark} {\rm A very slight variation of the proof would yield the improved condition
\vskip2pt\noindent
(c')~~~$\displaystyle\forall\lambda\in\bR_+,~~\int_X(e^{-\lambda\varphi}-e^{-\lambda\widetilde\varphi_m})dV_\omega\le 2^{-m}$ for $m\ge m_0(\lambda)$,
\vskip2pt\noindent
and thus an equality $\cI(\lambda\widetilde\varphi_m)=\cI(\lambda\varphi)$ for $m\ge m_0(\lambda)$. We just need to replace estimate (\ref{int-decomp}) by
$$
\int_X\big(e^{-m\varphi}-e^{-m\widetilde\varphi_m}\big)dV_\omega\le
\sum_{k=m}^{+\infty}
\int_X{\bf 1}_{\varphi<(1+2^{-k})\varphi_{p_k}} e^{-k\varphi}dV_\omega
$$
and take $p_k$ so large that $2p_k\ge k(2^k+1)$ and
$\int_{\varphi<(1+2^{-k})\varphi_{p_k}} e^{-k\varphi}dV_\omega\le 2^{-k-1}$.\QED}
\endclaim

We also quote the following very simple consequence of Lemma~\ref{bergman-convergence-estimate}, which will be needed a bit later. Since $\varphi_m$ is less singular than~$\varphi$, we have of course an inclusion $\cI(\lambda\varphi)\subset\cI(\lambda\varphi_m)$ for all~$\lambda\in\bR_+$. Conversely~:

\claim{\label{multiplier-ideals-comparison}{1.12}.\ Corollary}
For every pair of positive real numbers $\lambda'>\lambda>0$, we have an inclusion of multiplier ideals
$$\cI(\lambda'\varphi_m)\subset\cI(\lambda\varphi)\qquad
\hbox{as soon as}~~m\ge \Big\lceil
\frac{1}{2}\frac{\lambda\lambda'}{\lambda'-\lambda}\Big\rceil.
$$
\endclaim

\proof{Proof}  If $f\in\cO_{X,x}$ and $U$ is a sufficiently small neighborhood of~$x$, the H\"older ine\-quality for conjugate exponents $p,q>1$ yields
$$
\int_U|f|^2e^{-\lambda\varphi}dV_\omega\le
\Big(\int_U|f|^2e^{-\lambda'\varphi_m}dV_\omega\Big)^{1/p}
\Big(\int_U|f|^2e^{\frac{q}{p}\lambda'\varphi_m-q\lambda\varphi}dV_\omega\Big)^{1/q}.
$$
Therefore, if $f\in\cI(\lambda'\varphi_m)_x$, we infer that
$f\in\cI(\lambda\varphi)_x$ as soon as the integral
$\int_Xe^{\frac{q}{p}\lambda'\varphi_m-q\lambda\varphi}dV_\omega$ is convergent. If we
select  $p\in{}]1,\lambda'/\lambda]$, this is implied by the condition
$\int_Xe^{q\lambda(\varphi_m-\varphi)}dV_\omega<+\infty$. If we further take $q\lambda=2m_0$ to
be an even integer so large that 
$$p=\frac{q}{q-1}=\frac{2m_0/\lambda}{2m_0/\lambda-1}\le\frac{\lambda'}{\lambda},\qquad\hbox{e.g.~~}m_0=m_0(\lambda,\lambda')=\Big\lceil
\frac{1}{2}\frac{\lambda\lambda'}{\lambda'-\lambda}\Big\rceil,$$
then we indeed have $\int_Xe^{2m_0(\varphi_m-\varphi)}dV_\omega\le
\int_Xe^{2m(\varphi_m-\varphi)}dV_\omega<+\infty$ for 
$m\ge m_0(\lambda,\lambda')$, thanks to Lemma~\ref{bergman-convergence-estimate}.\QED
\endproof

\claim{1.13.\ Remark} {\rm Without the monotonicity requirement (b) for the sequence $(\widetilde\varphi_m)$ in Theorem~\ref{equisingular-approx}, the strong openness conjecture proved in the next section would directly provide an equisingular sequence, simply by taking 
$$\widehat\varphi_m=\left(1+\frac{1}{m}\right)\varphi_m$$
where $\varphi_m$ is the Bergman approximation sequence. In fact all $\widehat\varphi_m$ have analytic singularities and Cor.~\ref{multiplier-ideals-comparison}
applied with $\lambda=1$ and $\lambda'=1+1/m$ shows that \hbox{$\cI(\widehat\varphi_m)\subset \cI(\varphi)$}. Since $\widehat\varphi_m\ge (1+\frac{1}{m})\varphi$, the equality $\cI(\widehat\varphi_m)=\cI(\varphi)$ holds for $m$ large by strong openness, and properties \ref{equisingular-approx}~(a), (c), (d), (e) can be seen to hold. However, the sequence $(\widehat\varphi_m)$ is not monotone.}
\endclaim

\section{\label{section-strong-openness}{2}}{Semi-continuity of psh singularities and proof of the strong openness conjecture}

In this section, we present a proof of the strong openness conjecture for multiplier ideal sheaves. Let $\Omega$ be a domain in $\bC^n$, $f\in\cO(\Omega)$ a holomorphic function, and $\varphi\in \PSH(\Omega)$ a psh function on~$\Omega$. For every holomorphic function $f$ on $\Omega$, we introduce the {\em weighted log canonical threshold} of $\varphi$ with weight $f$ at $z_0$
$$c_{f,z_0}(\varphi) = \sup\big\{ c>0:\ |f|^2 e^{-2c\,\varphi} \hbox{ is } L^1 \hbox{ on a neighborhood of }z_0\big\}\in{}]0,+\infty].$$
The special case $f=1$ yields the usual log canonical threshold $c_{z_0}(\varphi)$ that was defined in the introduction. The openness conjectures can be stated as follows.\vskip5pt

\claim{\label{openness-conj}{2.1}.\ Conjectures}
\vskip2pt
\item{\rm(a)} {\rm(openness conjecture, \cite{DK01})}\\ The set~~$\{c>0:\ e^{-2c\,\varphi}~\hbox{\em is }~L^1~\hbox{\em on a neighborhood of }z_0\big\}$ equals the open interval~~$]0,c_{z_0}(\varphi)[.$
\vskip3pt
\item{\rm(b)} {\rm(strong openness conjecture, \cite{Dem00})}\\ The set~~$\{c>0:\ |f|^2e^{-2c\,\varphi}~\hbox{\em is }~L^1~\hbox{\em on a neighborhood of }z_0\big\}$ equals the open interval~~$]0,c_{f,z_0}(\varphi)[.$\vskip0pt
\endclaim

The openness conjecture \ref{openness-conj}~(a) was first established by Favre and Jonsson (\cite{FJ05} in dimension $2$ (see also \cite{JM12}, \cite{JM14}), and 8 years later by Berndtsson \cite{Bern13} in arbitrary dimension. The strong form \ref{openness-conj}~(b), which is equivalent to Cor.~\ref{cor-upper-regularization}, was settled very recently by Guan and Zhou \cite{GZ13}. Their proof uses a sophisticated version of the $L^2$-extension theorem of Ohsawa and Takegoshi in combination with the curve selection lemma. They have also obtained related semi-continuity statements in  \cite{GZ14a} and ``effective versions'' in \cite{GZ14b}. A~simplified proof along the same lines has been given by Lempert in \cite{Lem14}.

Here, we follow Pham Hoang Hiep's approach \cite{Hiep14}, which is more straightforward and avoids the curve selection lemma. It is based on the original version \cite{OT87} of the $L^2$-extension theorem, applied to members of a standard basis for a multiplier ideal sheaf of holomorphic functions associated with a plurisubharmonic function $\varphi$. In this way, by means of a simple induction on dimension, one can obtain the strong openness conjecture, and give simultaneously an effective version of the semicontinuity theorem for weighted log canonical thresholds. The main results are contained in the following theorem.

\claim{\label{thm-hiep}{2.2}.\ Theorem} {\rm(\cite{Hiep14})} Let $f$ be a holomorphic function on an open set $\Omega$ in $\bC^n$ and let $\varphi$ be a psh function on $\Omega$.
\vskip2pt
\item{\rm(i)} {\rm(``Semicontinuity theorem'')} Assume that $\intn_{\Omega'}e^{-2c\,\varphi}dV_{2n}<+\infty$ on some open subset $\Omega'\subset\Omega$ and let $z_0\in\Omega'$. Then there exists $\delta=\delta(c,\varphi,\Omega',z_0)>0$ such that for every $\psi\in\PSH(\Omega')$, $\Vert \psi-\varphi\Vert_{L^1(\Omega')}\le \delta$ implies \hbox{$c_{z_0}(\psi)>c$}. Moreover, as $\psi$ converges to $\varphi$ in~$L^1(\Omega')$, the function $e^{-2c\,\psi}$ converges to $e^{-2c\,\varphi}$ in~$L^1$ on every relatively compact open subset $\Omega''\compact\Omega'$.
\vskip3pt
\item{\rm (ii)} {\rm(``Strong effective openness'')} Assume that $\intn_{\Omega'}|f|^2e^{-2c\,\varphi}dV_{2n}<+\infty$ on some open subset $\Omega'\subset\Omega$. When $\psi\in\PSH(\Omega')$ converges to $\varphi$ in $L^1(\Omega')$ with $\psi\le\varphi$, the function $|f|^2e^{-2c\,\psi}$ converges to $|f|^2e^{-2c\,\varphi}$ in $L^1$ norm on every relatively compact open subset $\Omega''\compact\Omega'$.
\vskip0pt
\endclaim

\claim{\label{cor-strong-openness}{2.3}.\ Corollary} {\rm(``Strong openness'')}. {\it For any plurisubharmonic function $\varphi$ on a neighborhood of a point $z_0\in\bC^n$, the set $\{ c>0:\ |f|^2 e^{-2c\,\varphi} \hbox{ is } L^1 \hbox{ on a neighborhood of }z_0\}$ is an open interval $(0,c_{f,z_0}(\varphi))$.}
\endclaim

\claim{\label{cor-convergence-below}{2.4}.\ Corollary} {\rm(``Convergence from below'')}. {\it If $\psi\le\varphi$ converges to $\varphi$ in a neighborhood of $z_0\in\bC^n$, then $c_{f,z_0}(\psi)\leq c_{f,z_0}(\varphi)$ converges to $c_{f,z_0}(\varphi)$.}
\endclaim

In fact, after subtracting a large constant to~$\varphi$, we can assume $\varphi\leq 0$ in both corollaries. Then Cor.~\ref{cor-strong-openness} is a consequence of assertion (ii) of the main theorem when we take $\Omega'$ small enough and $\psi=(1+\delta)\varphi$ with~$\delta\searrow 0$. In Cor.~\ref{cor-convergence-below}, we have by definition $c_{f,z_0}(\psi)\leq c_{f,z_0}(\varphi)$ for $\psi\leq\varphi$, but again (ii) shows that $c_{f,z_0}(\psi)$ becomes${}\geq c$ for any given value $c\in(0,c_{f,z_0}(\varphi))$, whenever $\Vert\psi-\varphi\Vert_{L^1(\Omega')}$ is sufficiently small.

\claim{2.5.\ Remark} One cannot remove condition $\psi\leq\varphi$ in assertion (ii) of the main theorem. Indeed, choose $f(z)=z_1$, $\varphi(z)=\log |z_1|$ and $\varphi_j(z) = \log |z_1+\frac {z_2} j|$, for $j\geq 1$.One has $\varphi_j\to\varphi$ in $L_{\rm loc}^1(\bC^n)$, however $c_{f,0}(\varphi_j)=1<c_{f,0}(\varphi)=2$ for all $j\geq 1$. On the other hand, condition (i) of Theorem~\ref{thm-hiep} does not require any given inequality between $\varphi$ and $\psi$. Modulo Berndtsson's solution of the openness conjecture, (i) follows from the effective semicontinuity result of \cite{DK01}, but (like Guan and Zhou) Hiep's technique will reprove both by a direct and easier method.
\endclaim

\noindent
{\bf 2.6.\ A few preliminaries.} According to standard techniques in the theory of Gr\"obner bases, one equips the ring $\cO_{\bC^n,0}$ of germs of holomorphic functions at $0$ with the homogeneous lexicographic order of monomials $z^\alpha = z_1^{\alpha_1} \ldots z_n^{\alpha_n}$, that is, $z_1^{ \alpha_1 } \ldots z_n^{ \alpha_n }<z_1^{ \beta_1 } \ldots z_n^{ \beta_n }$ if and only if $| \alpha |=\alpha_1+\ldots+\alpha_n<| \beta |=\beta_1+\ldots+\beta_n$ or $| \alpha |=| \beta |$ and $\alpha_i < \beta_i$ for the first index $i$ with $\alpha_i \not = \beta_i$. For each $f(z)=a_{\alpha^1} z^{\alpha^1}+a_{\alpha^2} z^{\alpha^2}+\ldots\;$ with $a_{\alpha^j}\not = 0$, $j\geq 1$ and $z^{ \alpha^1 } < z^{ \alpha^2 } < \ldots\;$, we define the {\it initial coefficient}, {\it initial monomial} and  {\it initial term} of $f$ to be respectively $\IC(f)=a_{\alpha^1}$, $\IM(f)=z^{\alpha^1}$, $\IT(f)=a_{\alpha^1} z^{\alpha^1}$, and the support of $f$ to be $\SUPP(f)=\{z^{ \alpha^1 }, z^{ \alpha^2 }, \ldots\}$. For any ideal $\cI$~of~$\cO_{\bC^n,0}$, we define $\IM(\cI)$ to be the ideal generated by $\{ \IM(f) \}_{\{ f\in\cI\} }$. First, we recall the division theorem of Hironaka and the concept of standard basis of an ideal.

\claim{2.7.\ Theorem} {\rm (Division theorem of Hironaka, \cite{Gal79}, \cite{Bay82}, \cite{BM87}, \cite{BM89}, \cite{Eis95})} Let $f, g_1,\ldots,g_k\in \cO_{\bC^n,0}$. Then there
exist $h_1,\ldots,h_k,s\in \cO_{\bC^n,0}$ such that
$$f=h_1g_1+\ldots+h_kg_k+s,$$
and $\SUPP(s)\cap \langle\IM(g_1),\ldots,\IM(g_k)\rangle =\emptyset$, where $\langle\IM(g_1),\ldots,\IM(g_k)\rangle$ denotes the ideal generated by the family~$(\IM(g_1),\ldots,\IM(g_k))$.
\endclaim

\noindent
{\bf 2.8.\ Standard basis of an ideal.} Let $\cI$ be an ideal of $\cO_{\bC^n,0}$ and let $g_1,\ldots,g_k\in\cI$ be such that $\IM(\cI)=\langle\IM(g_1),\ldots,\IM(g_k)\rangle$. Take $f\in\cI$. By the division theorem of Hironaka, there exist $h_1,\ldots,h_k,s\in \cO_{\bC^n,0}$ such that
$$f=h_1g_1+\ldots+h_kg_k+s,$$
and $\SUPP(s)\cap\IM(\cI) =\emptyset$. On the other hand, since $s=f-h_1g_1+\ldots+h_kg_k\in\cI$, we have \hbox{$\IM(s)\in\IM(\cI)$}. Therefore $s=0$ and the $g_j$'s are generators of~$\cI$. By permuting the $g_j$'s and performing ad hoc subtractions, we can always arrange that
$\IM(g_1)<\IM(g_2)<\ldots<\IM(g_k)$, and we then say
that $(g_1,\ldots,g_k)$ is a standard basis of $\cI$.
\vskip5pt

\noindent
Th.~\ref{thm-hiep} will be proved by induction on dimension $n$. All statements are trivial for \hbox{$n=0$}. Assume that the theorem holds for dimension $n-1$. Thanks to the $L^2$-extension theorem of Ohsawa and Takegoshi (\cite{OT87}), one obtains the following key lemma.

\claim{\label{lemma-hiep}{2.9}.\ Lemma} Let $\varphi\leq 0$ be a plurisubharmonic function and $f$ be a holomorphic function on the polydisc $\Delta_R^n$ of center $0$ and $($poly$)$radius $R>0$ in $\bC^n$, such that for some $c>0$
$$\intn_{\Delta_R^n} |f(z)|^2 e^{-2c\,\varphi(z)} d V_{2n}(z)<+\infty.$$
Let $\psi_j\leq 0$, $j\geq 1$, be a nequence of plurisubharmonic functions on
$\Delta_R^n$ with $\psi_j\to\varphi$ in $L^1_{\rm loc}(\Delta_R^n)$, and assume
that either $f=1$ identically or $\psi_j\leq\varphi$ for all~$j\geq 1$.
Then for every $r<R$ and $\varepsilon\in(0,\frac{1}{2}r]$, there exist a value $w_n\in\Delta_\varepsilon \smallsetminus\{0\}$, an index $j_0$, a constant $\tilde c>c$  and a sequence of holomorphic functions $F_j$ on $\Delta_{r}^n$, $j\geq j_0$, such that $\IM(F_j)\leq \IM(f)$, $F_j(z) = f(z) + (z_n-w_n)\sum a_{j,\alpha} z^{\alpha}$ with $| w_n| | a_{j,\alpha} |\leq r^{-|\alpha|}\varepsilon$ for~all~$\alpha\in\bN^n$, and
$$\intn_{ \Delta_{r}^n } |F_j(z)|^2 e^{-2\tilde c\,\psi_j(z)} d V_{2n}(z)\leq\frac{\varepsilon^2}{|w_n|^2}<+\infty,~~~\forall j\geq j_0.$$
Moreover, one can choose $w_n$ in a set of positive measure in the
punctured disc $\Delta_\varepsilon\smallsetminus\{0\}$ $($the index $j_0=j_0(w_n)$ and the constant $\tilde c=\tilde c(w_n)$ may then possibly depend on $w_n)$.
\endclaim

\proof{Proof} By Fubini's theorem we have
$$\intn_{\Delta_R}\bigg[\intn_{\Delta_R^{n-1}} |f(z',z_n)|^2 e^{-2c\,\varphi(z',z_n)} d V_{2n-2}(z')\bigg] dV_2(z_n)<+\infty.$$
Since the integral extended to a small disc $z_n\in \Delta_\eta$ tends to $0$ as $\eta\to 0$, it will become smaller than any preassigned value, say $\varepsilon_0^2>0$, for $\eta\leq\eta_0$ small enough. Therefore we can choose a set of positive measure of values $w_n\in\Delta_\eta\smallsetminus \{ 0\}$ such that
$$\intn_{\Delta_R^{n-1}} |f(z',w_n)|^2e^{-2c\,\varphi(z',w_n)} d V_{2n-2}(z')\leq\frac { \varepsilon_0^2}{\pi\eta^2}<\frac {\varepsilon_0^2}{|w_n|^2}.$$
Since the main theorem is assumed to hold for $n-1$, for any $\rho<R$ there exist $j_0=j_0(w_n)$ and $\tilde c=\tilde c(w_n)>c$ such that
$$\intn_{\Delta_{\rho}^{n-1}} |f (z',w_{n}) |^2 e^{-2\tilde c\,\psi_j(z',w_n)} d V_{2n-2}(z')<\frac {\varepsilon_0^2}{|w_n|^2},~~~\forall j\geq j_0.$$
(For this, one applies part (i) in case $f=1$, and part (ii) in case $\psi_j\leq\varphi$, using the fact that $\psi=\frac{\tilde c}{c}\,\psi_j$ converges to $\varphi$ as $\tilde c\to c$ and $j\to+\infty$). Now, by the $L^2$-extension theorem of Ohsawa and Takegoshi (see \cite{OT87}), there exists a holomorphic function $F_j$ on $\Delta_{\rho}^{n-1}\times\Delta_R$ such that $F_j(z',w_n)=f(z',w_n)$ for all $z'\in\Delta_{\rho}^{n-1}$,
and
$$
\eqalign{
\intn_{\Delta_{\rho}^{n-1}\times \Delta_R} |F_j(z)|^2 e^{-2\tilde c\,\psi_j(z)} d V_{2n}(z)&\leq C_nR^2\intn_{\Delta_{\rho}^{n-1}} |f(z',w_n)|^2 e^{-2\tilde c\,\psi_j(z',w_n)} d V_{2n-2}(z')\cr
&\leq\frac { C_n R^2 \varepsilon_0^2}{ |w_n|^2 },\cr}
$$
where $C_n$ is a constant which only depends on $n$ (the constant is universal for $R=1$ and is rescaled by $R^2$ otherwise). By the mean value inequality for the plurisubharmonic function $|F_j|^2$, we get
$$
\eqalign{
|F_j(z)|^2&\leq \frac 1 {\pi^n(\rho -|z_1|)^2\ldots(\rho -|z_n|)^2}\intn_{ \Delta_{\rho -|z_1|} (z_1)\times\ldots\times\Delta_{\rho -|z_n|} (z_n) } |F_j|^2 dV_{2n}\cr
&\leq \frac {C_n R^2 \varepsilon_0^2}{\pi^n(\rho -|z_1|)^2\ldots(\rho -|z_n|)^2|w_n|^2},\cr}
$$
where $\Delta_\rho(z)$ is the disc of center $z$ and radius $\rho$. Hence, for any $r<R$, by taking $\rho=\frac{1}{2}(r+R)$ we infer
$$
\Vert F_j\Vert _{ L^\infty (\Delta_{r}^n) }\leq \frac {2^nC_n^{\frac 1 2}R\varepsilon_0}{\pi^{\frac n 2} (R-r)^n |w_n| }.
\leqno(\label{inequality-Fj}{2.10})
$$
Since $F_j(z',w_n)-f(z',w_n)=0$, $\forall z'\in\Delta_r^{n-1}$, we can write $F_j(z)=f(z)+(z_n-w_n)g_j(z)$ for some function $g_j(z)=\sum_{\alpha\in\bN^n} a_{j,\alpha}z^{\alpha}$ on $\Delta_r^{n-1}\times\Delta_R$. By (\ref{inequality-Fj}), we get
$$
\eqalign{
\Vert g_j\Vert _{ \Delta_{r}^n } = \Vert g_j\Vert _{ \Delta_{r}^{n-1}\times\partial\Delta_{r} }
&\leq \frac {1} {r-|w_n|} \Big(\Vert F_j\Vert _{ L^\infty (\Delta_{r}^n) }+\Vert f\Vert _{ L^\infty (\Delta_{r}^n) }\Big)\cr
&\leq \frac {1} {r-|w_n|} \Big(\frac {2^nC_n^{\frac 1 2}R \varepsilon_0} { \pi^{\frac n 2} (R-r)^n |w_n| }+\Vert f\Vert _{ L^\infty (\Delta_{r}^n) }\Big).\cr}
$$
Thanks to the Cauchy integral formula, we find
$$| a_{j,\alpha} |\leq \frac { \Vert g_j\Vert _{ \Delta_{r}^n } } { r^{ |\alpha | } }\leq \frac {1} { (r-|w_n|) r^{ |\alpha| } } \Big(\frac {2^nC_n^{\frac 1 2}R \varepsilon_0} { \pi^{\frac n 2} (R-r)^n |w_n| }+\Vert f\Vert _{ L^\infty (\Delta_{r}^n) }\Big).$$
We take in any case $\eta\leq\varepsilon_0\leq\varepsilon\leq\frac{1}{2}r$. As $|w_n|<\eta\leq\frac{1}{2}r$, this implies
$$| w_n | | a_{j,\alpha} |\,r^{|\alpha |}\leq \frac {2}{r} \Big(\frac {2^nC_n^{\frac 1 2}R \varepsilon_0} { \pi^{\frac n 2} (R-r)^n }+\Vert f\Vert _{ L^\infty (\Delta_{r}^n) } |w_n| \Big)\leq C'\varepsilon_0,$$
for some constant $C'$ depending only on $n,\,r,\,R$ and~$f$. This yields the estimates of Lemma~\ref{lemma-hiep} for $\varepsilon_0:=C''\varepsilon$ with $C''$ sufficiently small.
Finally, we prove that $\IM(F_j)\leq \IM(f)$. Indeed, if $\IM(g_j)\geq\IM(f)$, since $|w_n\Vert a_{j,\alpha}|\leq r^{-|\alpha| }\varepsilon$, we can choose $\varepsilon$ small enough such that $\IM(F_j) = \IM(f)$ and $\Big|\frac {\textstyle\IC(F_j)}{\textstyle\IC(f) }\Big|\in (\frac 1 2, 2)$. Otherwise, if $\IM(g_j)<\IM(f)$, we have $\IM(F_j)=\IM(g_j)<\IM(f)$.\QED
\endproof

\noindent
{\bf 2.11.\ Proof of Theorem~\ref{thm-hiep}.} By well-known properties of (pluri)potential theory, the $L^1$ convergence of $\psi$ to $\varphi$ implies that $\psi\to\varphi$ almost everywhere, and the assumptions guarantee that $\varphi$ and $\psi$ are uniformly bounded on every relatively compact subset of~$\Omega'$. In particular, after shrinking $\Omega'$ and subtracting constants, we can assume that $\varphi\le 0$ on~$\Omega$. Also, since the $L^1$ topology is metrizable, it is enough to work with a sequence $(\psi_j)_{j\geq 1}$ converging to $\varphi$ in $L^1(\Omega')$. Again, we can assume that $\psi_j\leq 0$ and that $\psi_j\to\varphi$ almost everywhere on $\Omega'$. By a trivial compactness argument, it is enough to show (i) and (ii) for some neighborhood $\Omega''$ of a given point $z_0\in\Omega'$. We assume here $z_0=0$ for simplicity of notation, and fix a polydisc $\Delta_R^n$ of center $0$ with $R$ so small that $\Delta_R^n\subset\Omega'$. Then \hbox{$\psi_j(\bu,z_n)\to\varphi (\bu,z_n)$} in the topology of $L^1(\Delta_R^{n-1})$ for almost every $z_n\in\Delta_R$.
\vskip3pt

\noindent
{\bf 2.11~(i).\ Proof of statement {\rm(i)} in Theorem~\ref{thm-hiep}}. We have here $\intn_{\Delta_R^n}e^{-2c\,\varphi}dV_{2n}<+\infty$ for $R>0$ small enough. By Lemma~\ref{lemma-hiep} with $f=1$, for every $r<R$ and $\varepsilon>0$, there exist $w_n\in\Delta_{\varepsilon}\smallsetminus\{0\}$, an index $j_0$, a number $\tilde c>c$ and a sequence of holomorphic functions $F_j$ on $\Delta_r^n$, $j\geq j_0$, such that $F_j(z)=1+(z_n-w_n)\sum a_{j,\alpha}z^\alpha$, $|w_n||a_{j,\alpha}|\,r^{-|\alpha|}\le\varepsilon$ and
$$\intn_{\Delta_r^n}|F_j(z)|^2e^{-2\tilde c\,\psi_j(z)}dV_{2n}(z)\le\frac{\varepsilon^2}{|w_n|^2},~~~\forall j\geq j_0.$$
For $\varepsilon\leq\frac{1}{2}$, we conclude that $|F_j(0)|=|1-w_na_{j,0}|\geq\frac{1}{2}$ hence $c_0(\psi_j)\geq\tilde c>c$ and the first part of (i) is proved. In fact, after fixing such $\varepsilon$ and $w_n$, we even obtain the existence of a neighborhood $\Omega''$ of $0$ on which $|F_j|\geq\frac{1}{4}$, and thus get a uniform bound $\int_{\Omega''}e^{-2\tilde c\,\psi_j(z)}dV_{2n}(z)\le M<+\infty$. The second assertion of (i) then follows from the estimate
$$
\eqalign{
\intn_{\Omega''}\big|e^{-2c\,\psi_j(z)}-e^{-2c\,\varphi(z)}\big|
dV_{2n}(z)&\leq
\intn_{\Omega''\cap\{|\psi_j|\leq A\}}\big|e^{-2c\,\psi_j(z)}-e^{-2c\,\varphi(z)}\big|dV_{2n}(z)\cr
&\qquad{}+\intn_{\Omega''\cap\{\psi_j<-A\}}e^{-2c\,\varphi(z)}dV_{2n}(z)\cr
&\qquad{}+e^{-2(\tilde c-c)A}\intn_{\Omega''\cap\{\psi_j<-A\}}e^{-2\tilde c\,\psi_j(z)}dV_{2n}(z).\cr}
$$
In fact the last two terms converge to $0$ as $A\to+\infty$, and, for $A$ fixed, the first integral in the right hand side converges to $0$ by Lebesgue's bounded
convergence theorem, since $\psi_j\to\varphi$ almost everywhere on
$\Omega''$.
\vskip3pt

\noindent
{\bf 2.11~(ii).\ Proof of statement {\rm (ii)} in Theorem~\ref{thm-hiep}.} Take $f_1,\ldots,f_k\in\cO_{\bC^n, 0}$ such that $(f_1,\ldots,f_k)$ is a standard basis of $\cI(c\,\varphi)_0$ with $\IM(f_1)<\ldots<\IM(f_k)$, and $\Delta_R^n$ a polydisc so small that
$$\intn_{\Delta_R^n}|f_l(z)|^2e^{-2c\,\varphi(z)}dV_{2n}(z)<+\infty,~~~l=1,\ldots,k.$$
Since the germ of $f$ at $0$ belongs to the ideal $(f_1,\ldots,f_k)$, we can essentially argue with the $f_l$'s instead of~$f$. By Lemma~\ref{lemma-hiep}, for every $r<R$ and $\varepsilon_l>0$, there exist $w_{n,l}\in\Delta_{\varepsilon_l}\smallsetminus\{0\}$, an index $j_0=j_0(w_{n,l})$, a number $\tilde c=\tilde c(w_{n,l})>c$ and a sequence of holomorphic functions $F_{j,l}$ on $\Delta_r^n$, $j\geq j_0$, such that $F_{j,l}(z)=1+(z_n-w_{n,l})\sum a_{j,l,\alpha}z^\alpha$, $|w_{n,l}||a_{j,l,\alpha}|\,r^{-|\alpha|}\le\varepsilon_l$ and
$$
\intn_{\Delta_r^n}|F_{j,l}(z)|^2e^{-2\tilde c\,\psi_j(z)}dV_{2n}(z)\le\frac{\varepsilon_l^2}{|w_{n,l}|^2},~~~\forall l=1,\ldots,k,~~\forall j\geq j_0.
\leqno(\label{integral-bound}{2.12)}
$$
Since $\psi_j\leq\varphi$ and $\tilde c>c$, we get $F_{j,l}\in\cI(\tilde c\,\psi_j)_0\subset \cI(c\,\varphi)_0$. The next step of the proof consists in modifying $(F_{j,l})_{1\leq l\leq k}$ in order to obtain a standard basis of $\cI(c\,\varphi)_0$. For this, we proceed by selecting successively $\varepsilon_1\gg\varepsilon_2\gg\ldots\gg\varepsilon_k$ (and suitable $w_{n,l}\in\Delta_{\varepsilon_l}\smallsetminus\{0\}$). We have $\IM(F_{j,1}),\ldots,\IM(F_{j,k})\in\IM(\cI(c\,\varphi)_0$, in particular $\IM(F_{j,1})$ is divisible by $\IM(f_l)$ for some $l=1,\ldots,k$. Since $\IM(F_{j,1})\leq\IM(f_1)<\ldots<\IM(f_k)$, we must have $\IM(F_{j,1})=\IM(f_1)$ and thus $\IM(g_{j,1})\geq\IM(f_1)$. As $|w_{n,1}||a_{j,1,\alpha}|\le\varepsilon_1$, we will have $\Big|\frac {\textstyle\IC(F_{j,1}) } {\textstyle\IC(f_1) }\Big|\in (\frac 1 2, 2)$ for $\varepsilon_1$ small enough. Now, possibly after changing $\varepsilon_2$ to a smaller value, we show that there exists a polynomial $P_{j,2,1}$ such that the degree and coefficients of $P_{j,2,1}$ are uniformly bounded, with $\IM(F_{j,2}-P_{j,2,1} F_{j,1}) = \IM(f_2)$ and $\frac {\textstyle|\IC (F_{j,2}-P_{j,2,1} F_{j,1})|} {\textstyle|\IC(f_2)|}\in (\frac 1 2, 2)$. We consider two cases:
\vskip5pt
\noindent
{\bf Case 1}\/: If $\IM(g_{j,2})\geq\IM(f_2)$, since $|w_{n,2}||a_{j,2,\alpha}|\leq r^{-|\alpha|}\varepsilon_2$, we can choose $\varepsilon_2$ so small that $\IM(F_{j,2}) = \IM(f_2)$ and $\frac {\textstyle|\IC(F_{j,2})|}{\textstyle|\IC(f_2)|}\in (\frac 1 2, 2)$. We then take $P_{j,2,1}=0$.
\vskip5pt
\noindent
{\bf Case 2}\/: If $\IM(g_{j,2})<\IM(f_2)$, we have $\IM(g_{j,2}) = \IM(F_{j,2})\in \IM(\cI(c\,\varphi)_0)$. Hence $\IM(g_{j,2})$ is divisible by $\IM(f_l)$ for some $l=1,\ldots,k$. However, since $\IM(g_{j,2})<\IM(f_2)<\ldots<\IM(f_k)$, the only possibility is that $\IM(g_{j,2})$ be divisible by $\IM(f_1)$. Take $b\in\bC$ and $\beta,\gamma\in\bN^n$ such that $\IT(g_{j,2}):=a_{j,2,\gamma}z^{\gamma}=bz^{\beta}\,\IT(F_{j,1})$. We have $z^{\beta}\leq z^{\gamma}=\IM(g_{j,2})<\IM(f_2)$ and
$$|w_{n,2}||b|=|w_{n,2}|\frac {|\IC(g_{j,2})|} {|\IC(F_{j,1})|}\leq
\frac {2|w_{n,2}||a_{j,2,\gamma}|} {|\IC(f_1)|}
\leq \frac{2r^{-|\gamma|}\varepsilon_2}{|\IC(f_1)|}$$
can be taken arbitrarily small.
Set $\tilde g_{j,2}(z)=g_{j,2}(z)- bz^{\beta} F_{j,1}(z)=\sum \tilde a_{j,2,\alpha} z^{\alpha}$ and
$$\tilde F_{j,2}(z)=f_2(z)+(z_n-w_{n,2})\tilde g_{j,2}(z)=F_{j,2}(z)-b(z_n-w_{n,2})z^{\beta} F_{j,1}(z).$$
We have $\IM(\tilde g_{j,2}) > \IM(g_{j,2})$. Since $|w_{n,2}||b|=O(\varepsilon_2)$ and $|w_{n,2}|| a_{j,2,\alpha}|=O(\varepsilon_2)$, we get $|w_{n,2}||\tilde a_{j,2,\alpha}|=O(\varepsilon_2)$ as well. Now, we consider two further cases. If $\IM(\tilde g_{j,2}) \geq \IM(f_2)$, we can again change $\varepsilon_2$ for a smaller value so that $\IM(\tilde F_{j,2}) = \IM(f_2)$ and $\frac {\textstyle|\IC (\tilde F_{j,2})|}{\textstyle|\IC (f_2)|}\in (\frac 1 2, 2)$. Otherwise, if $\IM(\tilde g_{j,2}) < \IM(f_2)$, we have $\IM(F_{j,2}) = \IM(g_{j,2}) < \IM( \tilde F_{j,2} ) = \IM( \tilde g_{j,2} ) < \IM(f_2)$. Notice that \hbox{$\{ z^{ \gamma }:\ z^{ \gamma }<\IM(f_2) \}$} is a finite set. By using similar arguments a finite number of times, we find $\varepsilon_2$ so small\break that $\IM(F_{j,2}-P_{j,2,1}F_{j,1}) = \IM(f_2)$ and $\frac {\textstyle|\IC(F_{j,2}-P_{j,2,1} F_{j,1})|}{\textstyle|\IC(f_2)|}\in (\frac 1 2, 2)$ for some polynomial $P_{j,2,1}$.
\noindent
Repeating the same arguments for $F_{j,3},\ldots,F_{j,k}$, we select inductively $\varepsilon_l$, $l=1,\ldots,k$, and construct linear combinations 
$$F'_{j,l}=F_{j,l}-\sum_{1\leq m\leq l-1}P_{j,l,m}F'_{j,m}$$
with polynomials $P_{j,l,m}$, $1\leq m<l\le k$, possessing uniformly bounded coefficients and degrees, such that $\IM(F'_{j,l})=\IM(f_l)$ and
$\frac {\textstyle|\IC (F'_{j,l})|} {\textstyle|\IC(f_l)|}\in (\frac 1 2, 2)$ for all $l=1,\ldots,k$ and $j\ge j_0$. This implies that
$(F'_{j,1},\ldots,F'_{j,k})$
is also a standard basis of $\cI(c\,\varphi)_0$. By Theorem 1.2.2 in \cite{Gal79}, we can find $\rho$,~\hbox{$K>0$} so small that there exist holomorphic functions $h_{j,1},\ldots,h_{j,k}$ on $\Delta_{\rho}^n$ with $\rho<r$, such that
$$f=h_{j,1}F'_{j,1}+h_{j,2}F'_{j,2}+\ldots+h_{j,k}F'_{j,k}~~~
\hbox{on $\Delta_{\rho}^n$}$$
and $\Vert h_{j,l}\Vert _{L^\infty (\Delta_{\rho}^n)}\leq K\Vert f\Vert _{L^\infty (\Delta_{r}^n)}$, for all $l=1,\ldots,k$ ($\rho$ and $K$ only depend on $f_1,\ldots,f_k$). By (\ref{integral-bound}), this implies a uniform bound
$$\intn_{ \Delta_\rho^n } |f(z)|^2e^{-2\tilde c\,\psi_j(z)} d V_{2n}(z)\le M<+\infty$$
for some $\tilde c>c$ and all $j\geq j_0$. Take $\Omega''=\Delta_\rho^n$.
We obtain the $L^1$ convergence of $|f|^2e^{-2c\,\psi_j}$ to $|f|^2e^{-2c\,\varphi}$ almost exactly as we argued for the second assertion of part~(i), by using the estimate
$$
\eqalign{
\intn_{\Omega''}|f|^2\big|e^{-2c\,\psi_j(z)}-e^{-2c\,\varphi(z)}\big|dV_{2n}(z)
&\leq
\intn_{\Omega''\cap\{|\psi_j|\leq A\}}|f|^2\big|e^{-2c\,\psi_j(z)}-e^{-2c\,\varphi(z)}\big|dV_{2n}(z)\cr
&\qquad{}+\intn_{\Omega''\cap\{\psi_j<-A\}}|f|^2e^{-2c\,\varphi(z)}dV_{2n}(z)\cr
&\qquad{}+e^{-2(\tilde c-c)A}\intn_{\Omega''\cap\{\psi_j<-A\}}|f|^2e^{-2\tilde c\,\psi_j(z)}dV_{2n}(z).\cr}
$$

\section{\label{section-hard-lefschetz}{3}}{Hard Lefschetz theorem for pseudoeffective line bundles}

\subsection{3.1}{A variant of the Bochner formula}
We first recall a variation of the Bochner formula that is required in the proof of the Hard Lefschetz Theorem with values in a positively curved (and therefore non flat) line bundle $(L,h)$. Here the base manifold is a K\"ahler (non necessarily compact) manifold $(Y,\omega)$.  We denote by
$|~~|=|~~|_{\omega,h}$ the pointwise Hermitian norm on
\hbox{$\Lambda^{p,q}T^*_Y\otimes L$} associated with $\omega$ and $h$,
and by $\Vert~~\Vert=\Vert~~\Vert_{\omega,h}$ the global $L^2$ norm
$$
\Vert u\Vert^2 = \int_Y |u|^2 dV_\omega \qquad
\hbox{where}\quad  dV_\omega=\frac{\omega^n}{n!}
$$
We consider the $\dbar$ operator acting on $(p,q)$-forms with values in
$L$, its adjoint $\dbar_h^*$ with respect to $h$ and the complex
Laplace-Beltrami operator $\Delta''_h=\dbar\dbar_h^*+\dbar_h^*\dbar$.
Let $v$ be a smooth $(n-q,0)$-form with compact support
in~$Y$. Then $u=\omega^q\wedge v$ satisfies
$$\Vert\dbar u\Vert^2+\Vert\dbar^*_h u\Vert^2=
\Vert\dbar v\Vert^2+\int_Y\sum_{I,J}\Big(\sum_{j\in J}\lambda_j\Big)|u_{IJ}|^2
\leqno(\label{bochner-variant}{3.1.1})$$
where $\lambda_1\le\cdots\le\lambda_n$ are the curvature
eigenvalues of $\Theta_{L,h}$ expressed in an orthonormal frame
$(\partial/\partial z_1,\ldots,\partial/\partial z_n)$ (at some fixed
point $x_0\in Y$), in such a way that
$$
\omega_{x_0}=\ii\sum_{1\le j\le n}dz_j\wedge d\overline z_j,\qquad
(\Theta_{L,h})_{x_0}=dd^c\varphi_{x_0}=
\ii\sum_{1\le j\le n} \lambda_jdz_j\wedge d\overline z_j.
$$
Formula (\ref{bochner-variant}) follows from the more or less
straightforward identity
$$
(\dbar^*_\varphi\,\dbar+\dbar\,\dbar^*_\varphi)(v\wedge\omega^q)-
(\dbar^*_\varphi\,\dbar v)\wedge\omega^q=
q\,\ii\ddbar\varphi\wedge\omega^{q-1}\wedge v,
$$
by taking the inner product with $u=\omega^q\wedge v$ and integrating
by parts in the left hand side (we leave the easy details to the reader). 
Our formula is thus established when 
$v$ is smooth and compactly supported. In general, we have:

\claim{\label{bochner-variant-complete}{3.1.2}.\ Proposition}
Let $(Y,\omega)$ be a {\rm complete} K\"ahler
manifold and $(L,h)$ a smooth Hermitian line bundle such that the
curvature possesses a uniform lower bound $\Theta_{L,h}\ge -C\omega$.
For every measurable $(n-q,0)$-form $v$ with $L^2$ coefficients and
values in $L$ such that $u=\omega^q\wedge v$ has differentials
$\dbar u$, $\dbar^* u$ also in $L^2$, we have
$$
\Vert\dbar u\Vert^2+\Vert\dbar^*_h u\Vert^2=
\Vert\dbar v\Vert^2+\int_Y\sum_{I,J}\Big(\sum_{j\in J}\lambda_j\Big)|u_{IJ}|^2
$$
$($here, all differentials are computed in the sense of distributions$)$.
\endclaim

\proof{Proof}  Since $(Y,\omega)$ is assumed to be complete,
there exists a sequence of smooth forms $v_\nu$ with compact support in $Y$
(obtained by truncating $v$ and taking the convolution with a regularizing
kernel) such that $v_\nu\to v$ in $L^2$ and such that
$u_\nu=\omega^q\wedge v_\nu$ satisfies $u_\nu\to u$, $\dbar u_\nu\to\dbar u$,
$\dbar^* u_\nu\to\dbar^* u$ in~$L^2$. By the curvature assumption,
the final integral in the right hand side of (\ref{bochner-variant})
must be under control
(i.e.\ the integrand becomes nonnegative if we add a term $C\Vert u\Vert^2$
on both sides, $C\gg 0$). We thus get the equality by passing to the limit
and using Lebesgue's monotone convergence theorem.\QED
\endproof

\subsection{3.2}{Proof of Theorem \ref{hard-lefschetz}}

Here $X$ denotes a compact K\"ahler manifold equipped with a K\"ahler
metric~$\omega$, and $(L,h)$ is a pseudoeffective line bundle on~$X$. 
To fix the ideas, we first indicate the~proof in the
much simpler case when $(L,h)$ has a smooth metric $h$ (so that
$\cI(h)=\cO_X$), and then treat the general case.

\subsubsection{3.2.1}{Special Case: $(L,h)$ is Hermitian semipositive (with a smooth metric).}

Let $\{\beta\}\in H^q(X,\Omega^n_X\otimes L)$ be an arbitrary cohomology class.
By standard $L^2$ Hodge theory, $\{\beta\}$ can be represented by a smooth
harmonic $(0,q)$-form $\beta$ with values in $\Omega^n_X\otimes L$. We can
also view $\beta$ as a $(n,q)$-form with values in $L$. The pointwise
Lefschetz isomorphism produces a unique $(n-q,0)$-form $\alpha$ such
that $\beta=\omega^q\wedge\alpha$. Proposition~\ref{bochner-variant-complete}
then yields
$$
\Vert\dbar\alpha\Vert^2+\int_Y\sum_{I,J}\Big(\sum_{j\in J}\lambda_j\Big)
|\alpha_{IJ}|^2=\Vert\dbar\beta\Vert^2+\Vert\dbar^*_h \beta\Vert^2=0,
$$
and the curvature eigenvalues $\lambda_j$ are nonnegative by our assumption.
Hence $\dbar\alpha=0$ and $\{\alpha\}\in H^0(X,\Omega^{n-q}_X\otimes L)$
is mapped to $\{\beta\}$ by $\Phi^q_{\omega,h}=\omega^q\wedge\bu~$.

\subsubsection{3.2.2}{General Case.}

There are several difficulties. The first difficulty is that the
metric $h$ is no longer smooth and we cannot directly represent
cohomology classes by harmonic forms. We circumvent this problem by
smoothing the metric on an (analytic) Zariski open subset and by
avoiding the remaining poles on the complement. However, some careful
estimates have to be made in order to take the error terms into account.

Fix $\varepsilon=\varepsilon_\nu$ and let $h_\varepsilon=h_{\varepsilon_\nu}$
be an approximation of~$h$, such that $h_\varepsilon$ is smooth on
$X\ssm Z_\varepsilon$ ($Z_\varepsilon$ being an analytic subset of $X$),
$\Theta_{L,h_\varepsilon}\ge -\varepsilon\omega$,
$h_\varepsilon\le h$ and $\cI(h_\varepsilon)=\cI(h)$. This is possible
by Th.~\ref{equisingular-approx}. Now, we can find a family
$$
\omega_{\varepsilon,\delta}=\omega+\delta(\ii\ddbar\psi_\varepsilon+\omega),
\qquad \delta>0
$$
of {\it complete K\"ahler} metrics on $X\ssm Z_\varepsilon$, where
$\psi_\varepsilon$ is a quasi-psh function on $X$ with
$\psi_\varepsilon=-\infty$ on $Z_\varepsilon$,
$\psi_\varepsilon$ smooth on $X\ssm Z_\varepsilon$ and $\ii\ddbar
\psi_\varepsilon+\omega\ge 0$ (see e.g.\ \cite{Dem82}, Th\'eor\`eme~1.5).
By construction, $\omega_{\varepsilon,\delta}\ge\omega$ and
$\lim_{\delta\to 0}\omega_{\varepsilon,\delta}=\omega$.
We look at the $L^2$ Dolbeault complex $K^\bu_{\varepsilon,\delta}$
of $(n,\bu)$-forms on $X\ssm Z_\varepsilon$, where the $L^2$ norms are
induced by $\omega_{\varepsilon,\delta}$ on differential forms and by
$h_\varepsilon$ on elements in~$L$. Specifically
$$
K^q_{\varepsilon,\delta}=\Big\{u{:}X\ssm Z_\varepsilon{\to}\Lambda^{n,q}
T^*_X\otimes L;\int_{X\ssm Z_\varepsilon}\kern-15pt
(|u|^2_{\Lambda^{n,q}\omega_{\varepsilon,\delta}
\otimes h_\varepsilon}+|\dbar u|^2_{\Lambda^{n,q+1}\omega_{\varepsilon,\delta}
\otimes h_\varepsilon})dV_{\omega_{\varepsilon,\delta}}<\infty\Big\}.
$$
Let $\cK^q_{\varepsilon,\delta}$ be the corresponding sheaf of germs
of locally $L^2$ sections on $X$ (the local $L^2$ condition
should hold on $X$, not only on $X\ssm Z_\varepsilon\,$!). Then,
for all $\varepsilon>0$ and $\delta\ge 0$,
$(\cK^q_{\varepsilon,\delta},\dbar)$ is a resolution of the sheaf
$\Omega^n_X\otimes L\otimes\cI(h_\varepsilon)=
\Omega^n_X\otimes L\otimes\cI(h)$. This is
because $L^2$ estimates hold locally on small Stein open sets, and the
$L^2$ condition on $X\ssm Z_\varepsilon$ forces holomorphic sections
to extend across~$Z_\varepsilon$ (\cite{Dem82}, Lemma 6.9).

Let $\{\beta\}\in H^q(X,\Omega^n_X\otimes L\otimes\cI(h))$ be a
cohomology class represented by a smooth form with values in
$\Omega^n_X\otimes L\otimes\cI(h)$ (one can use a \v Cech cocycle
and convert it to an element in the $\cC^\infty$ Dolbeault complex by
means of a partition of unity, thanks to the usual De Rham-Weil
isomorphism, see also the final proof in Section~\ref{section-cao-proof}
for more details). Then
$$
\Vert\beta\Vert_{\varepsilon,\delta}^2\le \Vert\beta\Vert^2=
\int_X|\beta|^2_{\Lambda^{n,q}\omega\otimes h}dV_\omega<+\infty.
$$
The reason is that $|\beta|^2_{\Lambda^{n,q}\omega\otimes h}dV_\omega$
decreases as $\omega$ increases. This is just an easy calculation,
shown by comparing two metrics $\omega$, $\omega'$ which are
expressed in diagonal form in suitable coordinates; the norm
$|\beta|^2_{\Lambda^{n,q}\omega\otimes h}$ turns out to decrease
faster than the volume $dV_\omega$ increases; see e.g.\ \cite{Dem82}, Lemma 3.2;
a special case is $q=0$, then $|\beta|^2_{\Lambda^{n,q}\omega\otimes h}
dV_\omega=i^{n^2}\beta\wedge\overline\beta$ with the identification
$L\otimes\overline L\simeq\bC$ given by the metric $h$, hence the
integrand is even independent of $\omega$ in that case.

By the proof of the De Rham-Weil isomorphism, the map
$\alpha\mapsto\{\alpha\}$ from the cocycle space
$Z^q(\cK^\bu_{\varepsilon,\delta})$ equipped with its $L^2$ topology,
into $H^q(X,\Omega^n_X\otimes L\otimes\cI(h))$ equipped with its
finite vector space topology, is continuous.  Also, Banach's open
mapping theorem implies that the coboundary space
$B^q(\cK^\bu_{\varepsilon,\delta})$ is closed in
$Z^q(\cK^\bu_{\varepsilon,\delta})$. This is true for all $\delta\ge
0$ (the limit case $\delta=0$ yields the strongest $L^2$ topology in
bidegree $(n,q)$). Now, $\beta$ is a $\dbar$-closed form in the Hilbert space
defined by $\omega_{\varepsilon,\delta}$ on $X\ssm Z_\varepsilon$, so there
is a $\omega_{\varepsilon,\delta}$-harmonic form
$u_{\varepsilon,\delta}$ in the same cohomology class
as $\beta$, such that
$$\Vert u_{\varepsilon,\delta}\Vert_{\varepsilon,\delta}\le
\Vert\beta\Vert_{\varepsilon,\delta}.
\leqno(\label{harmonic-rep}{3.2.3})$$
Let $v_{\varepsilon,\delta}$ be the unique $(n-q,0)$-form such that
$u_{\varepsilon,\delta}=v_{\varepsilon,\delta}\wedge
\omega_{\varepsilon,\delta}^q$ ($v_{\varepsilon,\delta}$ exists by the
pointwise Lefschetz isomorphism). Then
$$
\Vert v_{\varepsilon,\delta}\Vert_{\varepsilon,\delta}=
\Vert u_{\varepsilon,\delta}\Vert_{\varepsilon,\delta}\le
\Vert\beta\Vert_{\varepsilon,\delta}\le\Vert\beta\Vert.
$$
As $\sum_{j\in J}\lambda_j\ge -q\varepsilon$ by the assumption on
$\Theta_{L,h_\varepsilon}$, the Bochner formula yields
$$
\Vert\dbar v_{\varepsilon,\delta}\Vert_{\varepsilon,\delta}^2\le
q\varepsilon\Vert u_{\varepsilon,\delta}\Vert_{\varepsilon,\delta}^2
\le q\varepsilon\Vert\beta\Vert^2.
$$
These uniform bounds imply that there are subsequences $u_{\varepsilon,
\delta_\nu}$ and $v_{\varepsilon,\delta_\nu}$ with \hbox{$\delta_\nu\to 0$},
possessing weak-$L^2$ limits $u_\varepsilon=
\lim_{\nu\to+\infty}u_{\varepsilon,\delta_\nu}$
and $v_\varepsilon=\lim_{\nu\to+\infty}v_{\varepsilon,\delta_\nu}$.
The limit $v_\varepsilon=
\lim_{\nu\to+\infty}v_{\varepsilon,\delta_\nu}$ is with respect to
$L^2(\omega)=L^2(\omega_{\varepsilon,0})$. To check this, notice that
in bidegree $(n-q,0)$, the space $L^2(\omega)$ has the weakest topology
of all spaces~$L^2(\omega_{\varepsilon,\delta})$; indeed, an easy calculation
made in (\cite{Dem82}, Lemma 3.2) yields
$$
|f|^2_{\Lambda^{n-q,0}\omega\otimes h}dV_\omega\le
|f|^2_{\Lambda^{n-q,0}\omega_{\varepsilon,\delta}\otimes h}
dV_{\omega_{\varepsilon,\delta}}\qquad
\hbox{if $f$ is of type $(n-q,0)$}.
$$
On the other hand, the limit
$u_\varepsilon=\lim_{\nu\to+\infty}u_{\varepsilon,\delta_\nu}$
takes place in all spaces $L^2(\omega_{\varepsilon,\delta})$, $\delta>0$,
since the topology gets stronger and stronger as $\delta\downarrow 0$
[$\,$possibly not in $L^2(\omega)$, though, because in bidegree $(n,q)$
the topology of $L^2(\omega)$ might be strictly stronger than that
of all spaces $L^2(\omega_{\varepsilon,\delta})\,$].
The above estimates yield
$$\eqalign{
&\Vert v_\varepsilon\Vert^2_{\varepsilon,0}=
\int_X|v_\varepsilon|^2_{\Lambda^{n-q,0}\omega\otimes h_\varepsilon}
dV_\omega\le\Vert\beta\Vert^2,\cr
&\Vert \dbar v_\varepsilon\Vert^2_{\varepsilon,0}\le q\varepsilon
\Vert\beta\Vert^2_{\varepsilon,0},\cr
\noalign{\vskip3pt}\cr
&u_\varepsilon=\omega^q\wedge v_\varepsilon\equiv\beta
\qquad\hbox{in}~~H^q(X,\Omega^n_X\otimes L\otimes\cI(h_\varepsilon)).\cr}
$$
Again, by arguing in a given Hilbert space $L^2(h_{\varepsilon_0})$,
we find $L^2$ convergent subsequences $u_\varepsilon\to u$,
$v_\varepsilon\to v$ as $\varepsilon\to 0$, and in this way get
$\dbar v=0$ and
$$\eqalign{
&\Vert v\Vert^2\le \Vert \beta\Vert^2,\cr
&u=\omega^q\wedge v\equiv \beta
\qquad\hbox{in}~~H^q(X,\Omega^n_X\otimes L\otimes\cI(h)).\cr}
$$
Theorem~\ref{hard-lefschetz} is proved. Notice that the equisingularity property
$\cI(h_\varepsilon)=\cI(h)$ is crucial in the above proof,
otherwise we could not infer that $u\equiv \beta$ from the fact
that $u_\varepsilon\equiv \beta$. This is true only because
all cohomology classes $\{u_\varepsilon\}$ lie in the same fixed
cohomology group \hbox{$H^q(X,\Omega^n_X\otimes L\otimes\cI(h))$}, whose
topology is induced by the topology of $L^2(\omega)$ on $\dbar$-closed
forms (e.g.\ through the De Rham-Weil isomorphism).\QED

\claim{\label{nef-weird-example}{3.2.4}.\ Remark} {\rm
In (\ref{harmonic-rep}), the existence of a harmonic 
representative holds true only for~$\omega_{\varepsilon,\delta}$, $\delta>0$,
because we need to have a complete K\"ahler metric on $X\ssm Z_\varepsilon$.
The trick of employing $\omega_{\varepsilon,\delta}$ instead of a fixed
metric $\omega$, however, is not needed when $Z_\varepsilon$ is (or can
be taken to be) empty. This is the case if $(L,h)$ is such that
$\cI(h)=\cO_X$ and $L$ is nef. Indeed, by definition, $L$ is nef iff
there exists a sequence of smooth metrics $h_\nu$ such that
$\ii\Theta_{L,h_\nu}\ge-\varepsilon_\nu\omega$, so we can take the 
$\varphi_\nu$'s to be everywhere smooth in~Th.~\ref{equisingular-approx}.
However, multiplier ideal sheaves are needed in the surjectivity statement
even in case $L$ is nef, as it may happen that $\cI(h_{\min})\ne\cO_X$ even
then, and $h:=\lim h_\nu$ is anyway always more singular than $h_{\min}$.
Let us recall a standard example (see \cite{DPS94}, \cite{DPS01}).
Let $B$ be an elliptic curve and let $V$ be the rank $2$ vector bundle
over $B$ which is defined as the (unique) non split
extension
$$0 \to \cO_B \to V \to \cO_B \to 0.$$
In particular, the bundle $V$ is numerically flat, i.e.\ $c_1(V)=0$,
$c_2(V)=0$. We consider the ruled surface $X=\bP(V)$. On that surface
there is a unique section $C=\bP(\cO_B) \subset X$ with $C^2=0$ and
$$\cO_X (C)= \cO_{\bP (V)} (1)$$
is a nef line bundle. One can check that $L=\cO_{\bP (V)}(3)$ leads
to a {\it zero} Lefschetz map 
$$\omega\wedge\bu:~~H^0(X,\Omega^1_X\otimes L)\longrightarrow
H^1(X,K_X\otimes L)\simeq\bC,
$$
so this is a counterexample to Cor.~\ref{lefschetz-semipositive} in
the nef case. Incidentally, this also shows (in a somewhat sophisticated
way) that $\cO_{\bP (V)}(1)$ is nef but not semipositive, a fact that 
was first observed in \cite{DPS94}.}
\endclaim

\section{\label{section-num-dim}{4}}{Numerical dimension of currents}

A large part of this section borrows ideas from S.~Boucksom's \cite{Bou02}, \cite{Bou04} and Junyan Cao's \cite{JC14} PhD theses. We try however to give here a slightly more formal exposition. The main difference with S.~Boucksom's approach is that we insist on keeping track of singularities of currents and leaving them unchanged, instead of trying to minimize them in each cohomology class.

\subsection{\label{subsection-equising-approx}{4.1}}{Monotone asymptotically equisingular approximations}

Let $X$ be a compact complex $n$-dimensional manifold. 
We consider the closed convex cone of {\it pseudoeffective
classes}, namely the set $\cE(X)$ of cohomology classes
$\{\alpha\}\in H^{1,1}(X,\bR)$ containing a closed positive 
$(1,1)$-current $T=\alpha+dd^c\varphi$ (in the non K\"ahler case
one should use Bott-Chern cohomology groups here, but we will be mostly 
concerned with the K\"ahler case in the sequel). We also introduce
the set $\cS(X)$ of singularity equivalence classes of closed positive 
$(1,1)$-currents \hbox{$T=\alpha+dd^c\varphi$} (i.e., $\alpha$ 
being fixed, up to equivalence of singularities of the potentials~$\varphi$, 
cf.\ Def.~\ref{sing-ordering}). Clearly, there is a fibration
$$\pi:\cS(X)\to\cE(X),\qquad T\mapsto\{\alpha\}\in\cE(X)\subset H^{1,1}(X,\bR).
\leqno(\label{sing-equiv-fibration}{4.1.1})$$
We will denote by $\cS_\alpha(X)$ the fiber $\pi^{-1}(\{\alpha\})$ of $\cS(X)$
over a given cohomology class $\{\alpha\}\in\cE(X)$. Observe that the base $\cE(X)$ is a closed convex cone in a finite dimensional vector space, but in general the fiber $\cS_\alpha(X)$ must be viewed as a very complicated infinite dimensional space~: if we take e.g.\
$\{\alpha_1\}\in 
H^{1,1}(\bP^n,\bR)$ to be the unit class $c_1(\cO(1))$, then any current
$T=\frac{1}{d}[H]$ where $H_d$ is an irreducible hypersurface of degree $d$
defines a point in $\cS_{\alpha_1}(\bP^n)$, and these points are all 
distinct. The set $\cS(X)$ is nevertheless equipped in a natural way
with an addition law $\cS(X)\times\cS(X)\to \cS(X)$ that maps 
$\cS_\alpha(X)+\cS_\beta(X)$ into $\cS_{\alpha+\beta}(X)$, a scalar multiplication
$\bR_+\times\cS(X)\to\cS(X)$ that takes $\lambda\cdot\cS_\alpha(X)$
to the fiber~$\cS_{\lambda\alpha}(X)$. In this way, $\cS(X)$ should 
be viewed as some sort
of infinite dimensional convex cone. The fibers $\cS_\alpha(X)$ also possess
a partial ordering${}\preccurlyeq$ (cf.\ Def.~\ref{sing-ordering}) such
that $\forall j,~S_j\preccurlyeq T_j\Rightarrow\sum S_j\preccurlyeq\sum T_j$, and a fiberwise ``min'' operation
$$\leqalignno{
&\min:\cS_\alpha(X)\times \cS_\alpha(X)\longrightarrow\cS_\alpha(X),\cr
&(T_1,T_2)=(\alpha+dd^c\varphi_1,\alpha+dd^c\varphi_2)\longmapsto
T=\alpha+dd^c\max(\varphi_1,\varphi_2),
&(\label{min-operation}{4.1.2})\cr}
$$
with respect to which the addition is distributive, i.e.\
$$\min(T_1+S,T_2+S)=\min(T_1,T_2)+S.$$
Notice that when $T_1=\frac{1}{d}[H_1]$,
$T_2=\frac{1}{d}[H_2]$ are effective $\bQ$-divisors, all these operations $+$,~$\cdot\;$, $\min(\bu)$ and the ordering $\preccurlyeq$ coincide with the usual ones known for divisors. Following Junyan Cao \cite{JC14} (with slightly more 
restrictive requirements that do not produce much change in practice), we introduce

\claim{\label{asympt-equisingular-def}{4.1.3}.\ Definition}
Let $T=\alpha+dd^c\varphi$ be a closed positive $(1,1)$-current on~$X$,
where $\alpha$ is a smooth closed $(1,1)$-form and $\varphi$ is a quasi-psh 
function on~$X$. We say that the sequence of currents 
$T_k=\alpha+dd^c\psi_k$, $k\in\bN$,
is a ``monotone asymptotically equisingular approximation of~$T$ by currents with analytic 
singularities'' if the sequence of potentials $(\psi_k)$ satisfies the 
following properties\/$:$
\vskip2pt
\item{\rm(a)} {\rm(monotonicity)}
The sequence $(\psi_k)$ is non-increasing and
converges to~$\varphi$ at every point of~$X$.
\vskip2pt
\item{\rm(b)}
The functions $\psi_k$ have analytic singularities $($and $\psi_k\preccurlyeq\psi_{k+1}$ by {\rm(a)}$)$.
\vskip2pt
\item{\rm(c)} {\rm(lower bound of positivity)}
$$\alpha+dd^c\psi_k\geq  -\varepsilon_k\cdot \omega\qquad
\hbox{with}~~\lim_{k\to+\infty}\varepsilon_k=0$$
for any given smooth positive hermitian $(1,1)$-form $\omega$ on~$X$.
\vskip2pt
\item{\rm(d)} {\rm(asymptotic equisingularity)}
For every pair of positive numbers $\lambda'>\lambda>0$, there exists
an integer $k_{0}(\lambda,\lambda')\in \bN$ such that
$$\cI(\lambda'\psi_k) \subset \cI(\lambda\varphi)
\qquad\hbox{for}~~k\geq k_{0} (\lambda,\lambda').$$
\endclaim

\claim{4.1.4.\ Remark} {\rm  Without loss of generality, one can always assume that the
quasi-psh potentials $\varphi_k=c_k\log|g_k|^2+O(1)$ have rational coefficients
$c_k\in\bQ_+\,$; here again, $g_k$ is a tuple of locally defined holomorphic
functions. In fact, after subtracting constants, one can achieve
that $\varphi\le 0$ and $\psi_k\le 0$ for all~$k$. If the $c_k$ are 
arbitrary nonnegative real numbers, one can always replace $\psi_k$ by 
$\psi_k'=(1-\delta_k)\psi_k$ with a decreasing sequence $\delta_k\in{}]0,1[$ 
such that $\lim\delta_k=0$ and \hbox{$(1-\delta_k)c_k\in\bQ_+$}. 
Then (a), (b), (d) are still valid, and (c) holds with $\varepsilon'_k=(1-\delta_k)\varepsilon_k+C\delta_k$ and $C$ a constant such that $\alpha\ge -C\omega$.\QED}
\endclaim

\noindent
The fundamental observation is:

\claim{\label{theorem-bergman-approx-fund-observ}{4.1.5}.\ Theorem}  If $\psi_k:=\varphi_{m_k}$ is the sequence of potentials obtained by the Bergman kernel approximation of $\,T=\alpha+dd^c\varphi$ given in the proof of Theorem~\ref{current-bergman-approx} and $(m_k)$ is a multiplicative sequence, then the $\psi_k$ can be arranged to satisfy the positivity, monotonicity and asymptotic equisingularity properties of Definition~\ref{asympt-equisingular-def}. Moreover, if we start with
currents $T\preccurlyeq T'$ in the same cohomology class $\{\alpha\}$, we obtain corresponding approximations that satisfy $\psi_k\preccurlyeq\psi_k'$.
\endclaim

\proof{Proof}  By Cor.~\ref{multiplier-ideals-comparison}, the asymptotic equisingularity property~(d) in Def.~\ref{asympt-equisingular-def} is satisfied for $m_k\ge\lceil\frac{1}{2}\frac{\lambda\lambda'}{\lambda'-\lambda}\rceil$. The other properties are already known or obvious, especially the coefficients $\smash{c_k=\frac{1}{m_k}}$ are just inverses of integers in that case.\QED
\endproof

The following proposition provides a precise comparison of analytic singula\-rities of potentials when their multiplier ideal sheaves satisfy inclusion relations.

\claim{\label{sing-comparison}{4.1.6}.\ Proposition}
Let $\varphi$, $\psi$ be quasi-psh functions with analytic singularities, let $c>0$ be the constant such that $\varphi$ can be expressed as
$c\log\sum|g_j|^2+O(1)$ with holomorphic functions $g_j$, and 
let~\hbox{$\lambda\in\bR_+$}. Denoting $t_+:=\max(t,0)$, we have the implications\vskip2pt
\item{\rm(a)} $\displaystyle\forall f\in\cO_{X,x},~~
\int_{B_x\ni x} |f|^2e^{-\lambda\varphi}dV<+\infty~~~\Rightarrow~~~
\log|f|^2\succcurlyeq {\textstyle\frac{1}{c}}
\big(\lambda c-n\big)_+\varphi,$
\vskip2pt
\item{\rm(b)} $\displaystyle\cI(\psi)\subset\cI(\lambda\varphi)~~\Rightarrow~~
\int e^{\psi-\lambda\varphi}dV<+\infty~~\hbox{and}~~
\psi\succcurlyeq {\textstyle\frac{1}{c}}
\big(\lambda c-n\big)_+\varphi~~
\hbox{$($locally$)$}.$\vskip0pt
\endclaim

\proof{Proof} Since everything is local, we may assume that $\varphi$, $\psi$ are psh functions on a small ball $B\subset\bC^n$, and
$\varphi(z)=c\log|g|^2=c\log\sum_{1\le j\le N}|g_j(z)|^2$. 
\vskip3pt\noindent
(a) The convergence of the integral on a small ball $B_x$ of center $x$ implies
$$
\int_{B_x}|f|^2|g|^{-2\lambda c}dV\le \hbox{Const}
\int_{B_x}|f|^2e^{-\lambda\varphi}dV<+\infty
$$
By the openness of convergence exponents, one gets
$$
\int_{B_x}|f|^2|g|^{-2\lambda+\varepsilon}dV<+\infty
$$
for $\varepsilon>0$ small enough (this can be seen e.g.\ by using a log resolution of the ideal sheaf $(f,g_j)$). Now, if $\lambda c\ge n$,
Skoda's division theorem \cite{Sko72a} implies that each $f$ can be written
$f=\sum h_jg_j$ where $h_j$ satisfies a similar estimate
where the exponent of $|g|^{-2}$ is decreased by~$1$. An iteration of the Skoda division theorem for the $h_j$ yields $f\in(g_j)^k$ where 
$k=(\lfloor \lambda c\rfloor-(n-1))_+\ge (\lambda c-n)_+$. Hence
$$\log|f|^2\le k\log|g|^2 +C\le \frac{k}{c}\varphi+C'$$
and (a) is proved.
\vskip3pt\noindent
(b) If $(f_\ell)_{\ell\in\bN}$ is a Hilbert basis of the space of $L^2$ holomorphic functions $f$ with $\int_B|f|^2e^{-\psi}dV<+\infty$, the proof of Th.~\ref{psh-bergman-approx} yields $\psi\le C+\log\sum|f_\ell|^2$ (and locally the singularity is achieved by a finite sum of $f_\ell$'s by the Noetherian property). After possibly shrinking $B$, the relations $f_\ell\in\cI(\psi)\subset\cI(\lambda\varphi)$ imply 
$$\int_B|f_\ell|^2e^{-\lambda\varphi}dV<+\infty,$$
hence $\int e^{\psi-\lambda\varphi}dV<+\infty$ locally by taking the sum 
over~$\ell$. The inequality proved in (a) for each $f=f_\ell$ also yields
$$\psi\le \log\sum|f_\ell|^2+C\le \frac{1}{c}\big(\lambda c-n\big)_+\varphi+C',$$
and our singularity comparison relation follows.\QED
\endproof

\claim{4.1.7.\ Corollary}  If $T=\alpha+dd^c\varphi$ is a closed positive $(1,1)$-current
and $(\psi_k)$, $(\psi_k')$ are two monotone asymptotically equisingular approximations of $\varphi$ with analytic singularities, then for every $k$ and every $\varepsilon>0$, there exists $\ell$ such that $(1-\varepsilon)\psi_k\preccurlyeq\psi'_\ell$ $($and vice versa by exchanging the roles of 
$(\psi_k)$ and $(\psi_k'))$.
\endclaim

\proof{Proof}  Let $c>0$ be the constant occurring in the logarithmic poles 
of $\psi_k$ ($k$ being fixed). By condition (d) in Def.~\ref{asympt-equisingular-def}, for $\lambda'>\lambda\gg 1$ 
we have $\cI(\lambda'\psi'_\ell)\subset \cI(\lambda\varphi)\subset\cI(\lambda\psi_k)$
for $\ell\ge\ell_0(\lambda,\lambda')$ large enough. Proposition~\ref{sing-comparison} implies the singularity estimate $\psi'_\ell\succcurlyeq\frac{1}{c\lambda'}(c\lambda-n)_+\psi_k$, and the final constant in front of $\psi_k$ can be taken arbitrary close to~$1$.\QED
\endproof

Our next observation is that the $\min(\bu)$ procedure defined above for currents is well behaved in terms of asymptotic equisingular approximations.

\claim{\label{prop-bergman-max-relation}{4.1.8}.\ Proposition} Let $T=\alpha+dd^c\varphi$ and
$T'=\alpha+dd^c\varphi'$ be closed positive $(1,1)$-currents in the same cohomology class~$\{\alpha\}$. Let $(\psi_k)$ and $(\psi'_k)$ be respective monotone asymptotically equisingular approximations with analytic singularities and rational coefficients. Then $\max(\psi_k,\smash{\psi'_k})$ provides a monotone asymptotically equisingular approximation of $\min(T,T')=\alpha+dd^c\max(\varphi,\varphi')$ with analytic singularities and rational coefficients.
\endclaim

\proof{Proof}  If $\psi_k=c_k\log|g_k|^2+O(1)$ and $\psi'_k=c'_k\log|g'_k|^2+O(1)$,
we can write $c_k=p_k/q_k$, $c'_k=p'_k/q'_k$ and
$$
\max(\psi_k,\psi_k')=\frac{1}{q_kq'_k}\log\big(|g_k|^{2p_k}+|g'_k|^{2p'_k}\big)+O(1),
$$
hence $\max(\psi_k,\smash{\psi'_k})$ also has analytic singularities with rational coefficients (this would not be true with our definitions when the ratio
$c'_k/c_k$ is irrational, but of course we could just extend a little bit the definition of what we call analytic singularities, e.g.\ by allowing arbitrary positive real exponents, in order to avoid this extremely minor annoyance). It is well known that
$$\eqalign{
\alpha+dd^c\psi_k\ge-\varepsilon_k\omega,
&~~\alpha+dd^c\psi'_k\ge-\varepsilon'_k\omega\cr
&\Rightarrow~~\alpha+dd^c\max(\psi_k,\psi'_k)\ge-\max(\varepsilon_k,
\varepsilon'_k)\omega.\cr}
$$
Finally, if $\psi_{B,k}$ $($resp.\ $\psi'_{B,k}$ and $\widetilde\psi_{B,k}))$ comes from the Bergman approximation of $\varphi$ $($resp.\ of $\varphi'$ and
$\widetilde\varphi:=\max(\varphi,\varphi'))$, we have
$$\widetilde\varphi\ge\varphi~~\Rightarrow~~
\widetilde\psi_{B,k}\ge\psi_{B,k},\qquad
\widetilde\varphi\ge\varphi'~~\Rightarrow~~
\widetilde\psi_{B,k}\ge\psi'_{B,k}$$
hence $\widetilde\psi_{B,k}\ge\max(\psi_{B,k},\psi'_{B,k})$ and so
$\widetilde\psi_{B,k}\preccurlyeq\max(\psi_{B,k},\psi'_{B,k})$. However,
for every $\epsilon>0$, one has $(1-\varepsilon)\psi_{B_k}\preccurlyeq\psi_\ell$
and $(1-\varepsilon)\psi'_{B_k}\preccurlyeq\psi'_\ell$ for $\ell\ge\ell_0(k,\varepsilon)$ large, therefore $(1-\varepsilon)\widetilde\psi_{B,k}\preccurlyeq\max(\psi_\ell,\psi'_\ell)$. This shows that $\max(\psi_\ell,\psi'_\ell)$ has enough singularities (the ``opposite'' inequality $\max(\psi_\ell,\psi'_\ell)\ge\widetilde\varphi=\max(\varphi,\varphi')$, i.e.\ $\max(\psi_\ell,\psi'_\ell)\preccurlyeq\widetilde\varphi$, holds trivially).\QED
\endproof

Following Junyan Cao \cite{JC15}, we now investigate the additivity properties of the Bergman approximation procedure.

\claim{\label{theorem-bergman-additivity}{4.1.9}.\ Theorem} Let $T=\alpha+dd^c\varphi$ and $T'=\beta+dd^c\varphi'$ be closed $(1,1)$-currents in cohomology classes $\{\alpha\}$, $\{\beta\}\in\cE(X)$. Then for every multiplicative sequence $(m_k)$, the sum $\varphi_{m_k}+\varphi'_{m_k}$ of the Bergman approximations of $\varphi$, $\varphi'$ gives a monotone asymptotically equisingular approximation of $\varphi+\varphi'$ and $T+T'$. 
\endclaim

\proof{Proof}
Let $\widetilde\varphi_m$ be the Bergman kernel approximations of $\widetilde\varphi=\varphi+\varphi'$. By the sub\-addi\-tivity property of ideal sheaves $\cI(m\varphi+m\varphi') \subset \cI(m\varphi) \cI(m\varphi') $ (\cite{DEL00}, Th.~2.6), hence we have $\varphi_m + \varphi'_m  \preccurlyeq \widetilde\varphi_m $.
By Def.~\ref{asympt-equisingular-def}~(d), Th.~\ref{theorem-bergman-approx-fund-observ} and Cor.~\ref{corol-equiv-asympt-sing}, to~prove Th.~\ref{theorem-bergman-additivity}, it is sufficient to prove that for every $m\in\bN$ fixed, there exists a positive sequence $\lim\limits_{p\rightarrow +\infty}\varepsilon_p = 0$ such that 
$$
(1-\varepsilon_p) \widetilde\varphi_m \preccurlyeq \varphi_p + \varphi'_p \qquad\hbox{for  every }p \gg 1.
\leqno(\label{equivalproof}{4.1.10})
$$
For every $m\in\bN$ fixed, there exists a bimeromorphic map $\pi: \widetilde{X}\rightarrow X$, such that 
$$
\widetilde\varphi_m\circ \pi = \sum_i c_i \ln |s_i| + C^{\infty} \qquad\hbox{for some }c_i >0,
\leqno(\label{add2}{4.1.11})
$$
and the effective divisor $\sum_i \Div (s_i)$ is normal crossing. 
By the construction of $\widetilde\varphi_m$, we have $ \widetilde\varphi_m\preccurlyeq \varphi +\varphi'$. Therefore 
$$
\widetilde\varphi_m \circ  \pi \preccurlyeq (\varphi +\varphi')\circ \pi.
\leqno(\label{add1}{4.1.12})
$$
By Siu's decomposition formula for closed positive currents applied to $dd^c (\varphi\circ \pi)$, $dd^c (\varphi'\circ \pi)$ respectively, the divisorial parts add up to produce a divisor that is at least equal to the divisorial part in $dd^c(\widetilde\varphi_m\circ\pi)$, thus (\ref{add1}) and (\ref{add2}) imply the existence of numbers $a_i,\,b_i \geq 0$ satisfying
\vskip3pt\noindent
(i)~~ $a_i +b_i = c_i$ for every $i$,
\vskip3pt\noindent
(ii)~ $\displaystyle\sum\limits_i a_i \ln |s_i|  \preccurlyeq   \varphi\circ \pi$~ and~ $\displaystyle\sum\limits_i b_i \ln |s_i| \preccurlyeq   \varphi'\circ \pi$.
\vskip3pt\noindent
Let $p\in \bN$ be an integer,  $J$ be the Jacobian of $\pi$, $f\in \cI (p\varphi)_x$ and $g\in \cI (p\varphi')_x$ for some $x\in X$. The inequalities in (ii) and a change of variables $w=\pi(z)$ in the $L^2$ integrals yield
$$
\int_{\pi^{-1}(U_x)} \frac{|f\circ \pi|^2 |J|^2}{\prod\limits_i |s_i|^{2 p a_i }} < +\infty \quad\hbox{and}\quad \int_{\pi^{-1}(U_x)} \frac{|g\circ\pi|^2 |J|^2}{\prod\limits_i |s_i|^{2 p b_i}} < +\infty
\leqno(\label{integral-condition}{4.1.13})
$$
for some small open neighborhood $U_x$ of $x$.
Since $\sum\limits_i \Div (s_i )$ is normal crossing, (\ref{integral-condition})~implies that
$$\sum_i (p a_i -1)\ln |s_i|  \preccurlyeq \ln (|f\circ \pi  |) +\ln |J| \quad\hbox{and}\quad
\sum_i (p b_i -1)\ln |s_i|  \preccurlyeq \ln (|g\circ \pi  |)+\ln |J|.$$
Combining this with (i), we get
$$
\sum_i (p c_i -2)\ln |s_i|  \preccurlyeq \ln (|(f\cdot g)\circ \pi  |) +2\ln |J| .
\leqno(\label{add3}{4.1.14})
$$
Note that $J$ is independent of $p$, and $c_i >0$. (\ref{add3}) implies thus that,
when $p \rightarrow +\infty$, we can find a sequence $\varepsilon_p \rightarrow 0^+$, such that
$$
\sum_i p c_i (1-\varepsilon_p)\ln |s_i|  \preccurlyeq \ln |(f\cdot g)\circ \pi  |  .
\leqno(\label{sing-compare}{4.1.15})
$$
Since $f$ (respectively $g$) is an arbitrary element in $\cI (p\varphi)$ (respectively $\cI(p\varphi')$), by the construction of $\varphi_p$ and $\varphi'_p$, 
(\ref{sing-compare}) implies that
$$ \sum_i c_i  (1-\varepsilon_p) \ln |s_i|    \preccurlyeq  (\varphi_p +\varphi'_p ) \circ \pi.$$
Combining this with the fact that $ (1-\varepsilon_p) \widetilde\varphi_m \circ \pi\sim\sum\limits_i c_i  (1-\varepsilon_p) \ln |s_i| $, we get 
$$ (1-\varepsilon_p) \widetilde\varphi_m \circ \pi \preccurlyeq  (\varphi_p +\varphi'_p ) \circ \pi .$$
Therefore $ (1-\varepsilon_p ) \widetilde\varphi_m   \preccurlyeq  \varphi_p +\varphi'_p $ and (\ref{equivalproof}) is proved.\QED
\endproof

This motivates the following formal definition.

\claim{\label{weak-sing-relation}{4.1.16}.\ Definition} For each class $\{\alpha\}\in\cE(X)$, we define $\hatcS_\alpha(X)$ as a set of equivalence classes of sequences of quasi-positive currents $T_k=\alpha+dd^c\psi_k$ such that
\item{\rm(a)} $T_k=\alpha+dd^c\psi_k\geq  -\varepsilon_k\cdot \omega$
with $\lim_{k\to+\infty}\varepsilon_k=0$,
\vskip2pt
\item{\rm(b)}
the functions $\psi_k$ have analytic singularities and $\psi_k\preccurlyeq\psi_{k+1}$ for all $k$.
\noindent
We say that $(T_k)$ is weakly less singular than $(T_k')$ in 
$\hatcS_\alpha(X)$, and write $(T_k)\preccurlyeq_W(T_k')$,
 if for every $\varepsilon>0$ and $k$, there exists
$\ell$ such that~ $(1-\varepsilon)T_k\preccurlyeq T_\ell'$. Finally,
we write $(T_k)\sim_W(T_k')$ when we have $(T_k)\preccurlyeq_W(T_k')$
and $(T'_k)\preccurlyeq_W(T_k)$, and define $\hatcS_\alpha(X)$ to be 
the quotient space by this equivalence relation.
\endclaim

The set
$$
\hatcS(X)=\bigcup_{\{\alpha\}\in\cE(X)}\hatcS_\alpha(X)
\leqno(4.1.17)$$
is by construction a fiber space $\hat\pi:\hatcS(X)\to\cE(X)$,
and, by fixing a multiplicative sequence such as $m_k=2^k$, we find
a natural ``Bergman approximation functional''
$$\bfB:\cS(X)\to\hatcS(X),\qquad T=\alpha+dd^c\varphi\longmapsto(T_{B,k}),~~
T_k=\alpha+dd^c\psi_{B,k}
\leqno(\label{bergman-approx-map}{4.1.18})$$
where $\psi_{B,k}:=\varphi_{m_k}$ is the corresponding subsequence of the sequence of Bergman approximations $(\varphi_m)$.\vskip4pt

The set $\hatcS(X)$ is equipped with a natural addition $(T_k)+(T_k')=(T_k+T_k')$, with a scalar multiplication $\lambda\cdot(T_k)=(\lambda T_k)$ for $\lambda\in\bR_+$, as well as with the $\min(\bu)$ opera\-tion
$\min((T_k),(T_k'))=(\min(T_k,T_k'))$ obtained by taking $\max(\psi_k,\psi'_k)$ of the corres\-ponding potentials. By Prop.~\ref{prop-bergman-max-relation}, $\bfB$ is a morphism for the $\min(\bu)$ operation, and by Th.~\ref{theorem-bergman-additivity}, $\bfB$ is also a morphism for addition. Accordingly, it is natural to define a weak equivalence of singularities for closed positive currents by
$$
\leqalignno{
&T\preccurlyeq_W T' ~~\Longleftrightarrow_{\rm def}~~(T_{B,k})
\preccurlyeq_W (T'_{B,k}),&(\label{weak-inequality}{4.1.19})\cr
&T \sim_W T'~~\Longleftrightarrow~~
T\preccurlyeq_W T'~~\hbox{and}~~T'\preccurlyeq_W T.&(\label{weak-equivalence}{4.1.20})\cr}
$$
Related ideas are discussed in \cite{BFJ08} (especially \S~5), using the theory of valuations. One can summarize the above results in the following statement.

\claim{\label{bergman-singularity-isomorphism}{4.1.21}.\ Theorem} The Bergman approximation functional
$$\bfB:\cS(X)\to\hatcS(X),\qquad T=\alpha+dd^c\varphi\longmapsto(T_{B,k})$$
is a morphism for addition and for the $\min(\bu)$ opera\-tion on currents. Moreover $\bfB$~induces an injection $\cS(X)/{\sim_W}\to\hatcS(X)$.
\endclaim

\claim{4.1.22.\ Remark} {\rm It is easy to see that the induced map $\cS(X)/{\sim_W}\to\hatcS(X)$ is an isomorphism when $\dim X=1$. However, this map is not always surjective when $\dim X\ge 2$. In fact, Example~1.7 in \cite{DPS94} exhibits a ruled surface over an elliptic curve $\Gamma$ and a nef line bundle $L$ over $X$, such that $\alpha=c_1(L)$ contains a unique closed positive current $T=[C]$, for some curve $C\subset X$ that is a section of $X\to\Gamma$. Then the Bergman approximation is (up to equivalence of singularities) the constant sequence $T_{B,k}=T$, while $\hatcS_\alpha(X)$ also contains a sequence of smooth currents $T_k\ge -\varepsilon_k\omega$. This implies that $\cS(X)\to\hatcS(X)$ is not surjective in this situation. The following proposition shows however that the ``formal elements'' $(T_k)$ from $\hatcS(X)$ do not carry larger singularities than the closed positive current classes in $\cS(X)$ (the latter being constrained by the singularities of the ``limiting currents'' $T$ representing the class).}
\endclaim

\claim{4.1.23.\ Proposition} Let $T_k=\alpha+dd^c\psi_k$ be a sequence of closed $(1,1)$-currents representing an element in $\hatcS_\alpha(X)$. Then there exists a closed positive current
$T\in\alpha$ such that $(T_k)\preccurlyeq_W (T_{B,k})$.
\endclaim

\proof{Proof}  We have 
$T_k\ge-\varepsilon_k\omega$ and $\psi_k\preccurlyeq\psi_{k+1}$ for some 
decreasing sequence $\varepsilon_k\downarrow 0$. We replace $\psi_k$ by setting
$$
\widetilde\psi_k(x)=\sup\big\{\tau(x)\,;\;\sup_X\tau\le 0,~\alpha+dd^c\tau\ge-\varepsilon_k\omega,~\hbox{and}~\exists C>0,~\tau\le \psi_k+C\big\}.
$$
Then $(\widetilde\psi_k)$ is a decreasing sequence for the usual order relation~$\leq$ and $\widetilde\psi_k\sim\psi_k$ (the argument to prove the equivalence of singularities is similar to the one already used in the proof of Th.~\ref{theorem-bergman-additivity}, clearly $\widetilde\psi_k\ge\psi_k-M_k$ where $M_k=\sup_X\psi_k$, and the converse inequality $\widetilde\psi_k\le\psi_k+C_k$ is seen by using a blow-up to make the singularities of $\psi_k$ divisorial). We take
$$
\varphi=\lim_{k\to+\infty}\widetilde\psi_k\quad\hbox{and}\quad T=\alpha+dd^c\varphi.
$$
Since $\alpha+dd^c\widetilde\psi_k\ge-\varepsilon_k\omega$, we get in the limit
$T=\alpha+dd^c\varphi\ge 0$. Let $(\varphi_m)$ be the Bergman approximation sequence of $\varphi$. Since \hbox{$\varphi\le\widetilde\psi_\ell\le\psi_\ell+C_\ell$}, Prop.~\ref{sing-comparison}~(a) applied with $\lambda=2m$ shows that $\varphi_m\succcurlyeq \frac{1}{2mc_\ell}(2mc_\ell-n)_+\psi_\ell$ where $c_\ell>0$ is the coefficient of the log singularity of $\psi_\ell$. Therefore, if we take $T_{B,k}=\alpha+dd^c\varphi_{m_k}$, we get in the limit $(T_{B,k})\succcurlyeq_W (T_\ell)$.
\QED
\endproof

\claim{4.1.24.\ Remark} {\rm  When $X$ is projective algebraic and $\{\alpha\}$ belongs to the Neron-Severi space 
$${\rm NS}_\bR(X)=(H^{1,1}(X,\bC)\cap H^2(X,\bZ)/{\rm torsion})\otimes_\bZ\bR,$$
the fiber $\hatcS_\alpha(X)$ is essentially an algebraic object. In fact, we could
define $\hatcS_\alpha(X)$ as the set of 
suitable equivalence classes of ``formal limits''
$\lim_{c_1(D)\to\{\alpha\}}\lim_{k\to+\infty}\frac{1}{k}\fraka_k$
associated with sequences of graded
ideals $\fraka_k\subset H^0(X,\cO_X(kD)$ satisfying the subadditive 
property $\fraka_{k+\ell}\subset \fraka_k\fraka_\ell$, where $D$ are big
$\bQ$-divisors whose first Chern classes $c_1(D)$ approximate 
$\{\alpha\}\in{\rm NS}_\bR(X)$. Many related questions are
discussed in the algebraic setting in Lazarfeld's book
\cite{Laz04}. It is nevertheless an interesting point, even in the 
projective case, that one can ``extrapolate'' these concepts to all
trans\-cendental classes, and get in this way a global space $\hatcS(X)$ which
looks well behaved, e.g.\ semicontinuous, under variation of the
complex structure of~$X$.}
\endclaim

\subsection{\label{subsection-intersection-theory}{4.2}}{Intersection theory on $\cS(X)$ and $\widehat\cS(X)$}

Let $X$ be a compact K\"ahler $n$-dimensional manifold equipped with a K\"ahler metric~$\omega$. We consider closed positive $(1,1)$-currents
$T_j=\alpha_j+dd^c\varphi_j$, $1\le j\le p$. Let us first assume that
the functions $\varphi_j$ have analytic singularities, and let $Z\subset X$ 
be an analytic set such that the $\varphi_j$'s are locally bounded on 
$X\ssm Z$. The $(p,p)$-current 
$$
\Theta=\bfone_{X\ssm Z}T_1\wedge\ldots\wedge T_k
$$
is well defined on $X\ssm Z$, thanks to Bedford and Taylor \cite{BT76}, and
it is a closed positive current there. By \cite{BT76} such a current does not carry mass on any analytic set, so we can enlarge $Z$ without changing the total mass of~$\Theta$. In fact, $\Theta$ extends as a closed positive current on the whole of~$X$. To see this, let us take a simultaneous {\it log resolution} of the $T_j$'s, i.e.\ a modification
$$
\mu:\widehat X\to X
$$
such that if $\varphi_j=c_j\log\sum_\ell|g_{j,\ell}|^2+O(1)$, then
the pull-back of the ideals $(g_{j,\ell})_\ell$, namely 
$\mu^*(g_{j,\ell})_\ell=(g_{j,\ell}\circ \mu)_\ell$ is a purely divisorial
ideal sheaf $\cO_{\widehat X}(-D_j)$ on~$\widehat X$. Let $u_j=0$ be a local holomorphic equation of the divisor $D_j$ on $\widehat X$. Since
$\log\sum_\ell|g_{j,\ell}|^2=\log|u_j|^2+\log\sum_\ell|g_{j,\ell}/uj|^2=\log|u_j|^2+
v_j$, where $v_j\in C^\infty$ and $dd^c\log|u_j|^2=[D_j]$ by 
the Lelong-Poincar\'e equation, we find
$$
\mu^*T_j=\mu^*\alpha_j+dd^c(\varphi_j\circ\mu)=c_j[D_j]+\widehat T_j,\quad
\hbox{where}\quad
\widehat T_j=\mu^*\alpha_j+dd^c\widehat\varphi_j
\leqno(\label{log-res-current}{4.2.1})$$
and $\widehat\varphi_j$ is a locally bounded potential on~$\widehat X$ such that
$\widehat T_j\ge 0$. Now, if $E=\mu^{-1}(Z)$, we get
$$
\bfone_{X\ssm Z}T_1\wedge\ldots\wedge T_p=\mu_*(\bfone_{\widehat X\ssm E}
\widehat T_1\wedge\ldots\wedge \widehat T_p)=
\mu_*(\widehat T_1\wedge\ldots\wedge \widehat T_p).
\leqno(\label{direct-image-equation}{4.2.2})$$
Hence the right-hand side defines the desired extension of
$\bfone_{X\ssm Z}T_1\wedge\ldots\wedge T_p$ to $X$ as the direct
image of a closed positive current on $\widehat X$ carrying no mass on $E$.
An essential point is the following monotonicity lemma -- the reader will
find a more general version for non-pluripolar products in \cite{BEGZ},
Theorem 1.16.

\claim{\label{monotonicity-lemma}{4.2.3}.\ Lemma}
Assume that we have closed positive $(1,1)$-currents with
analytic singularities
$T_j$, $T'_j \in\{\alpha_j\}$ with $T_j\preccurlyeq T'_j$, $1\le j\le p$, 
and let $\gamma\ge 0$ be a closed positive smooth $(n-p,n-p)$-form on $X$.
If $Z$ is an analytic set containing the poles of all $T_j$ and $T'_j$, we have
$$
\int_X\bfone_{X\ssm Z}T_1\wedge\ldots\wedge T_p\wedge\gamma \ge
\int_X\bfone_{X\ssm Z}T'_1\wedge\ldots\wedge T'_p\wedge\gamma.
$$
\endclaim

\proof{Proof}  We take a log-resolution $\mu:\widehat X\to X$ that works
for all $T_j$ and $T'_j$ simultaneously. By (\ref{log-res-current}) and
(\ref{direct-image-equation}), 
we have $\mu^*T_j=c_j[D_j]+\widehat T_j$ where $\widehat T_j\ge 0$ 
has a locally bounded potential on $\widehat X$, and
$$
\int_X\bfone_{X\ssm Z}T_1\wedge\ldots\wedge T_p\wedge\gamma =
\int_{\widehat X}\widehat T_1\wedge\ldots\wedge\widehat T_p\wedge\mu^*\gamma.
$$
There are of course similar formulas $\mu^*T'_j=c_j[D'_j]+\widehat T'_j$
for the $T'_j$'s, and our assumption $T_j\preccurlyeq T'_j$ means that
the corresponding divisors satisfy $c_jD_j\le c'_jD'_j$, hence
$\Delta_j:=c'_jD'_j-c_jD_j\ge 0$. In terms of cohomology, we have
$$
\mu^*\{\alpha_j\}=\{\mu^*T_j\}=\{\widehat T_j\}+\{c_jD_j\}=
\{\mu^*T'_j\}=\{\widehat T'_j\}+\{c'_jD'_j\},
$$
hence $\{\widehat T_j\}=\{\widehat T'_j\}+\{\Delta_j\}$ in 
$H^2(\widehat X,\bR)$. By Stokes' theorem, we conclude that
$$\eqalign{
\int_{\widehat X}\widehat T_1\wedge\widehat T_2\wedge\ldots\wedge\widehat T_p\wedge\mu^*\gamma
&=\int_{\widehat X}(\widehat T'_1+\{\Delta_1\})\wedge\widehat T_2\wedge\ldots\wedge\widehat T_p\wedge\mu^*\gamma\cr
&\ge
\int_{\widehat X}\widehat T'_1\wedge\widehat T_2\wedge\ldots\wedge\widehat T_p\wedge\mu^*\gamma\cr}
$$
thanks to the positivity of our currents $\widehat T_j$, $\widehat T'_j$ and
the fact that the product of such currents with bounded potentials by the current of integration $[\Delta_j]$ is well defined and positive (\cite{BT76}). By replacing successively all terms $\{\widehat T_j\}$ by $\{\widehat T'_j\}+\{\Delta_j\}$ we infer
$$
\int_{\widehat X}\widehat T_1\wedge\ldots\wedge\widehat T_p\wedge\mu^*\gamma\ge
\int_{\widehat X}\widehat T'_1\wedge\ldots\wedge\widehat T'_p\wedge\mu^*\gamma.
\eqno\square
$$
\endproof

Now, assume that we have arbitrary closed positive $(1,1)$-currents
$T_1$, $\ldots$~, $T_p$. For each of them, we take a
sequence $T_{j,k}=\alpha_j+i\smash{\ddbar}\psi_{j,k}$ of monotone asymptotically equi\-singular approximations by currents with analytic singularities,
$T_{j,k}\ge -\varepsilon_{j,k}\omega$, $\lim_{k\to +\infty}\varepsilon_{j,k}=0$.
We have $T_{j,k}\preccurlyeq T_{j,k+1}$, and we may also assume without loss
of generality that
$\varepsilon_{j,k}\ge\varepsilon_{j,k+1}>0$ for all $j,k$. Let $Z_k$ be an 
analytic containing all poles of the $T_{j,k}$, $1\le j\le p$. It follows 
immediately from the above discussion and especially from
Lemma~\ref{monotonicity-lemma} that the integrals
$$
\int_X\bfone_{X\ssm Z_k}(T_{1,k}+\varepsilon_{1,k}\omega)\wedge\ldots\wedge
(T_{p,k}+\varepsilon_{p,k}\omega)\wedge\gamma \ge 0
$$
are well defined and nonincreasing in $k$ (the fact that $\varepsilon_{j,k}$
is non increasing even helps here). From this, we conclude

\claim{\label{intersection-product}{4.2.4}.\ Theorem} For every $p=1,\,2,\ldots,\,n$, there is a well defined $p$-fold intersection product
$$
\hatcS(X)\times\cdots\times\hatcS(X)\longrightarrow
H^{p,p}_+(X,\bR)
$$
which assigns to any $p$-tuple of equivalence classes of 
monotone sequences $(T_{j,k})$ in $\hatcS(X)$, \hbox{$1\le j\le p$},
the limit cohomology class
$$
\lim_{k\to+\infty}\big\{\bfone_{X\ssm Z_k}(T_{1,k}+\varepsilon_{1,k}\omega)\wedge\ldots\wedge(T_{p,k}+\varepsilon_{p,k}\omega)\big\}\in H^{p,p}_+(X,\bR)
$$
where $H^{p,p}_+(X,\bR)\subset H^{p,p}(X,\bR)$ denotes the cone of cohomology
classes of closed positive $(p,p)$-currents. This product is additive and 
homogeneous in each argument in the space $\hatcS(X)$. 
\endclaim

\claim{4.2.5.\ Corollary} By combining the above formal intersection product with the Bergman approximation operator $\bfB:\cS(X)\to \hatcS(X)$, we get an intersection product
$$
\cS(X)\times\cdots\times\cS(X)\longrightarrow
H^{p,p}_+(X,\bR)\quad\hbox{denoted}\quad
(T_1,\ldots,T_p)\longmapsto\langle T_1,\ldots,T_p\rangle^+,
$$
which is homogeneous and additive in each argument.
\endclaim

\proof{Proof of Th.~\ref{intersection-product}}  The existence of a limit 
in cohomology is seen by fixing a dual basis $(\{\gamma_j\})$ of 
$H^{n-p,n-p}(X)$, using the Serre duality pairing
$$
H^{p,p}(X,\bR)\times H^{n-p,n-p}(X)\to\bR,\qquad(\beta,\gamma)
\mapsto\int_X\beta\wedge\gamma.
$$
Since $X$ is K\"ahler, we can take $\gamma_1=\omega^{n-p}$ and replace if
necessary $\gamma_j$ by $\gamma_j+C\omega^{n-p}$, $C\gg 1$, to get $\gamma_j\ge 0$ for all $j\ge 2$. Then the integrals 
$$
\int_X\bfone_{X\ssm Z_k}(T_{1,k}+\varepsilon_{1,k}\omega)\wedge\ldots\wedge
(T_{p,k}+\varepsilon_{p,k}\omega)\wedge\gamma_j\ge 0
$$
are nonincreasing in $k$, and the limit must therefore exist by 
monotonicity.\QED
\endproof

\claim{4.2.6.\ Remark} {\rm It is natural to ask how the above
intersection product compares with the (cohomology class of the)
``non-pluripolar product'' $\langle T_1,\ldots,T_p\rangle$ defined in
\cite{BEGZ}, \S$\,$1. In fact, the above product only neglects
analytic parts of the currents involved. The simple example of a
probability measure $T$ without atoms supported on a polar set of a
compact Riemann surface $X$ yields e.g.\ $\langle T\rangle^+ = 1$,
while the non-pluripolar part $\langle T\rangle$ vanishes.}
\endclaim

\subsection{4.3}{K\"ahler definition of the numerical dimension}

Using the intersection product defined in Th.~\ref{intersection-product}, 
we can give a precise definition of the numerical dimension. 

\claim{\label{numerical-dimension}{4.3.1}.\ Definition}
Let $(X,\omega)$ be a compact K\"ahler
$n$-dimensional manifold. We define the numerical dimension $\nd(T)$ of a closed positive $(1,1)$-current~$T$ on $X$ to be the largest integer $p=0,1,\ldots,n$ such that
$\langle T^p\rangle^+\ne 0$, i.e.\ $\int_X\langle T^p\rangle^+\wedge\omega^{n-p}>0$.
\endclaim

Accordingly, if $(L,h )$ be a pseudoeffective line bundle on $X$, we define
its nume\-rical dimension to be 
$$
\nd(L,h)=\nd(\ii\Theta_{L,h}).
\leqno(\label{numerical-dimension-line-bundle}{4.3.2})$$
\noindent
By the results of the preceding subsection, $\nd(L,h)$ depends only on the 
weak equi\-valence class of singularities of the metric~$h$.

\claim{4.3.3.\ Remark} {\rm 
H.~Tsuji \cite{Tsu07} has defined a notion of numerical dimension by a
more algebraic method:}
\endclaim

\claim{4.3.4.\ Definition} 
Let $ X $ be a projective variety and $(L,h) $ a pseudo-effective line bundle.
When $V$ runs over all irreducible algebraic suvarieties of $X$, one defines
$$\nu_{\num}(L,h)=\sup\Big\{p=\dim V\,;\;
\limsup_{m\to \infty}\frac{h^{0}\big(\widetilde{V},\mu^*(L^{\otimes m})\otimes
\cI(\mu^*h^m)\big)}{m^p}>0\Big\}$$
where $\mu:\widetilde V\to V\subset X$ is an embedded desingularization of $V$
in~$X$.
\endclaim

Junyan Cao \cite{JC14} has shown that $\nu_{\num}(L,h)$ coincides with
$\nd(L,h)$ as defined in (\ref{numerical-dimension-line-bundle}). The idea is to make a reduction to the ``big'' case $\nd(L,h)=\dim X$ and to use holomorphic Morse inequalities \cite{Dem85b} in combination with a regularization procedure. We omit the rather technical details here.

\claim{4.3.5.\ Remark} {\rm If $L$ is pseudo-effective, there is also a natural 
concept of numerical dimension $\nd(L)$ that does not depend on the choice
of a metric~$h$ on~$L$. One can set e.g.
$$\eqalign{
\nd(L)=\max\Big\{p\in[0,n]\,;\;\exists c>0,~\forall \varepsilon&>0,~
\exists h_\varepsilon,~~
\Theta_{L,h_\varepsilon}\ge -\varepsilon\omega,~~\hbox{such that}\cr
&\int_{X\ssm Z_\varepsilon}(\ii\Theta_{L,h_\varepsilon}+\varepsilon\omega)^p\wedge\omega^{n-p}\ge c\Big\},\cr}
$$
where $h_\varepsilon$ runs over all metrics with analytic singularities on $L$.
It may happen in general that $\nd(L,h_{\min})<\nd(L)$, even when $L$ is nef; in that case the $h_\varepsilon$ can be taken to be smooth in the definition of $\nd(L)$, and therefore $\nd(L)$ is the largest integer $p$ such that $c_1(L)^p\ne 0$. In fact, for the line bundle $L$ already mentioned in Remark~\ref{nef-weird-example}, it is shown in \cite{DPS94} that there is unique positive current $T\in c_1(L)$, namely the current of integration $T=[C]$ on the negative curve $C\subset X$, hence
$\nd(L,h_{\min})=\nd([C])=0$, although we have $\nd(L)=1$ here.}
\endclaim

\section{\label{section-cao-proof}{5}}{Proof of Junyan Cao's vanishing theorem}

This section is a brief account and a simplified exposition of Junyan Cao's proof, as detailed in his PhD thesis \cite{JC13}. The key curvature and singularity estimates are contained in the following technical statement, which depends in a crucial way on Bergman regularization and on Yau's theorem \cite{Yau78} for solutions of Monge-Amp\`ere equations.

\claim{\label{curvature-estimates}{5.1}.\ Proposition}
Let $(L,h)$ be a pseudoeffective line bundle on a compact K\"ahler manifold~$(X,\omega)$. Let us write $T=\frac{\ii}{2\pi}\Theta_{L,h}=\alpha+dd^c\varphi$ where $\alpha$ is smooth and $\varphi$ is a quasi-psh potential. Let $p=\nd(L,h)$ be the numerical dimen\-sion of~$(L,h)$. Then, for every $\gamma\in{}]0,1]$ and $\delta\in{}]0,1]$, there exists a quasi-psh potential $\Phi_{\gamma,\delta}$ on $X$ satisfying the following properties$\,:$
\item{\rm(a)} $\Phi_{\gamma,\delta}$ is smooth in the complement $X\ssm Z_\delta$ of an analytic set $Z_\delta\subset X$.
\vskip2pt
\item{\rm(b)} $\alpha+\delta\omega+dd^c\Phi_{\gamma,\delta}\ge \frac{\delta}{2}(1-\gamma)\omega$ on $X$.
\vskip2pt
\item{\rm(c)} $(\alpha+\delta\omega+dd^c\Phi_{\gamma,\delta})^n\ge a\,\gamma^n\delta^{n-p}\omega^n$ on $X\ssm Z_\delta$.
\vskip2pt
\item{\rm(d)} $\Phi_{\gamma,\delta}\le (1+b\delta)\psi_{B,k}+C_{\gamma,\delta}$ where $\psi_{B,k}\ge \varphi$ is a Bergman approximation of $\varphi$ of sufficiently high index $k=k_0(\delta)$.
\vskip2pt
\item{\rm(e)} $\sup_X\Phi_{1,\delta}=0$, and for all $\gamma\in{}]0,1]$ there are estimates $\Phi_{\gamma,\delta}\le A$ and
$$\exp\big(-\Phi_{\gamma,\delta}\big)\le 
e^{-(1+b\delta)\varphi}\exp\big(A-\gamma\Phi_{1,\delta}\big)$$
\item{\rm(f)} For $\gamma_0,\,\delta_0>0$ small,
$\gamma\in{}]0,\gamma_0]$, $\delta\in{}]0,\delta_0]$ and $k=k_0(\delta)$ large
enough, we have
$$\cI(\Phi_{\gamma,\delta})=\cI_+(\varphi)=\cI(\varphi).$$
Here $a,\,b,\,A,\,\gamma_0,\,\delta_0,\,C_{\gamma,\delta}>0$ are suitable constants $(C_{\gamma,\delta}$ being the only one that depends on $\gamma$, $\delta)$.
\endclaim

\proof{Proof} 
Denote by $\psi_{B,k}$ the nonincreasing sequence of Bergman approximations of~$\varphi$ (obtained with denominators $m_k=2^k$, say). We have $\psi_{B,k}\ge\varphi$ for all $k$, the $\psi_{B,k}$ have analytic singularities and $\alpha+dd^c\psi_{B,k}\ge-\varepsilon_k\omega$ with $\varepsilon_k\downarrow 0$. Then $\varepsilon_k\le\frac{\delta}{4}$ for $k\ge k_0(\delta)$ large enough, and so
$$\eqalign{
\alpha+\delta\omega+dd^c\big((1+b\delta)\psi_{B,k}\big)
&\ge\alpha+\delta\omega-(1+b\delta)(\alpha+\varepsilon_k\omega)\cr
&\ge\delta\omega-(1+b\delta)\varepsilon_k\omega-b\delta\alpha\ge
{\textstyle\frac{\delta}{2}}\omega\cr}
$$
for $b>0$ small enough (independent of $\delta$ and $k$).
Let $\mu:\widehat X\to X$ be a log-resolution of $\psi_{B,k}$, so that 
$$\mu^*\big(\alpha+\delta\omega+dd^c((1+b\delta)\psi_{B,k})\big)=c_k[D_k]+\beta_k$$
where $\smash{\beta_k\ge\frac{\delta}{2}\mu^*\omega}\ge 0$ is a smooth closed $(1,1)$-form on $\widehat X$ that is${}>0$ in the complement $\widehat X\ssm E$ of the exceptional divisor, $c_k=\frac{1+b\delta}{m_k}>0$, and $D_k$ is a divisor that includes all components $E_\ell$ of $E$. The map $\mu$ can be obtained by Hironaka \cite{Hir64} as a composition of a sequence of blow-ups with smooth centers, and we can even achieve that $D_k$ and $E$ are normal crossing divisors. In this circumstance, it is well known that there exist arbitrary small numbers $\eta_\ell>0$ such that $\beta_k-\sum\eta_\ell[E_\ell]$ is a K\"ahler class on~$\widehat X$. Hence we can find a quasi-psh potential $\widehat\theta_k$ on $\widehat X$ such that $\widehat\beta_k:=\beta_k-\sum\eta_\ell[E_\ell]+dd^c\widehat\theta_k$ is a K\"ahler metric on~$\widehat X$, and by taking the $\eta_\ell$ small enough, we may assume that
$\smash{\int_{\widehat X}(\widehat\beta_k)^n\ge
\frac{1}{2}\int_{\widehat X}\beta_k^n}$. Now, we write
$$\eqalign{
\alpha+\delta\omega+dd^c\big((1+b\delta)\psi_{B,k}\big)
&\ge\alpha+\varepsilon_k\omega+dd^c\psi_{B,k}+(\delta-\varepsilon_k)\omega
-b\delta(\alpha+\varepsilon_k\omega)\cr
&\ge(\alpha+\varepsilon_k\omega+dd^c\psi_{B,k})+\textstyle{\frac{\delta}{2}}\omega\cr}
$$
for $k\ge k_0(\delta)$ and $b>0$ small (independent of $\delta$ and $k$). The assumption on the numerical dimension of $\frac{\ii}{2\pi}\Theta_{L,h}=\alpha+dd^c\varphi$ implies the existence of a constant $c>0$ such that, with $Z=\mu(E)\subset X$, we have
$$\eqalign{
\int_{\widehat X}\beta_k^n&=
\int_X\bfone_{X\ssm Z}\big(\alpha+\delta\omega+dd^c((1+b\delta)\psi_{B,k})\big)^n\cr
&\ge{n\choose p}\Big(\frac{\delta}{2}\Big)^{n-p}\int_{X\ssm Z}
\big(\alpha+\varepsilon_k\omega+dd^c\psi_{B,k}\big)^p\wedge\omega^{n-p}
\ge c\,\delta^{n-p}\int_X\omega^n\cr}
$$
for all $k\ge k_0(\delta)$. Therefore, we may assume
$$\int_{\widehat X}(\widehat\beta_k)^n\ge \frac{c}{2}\,\delta^{n-p}\int_X\omega^n.$$
By Yau's theorem \cite{Yau78}, there exists a quasi-psh potential $\widehat\tau_k$ on $\widehat X$ such that 
\hbox{$\widehat\beta_k+dd^c\widehat\tau_k$} is a K\"ahler metric on $\widehat X$ with a prescribed volume form $\widehat f>0$ such that $\int_{\widehat X}f=\int_{\widehat X}\widehat\beta_k^n$. By the above discussion, we can take here $\widehat f>\frac{c}{3}\delta^{n-p}\mu^*\omega^n$ everywhere on $\widehat X$.

Now, we consider $\theta_k=\mu_*\widehat\theta_k$ and
$\tau_k=\mu_*\widehat\tau_k\in L^1_{\loc}(X)$. Since $\widehat\theta_k$ was defined in such a way that
$dd^c\widehat\theta_k=\widehat\beta_k-\beta_k+\sum_\ell\eta_\ell[E_\ell]$, we get
$$\eqalign{
\mu^*\big(\alpha+\delta\omega+dd^c((1&+b\delta)\psi_{B,k}+
\gamma(\theta_k+\tau_k))\big)\cr
&=c_k[D_k]+(1-\gamma)\beta_k+\gamma\Big(
\sum_\ell\eta_\ell[E_\ell]+\widehat\beta_k+dd^c\widehat\tau_k\Big)\ge 0.\cr}
$$
This implies in particular that $\Phi_{\gamma,\delta}:=(1+b\delta)\psi_{B,k}+\gamma(\theta_k+\tau_k)$ is a quasi-psh potential on $X$ and that
$$\mu^*\big(\alpha+\delta\omega+dd^c\Phi_{\gamma,\delta}\big)\ge(1-\gamma)\beta_k\ge\frac{\delta}{2}(1-\gamma)\,\mu^*\omega,$$
thus condition (b) is satisfied. Putting $Z_\delta=\mu(|D_k|)\supset \mu(E)=Z$, we also have
$$
\mu^*\bfone_{X\ssm Z_\delta}\big(\alpha+\delta\omega+dd^c\Phi_{\gamma,\delta}\big)^n\ge\gamma^n\,\widehat\beta_k^n\ge \frac{c}{3}\,\gamma^n\delta^{n-p}\mu^*\omega^n,
$$
therefore condition (c) is satisfied as well with $a=c/3$. Property (a) is clear, and (d) holds since the quasi-psh function $\widehat\theta_k+\widehat\tau_k$ must be bounded from above on~$\widehat X$. We will actually adjust constants in $\widehat\theta_k+\widehat\tau_k$ (as we may), so that $\sup_X\Phi_{1,\delta}=0$. 
Since $\varphi\le\psi_{B,k}\le\psi_{B,0}\le A_0:=\sup_X\psi_{B,0}$ and
$$
\Phi_{\gamma,\delta}=(1+b\delta)\psi_{B,k}+\gamma
\big(\Phi_{1,\delta}-\psi_{B,k}\big)=
(1-\gamma+b\delta)\psi_{B,k}+\gamma\Phi_{1,\delta},
$$
we have
$$(1+b\delta)\varphi-\gamma(A_0-\psi_{B,k})\le
\Phi_{\gamma,\delta}\le(1-\gamma+b\delta)A_0$$
and the estimates in (e) follow with $A=(1+b)A_0$. The only remaining
property to be proved is~(f). Condition (d) actually implies 
$\cI(\Phi_{\gamma,\delta})\subset\cI((1+b\delta)\psi_{B,k})$, and 
Cor.~\ref{multiplier-ideals-comparison} also gives
$\cI((1+b\delta)\psi_{B,k})\subset \cI((1+b\delta/2)\varphi)$ if we take
$k\ge k_0(\delta)$ large enough, hence $\cI(\Phi_{\gamma,\delta})\subset
\cI_+(\varphi)$ for $\delta\le\delta_0$ small. In the opposite direction,
we observe that $\Phi_{1,\gamma}$ satisfies $\alpha+\omega+dd^c\Phi_{1,\delta}\ge 0$ and $\sup_X\Phi_{1,\delta}=0$, hence $\Phi_{1,\delta}$ belongs to a compact 
family  of quasi-psh functions. A standard result of potential theory then
shows the existence of a uniform small constant $c_0>0$ such that $\int_X
\exp(-c_0\Phi_{1,\delta})dV_\omega<+\infty$ for all $\delta\in{}]0,1]$.
If $f\in\cO_{X,x}$ is a germ of holomorphic function and $U$ a small
neighborhood of~$x$, the H\"older inequality combined with estimate~(e)
implies
$$
\int_U|f|^2\exp(-\Phi_{\gamma,\delta})dV_\omega\le
e^A\Big(\int_U|f|^2e^{-p(1+b\delta)\varphi}dV_\omega\Big)^{\frac{1}{p}}
\Big(\int_U|f|^2e^{-q\gamma\Phi_{1,\delta}}dV_\omega\Big)^{\frac{1}{q}}.
$$
We fix $\lambda_0>1$ so that $\cI(\lambda_0\varphi)=\cI_+(\varphi)$,
$p\in{}]1,\lambda_0[$ (say $p=1+\lambda_0)/2$),  and take
$$\gamma\le\gamma_0:=\frac{c_0}{q}=c_0\frac{\lambda_0-1}{\lambda_0+1}\quad
\hbox{and}\quad\hbox{$\delta \le\delta_0\in{}]0,1]$ so small that 
$p(1+b\delta_0)\le\lambda_0$.}$$
Then clearly $f\in\cI(\lambda_0\varphi)$ implies $f\in\cI(\Phi_{\gamma,\delta})$, and
(f) is proved.\QED
\endproof

The rest of the arguments proceeds along the lines of \cite{Dem82}, \cite{Mou95} and \cite{DP02}. Let $(L,h)$ be a pseuffective line bundle and $p=\nd(L,h)=\nd(\ii\Theta_{L,h})$. We equip $L$ be the hermitian metric $h_\delta$ defined by the quasi-psh weight $\Phi_\delta=\Phi_{\gamma_0,\delta}$ obtained in Prop.~\ref{curvature-estimates}, with $\delta\in{}]0,\delta_0]$. Since $\Phi_\delta$ is smooth on $X\ssm Z_\delta$, the well-known Bochner-Kodaira identity shows that for 
every smooth $(n,q)$-form $u$ with values in $K_X\otimes L$ that is 
compactly supported on $X\ssm Z_\delta$, one has
$$\Vert\dbar u\Vert_\delta^2+\Vert\dbar^* u\Vert_\delta^2\ge
2\pi\int_X(\lambda_{1,\delta}+\ldots+\lambda_{q,\delta}-q\delta)|u|^2e^{-\Phi_\delta}dV_\omega,$$
where $\Vert u\Vert_\delta^2:=\int_X|u|^2_{\omega,h_\delta}dV_\omega=
\int_X|u|^2e^{-\Phi_\delta}dV_\omega$ and
$$
0<\lambda_{1,\delta}(x)\le\ldots\le\lambda_{n,\delta}(x)
$$
are, at each point $x\in X$, the eigenvalues of $\alpha+\delta\omega+dd^c\Phi_\delta$ with respect to the base K\"ahler metric~$\omega$. Notice that the 
$\lambda_{j,\delta}(x)-\delta$ 
are the actual eigenvalues of $\frac{\ii}{2\pi}\Theta_{L,h_\delta}=\alpha+
dd^c\Phi_\delta$ with respect to $\omega$ and that the inequality $\lambda_{j,\delta}(x)\ge\frac{\delta}{2}(1-\gamma)>0$ is guaranted by Prop.~\ref{curvature-estimates}~(b). After dividing by $2\pi q$ (and neglecting that constant in the left hand side), we get
$$
\Vert\dbar u\Vert_\delta^2+\Vert\dbar^* u\Vert_\delta^2+\delta\Vert u\Vert_\delta^2\ge
\int_X(\lambda_{1,\delta}+\ldots+\lambda_{q,\delta})|u|^2e^{-\Phi_\delta}dV_\omega.
\leqno(5.2)$$
A standard Hahn-Banach argument in the $L^2$-theory of the $\dbar$-operator
then yields the following conclusion.

\claim{\label{Ltwo-estimate}{5.3}.\ Proposition}
For every $L^2$ section of $\Lambda^{n,q}T^*_X
\otimes L$ such that $\Vert f\Vert_\delta<+\infty$ and $\dbar f=0$ in the sense of distributions, there exists a $L^2$ section $v=v_\delta$ of $\Lambda^{n,q-1}T^*_X\otimes L$ and a $L^2$ section $w=w_\delta$ of $\Lambda^{n,q}T^*_X\otimes L$ such that $f=\dbar v+w$ with
$$
\Vert v\Vert_\delta^2+\frac{1}{\delta}\Vert w\Vert_\delta^2\le
\int_X\frac{1}{\lambda_{1,\delta}+\ldots+\lambda_{q,\delta}}|f|^2e^{-\Phi_\delta}dV_\omega.
$$
\endclaim
Because of the singularities of the weight on $Z_\delta$, one should in fact argue first on \hbox{$X\ssm Z_\delta$} and approximate the base K\"ahler metric $\omega$ by a metric $\widehat\omega_{\delta,\varepsilon}=\omega+\varepsilon\widehat\omega_\delta$ that is complete on $X\ssm Z_\delta$, exactly as explained in \cite{Dem82}; we omit the (by now standard) details here. A consequence of Prop.~\ref{Ltwo-estimate} is that the ``error term'' $w$ satisfies the $L^2$ bound
$$\int_X|w|^2e^{-\Phi_\delta}dV_\omega\le\int_X\frac{\delta}{\lambda_{1,\delta}+\ldots+\lambda_{q,\delta}}|f|^2e^{-\Phi_\delta}dV_\omega.
\leqno(\label{error-estimate}{5.4})$$

The idea for the next estimate is taken from Mourougane's PhD thesis \cite{Mou95}.

\claim{\label{ratio-lemma}{5.5}.\ Lemma}
The ratio $\rho_\delta(x):=\delta/(\lambda_{1,\delta}(x)+\ldots+\lambda_{q,\delta}(x))$ is uniformly bounded on~$X$ $($independently of $\delta)$, and, as soon as $q\ge n-\nd(L,h)+1$, there exists a subsequence $(\rho_{\delta\ell})$, $\delta_\ell\to 0$, that tends almost everywhere to $0$ on~$X$.
\endclaim

\proof{Proof} 
 By estimates (b,c) in Prop.~\ref{curvature-estimates}, we have $\lambda_{j,\delta}(x)\ge \frac{\delta}{2}(1-\gamma_0)$ and
$$\lambda_{1,\delta}(x)\ldots \lambda_{n,\delta}(x)\ge a\gamma_0^n\delta^{n-p}\quad
\hbox{where}~~p=\nd(L,h).
\leqno(\label{product-eigenvalues-bound}{5.6})$$
Therefore we already find $\rho_\delta(x)\le 2/q(1-\gamma_0)$. Now, we have
$$
\int_{X\ssm Z_\delta}\lambda_{n,\delta}(x)dV_\omega\le
\int_X(\alpha+\delta\omega+dd^c\Phi_\delta)\wedge\omega^{n-1}
=\int_X(\alpha+\delta\omega)\wedge\omega^{n-1}
\le{\rm Const},
$$
therefore the ``bad set'' $S_\varepsilon\subset X\ssm Z_\delta$ of points $x$ where $\lambda_{n,\delta}(x)>\delta^{-\varepsilon}$ has a volume Vol$(S_\varepsilon)\le C\delta^\varepsilon$ converging to~$0$ as $\delta\to 0$ (with a slightly more elaborate argument we could similarly control any elementary symmetric function in the $\lambda_{j,\delta}$'s, but this is not needed here). Outside of $S_\varepsilon$, the inequality (\ref{product-eigenvalues-bound}) yields
$$
\lambda_{q,\delta}(x)^q\delta^{-\varepsilon(n-q)}\ge \lambda_{q,\delta}(x)^q\lambda_{n,\delta}(x)^{n-q}\ge a\gamma_0^n\delta^{n-p}
$$
hence
$$
\lambda_{q,\delta}(x)\ge c\delta^{\frac{n-p+(n-q)\varepsilon}{q}}
\quad\hbox{and}\quad 
\rho_\delta(x)\le C\delta^{1-\frac{n-p+(n-q)\varepsilon}{q}}.
$$
If we take $q\ge n-p+1$ and $\varepsilon>0$ small enough, the exponent of $\delta$ in the final estimate is positive, and Lemma~\ref{ratio-lemma} follows.\QED
\endproof

\proof{Proof of Junyan Cao's Theorem, Th.~\ref{kvn-cao-theorem}} Let $\{f\}$ be a cohomology class in the group $H^q(X,K_X\otimes L\otimes\cI_+(h))$, $q\ge n-\nd(L,h)+1$. Consider a finite Stein open covering $\cU=(U_\alpha)_{\alpha=1,\ldots,N}$ by coordinate balls $U_\alpha$. There is an isomorphism between {\v C}ech cohomology $\check H^q(\cU,\cF)$ with values in the sheaf $\cF=\cO(K_X\otimes L)\otimes\cI_+(h)$ and the cohomology of the complex $(K^{\bu}_\delta,\dbar)$ of $(n,q)$-forms $u$ such that both $u$ and $\dbar u$ are $L^2$ with respect to the weight~$\Phi_\delta$, i.e.\ $\int_X|u|^2\exp(-\Phi_\delta)dV_\omega<+\infty$ and
$\int_X|\dbar u|^2\exp(-\Phi_\delta)dV_\omega<+\infty$. The isomorphism comes from Leray's theorem and from the fact that the sheafified complex $(\cK^{\bu}_\delta,\dbar)$ is a complex of $\cC^\infty$-modules that provides a resolution of the sheaf $\cF\,$: the main point here is that $\cI(\Phi_\delta)=\cI_+(\varphi)=\cI_+(h)$, as asserted by Prop.~\ref{curvature-estimates}~(f), and that we can locally solve $\dbar$-equations by means of H\"ormander's estimates \cite{H\"or66}.

Let $(\psi_\alpha)$ be a partition of unity subordinate to~$\cU$.
The explicit isomorphism between {\v C}ech cohomology and $L^2$ cohomology yields a smooth $L^2$ representative $f=\sum_{|I|=q}f_I(z)dz_1\wedge...\wedge dz_n\wedge d\overline z_I$ which is a combination
$$f=\sum_{\alpha_0}\psi_{\alpha_0}c_{\alpha_0\alpha_1\ldots\alpha_q}\dbar\omega_{\alpha_1}\wedge\ldots\wedge\dbar\psi_{\alpha_q}$$
of the components of the corresponding {\v C}ech cocycle
$$c_{\alpha_0\alpha_1\ldots\alpha_q}\in\Gamma\big(
U_{\alpha_0}\cap U_{\alpha_1}\cap\ldots\cap U_{\alpha_q},\cO(\cF)\big).
$$
Estimate (e) in Prop.~\ref{curvature-estimates} implies the
H\"older inequality
$$
\int_X\rho_\delta|f|^2\exp(-\Phi_{\delta})dV_\omega\le
e^A\Big(\int_X\rho_\delta^p|f|^2e^{-p(1+b\delta)\varphi}dV_\omega\Big)^{\!\frac{1}{p}}
\Big(\int_X|f|^2e^{-q\gamma_0\Phi_{1,\delta}}dV_\omega\Big)^{\!\frac{1}{q}}\kern-3pt.
$$
Our choice of $\delta\le\delta_0$, $\gamma_0$ and $p,\,q$ shows that the integrals in the right hand side are convergent, and especially
$\int_X|f|^2e^{-p(1+b\delta)\varphi}dV_\omega<+\infty$. 
Lebesgue's dominated convergence theorem combined with Lemma~\ref{ratio-lemma} implies that the $L^p$-part goes to~$0$ as $\delta=\delta_\ell\to 0$, hence the ``error term'' $w$ converges to $0$ in $L^2$ norm by estimate~(\ref{error-estimate}). If~we express the corresponding class $\{w\}$ in {\v C}ech cohomology and use H\"ormander's estimates on the intersections $U_\alpha=\bigcap\smash{U_{\alpha_j}}$, we see that $\{w\}$ will be given by a {\v C}ech cocycle $(\widetilde w_\alpha)$ such that $\int_{U_\alpha}|\widetilde w_\alpha|^2e^{-\Phi_\delta}dV_\omega\to 0$ as $\delta=\delta_\ell\to 0$ (we may lose here some fixed constants since $\Phi_\delta$ is just quasi-psh on our balls, but this is irrelevant thanks to the uniform lower bounds for the Hessian). The inequa\-lity $\Phi_\delta\le A$ in
Prop.~\ref{curvature-estimates}~(e) shows that we have as well
an unweighted $L^2$ estimate $\smash{\int_{U_\alpha}}|\widetilde w_\alpha|^2dV\to 0$.
However it is well-known that when one takes unweighted $L^2$ norms on spaces of {\v Cech} cocyles (or uniform convergence on compact subsets, for that purpose), the resulting topology on the finite dimensional space $\check H^q(\cU,\cF)$ is Hausdorff, so the subspace of coboundaries is closed in the space of cocycles. Hence we conclude from the above that $f$ is a coboundary, as desired.\QED
\endproof

\claim{5.7.\ Remark} {\rm In this proof, it is remarkable that one can control the error term $w$, but a priori completely lose control on the element $v$ such that $\dbar v\approx f$ when $\delta\to 0\,$!}
\endclaim

\section{6}{Compact K\"ahler threefolds without nontrivial subvarieties}

The bimeromorphic classification of compact K\"ahler manifolds leads to 
considering those, termed as ``simple'', that have as little internal 
structure as possible, and are somehow the elementary bricks needed to reconstruct all others through meromorphic fibrations (cf.\ \cite{Cam80}, \cite{Cam85}).

\claim{6.1.\ Definition} A compact K\"ahler manifold $X$ is said to be simple if there does not exist any irreducible analytic subvariety $Z$ with $0<\dim Z<\dim X$ through a very general point $x\in X$, namely a point $x$ in the complement $X\ssm\bigcup S_j$ of a countable union of analytic sets $S_j\subsetneq X$.
\endclaim

Of course, every one dimensional manifold $X$ is simple, but in higher dimensions $n>1$, one can show that a very general torus $X=\bC^n/\Lambda$ has no nontrivial analytic subvariety $Z$ at all (i.e.\ none beyond finite sets and $X$ itself), in any dimension~$n$. In even dimension, a very general Hyperk\"ahler manifold can be shown to be simple as well. It has been known since Kodaira that there are no other simple K\"ahler surfaces (namely only very general 2-dimensional tori and K3 surfaces). Therefore, the next dimension to be investigated is dimension~$3$. In this case, Campana, H\"oring and Peternell have shown in \cite{CHP14} that $X$ is bimeromorphically a quotient of a torus by a finite group (see Theorem~\ref{theorem-kummer} at the end). Following \cite{CDV13}, we give here a short self-contained proof for ``strongly simple'' K\"ahler threefolds, namely threefolds that do not possess any proper analytic subvariety.

\claim{6.2.\ Theorem}  {\rm(\cite{Bru10})} Let $X$ be a compact K\"ahler manifold with a $1$-dimensional holomorphic foliation $F$ given by a nonzero morphism of vector bundle $L\to T_X$, where $L$ is a line bundle on $X$, and $T_X$ is its holomorphic tangent bundle. If $L^{-1}$ is not pseudoeffective, the closures of the leaves of $F$ are rational curves, and $X$ is thus uniruled.
 \endclaim

We use this result in the form of the following corollary, which has been observed in \cite{HPR11}, Proposition~4.2.

\claim{6.3.\ Corollary}  If $X$ is a non uniruled $n$-dimensional compact K\"ahler manifold with\break $H^0(X,\Omega^{n-1}_X)\neq 0$, then $K_X$ is pseudoeffective.
\endclaim

\proof{Proof}  $\Omega^{n-1}_X$ is canonically isomorphic to $K_X\otimes T_X$. Any nonzero section of $\Omega^{n-1}_X$ thus provides a nonzero map $K_X^{-1}\to T_X$, and an associated foliation.\QED
\endproof

It follows from the above that the canonical line bundle $K_X$ of our simple threefold~$X$ must be pseudoeffective. We then use the following simple observation.

\claim{6.4.\ Proposition} Assume that $X$ is a strongly simple compact complex manifold. Then every pseudoeffective line bundle $(L,h)$ is nef, and all multiplier sheaves  $\cI(h^m)$ are trivial, i.e.\ $\cI(h^m)=\cO_X$. Moreover, we have $c_1(L)^n=0$.
\endclaim

\proof{Proof} Since there are not positive dimensional analytic subvarieties, the zero varieties of the ideal sheaves $\cI(h^m)$ must be finite sets of points, hence, by Skoda \cite{Sko72a}, the Lelong numbers $\nu(\ii\Theta_{L,h},x)$ are zero except on a countable set $S\subset X$. By \cite{Dem92}, this implies that $L$ is nef and $c_1(L)^n\ge \sum_{x\in S}\nu(\ii\Theta_{L,h},x)^n$. However, by the Grauert-Riemenschneider conjecture solved in \cite{Siu84}, \cite{Siu85} and \cite{Dem85b}, the positivity of $c_1(L)^n$ would imply that $a(X)=n$ (i.e.\ $X$ Moishezon, a contradiction). Therefore $c_1(L)^n=0$ and $S=\emptyset$.\QED
\endproof

\claim{\label{cohomology-bound}{6.5}.\ Proposition}
Let $X$ be a compact K\"ahler manifold of dimension $n>1$ without any non-trivial subvariety, and with $K_X$ pseudoeffective. Then 
$$h^j(X,K_X^{\otimes m})\leq h^0(X,\Omega^j_X\otimes K_X^{\otimes m})\leq
{n \choose j}\quad\hbox{for every $j\geq 0$},$$
and the Hilbert polynomial $P(m):=\chi(X,K_X^{\otimes m})$ is constant, equal to $\chi(X,\cO_X)$.
\endclaim

\proof{Proof} The inequality $h^j(X,K_X^{\otimes m})\leq h^0(X,\Omega^j_X\otimes K_X^{\otimes m})$ follows from the Hard Lefschetz Theorem~\ref{hard-lefschetz} applied with $L=K_X$ and the corresponding trivial multiplier ideal sheaf. Also,
for any holomorphic vector bundle $E$ on $X$, we have $h^0(X,E)\leq \hbox{\rm rank}(E)$, otherwise, some ratios of determinants of sections would produce a nonconstant meromorphic function, and thus $a(X)>0$, contradiction; here we take
$E= \Omega^j_X\otimes K_X^{\otimes m}$ and get rank$\,E=\smash{n \choose j}$.
The final claim is clear because a polynomial function $P(m)$ which remains bounded as $m\to+\infty$ is necessarily constant.\QED
\endproof

\claim{6.6.\ Corollary} Let $X$ be a strongly simple K\"ahler threefold. Let $h^{i,j}=\dim H^{i,j}(X,\bC)$ be the Hodge numbers. We have
$$c_1(X)^3=c_1(X)\cdot c_2(X)=0,\quad\chi(X,\cO_X)=0\quad\hbox{and}\quad
q:=h^{1,0}>0.$$
\endclaim

\proof{Proof}  The intersection number $K_X^3=-c_1(X)^3$ vanishes because it is the leading term of $P(m)$, up to the factor $3!$. The Riemann-Roch formula then gives 
$$P(m)=\frac{(1-12\,m)}{24}\,c_1(X)\cdot c_2(X).$$
The boundedness of $P(m)$ implies $\chi(X,\cO_X)=\frac{1}{24}c_1(X)\cdot c_2(X)=0$. Now, we write
$$0=\chi(X,\cO_X)=1-h^{1,0}+h^{2,0}-h^{3,0}.$$
By Kodaira's theorem, $h^{2,0}>0$ since $X$ is not projective, and $h^{3,0}\leq 1$ since $a(X)=0$. Thus $0=1-h^{1,0}+h^{2,0}-h^{3,0}\geq 1-q+1-1=1-q$, and $q>0$.\QED
\endproof

Everything is now in place for the final conclusion.

\claim{6.7.\ Theorem} For any strongly simple K\"ahler threefold~$X$, the
Albanese map $\alpha:X\to\Alb(X)$ is a biholomorphism of $3$-dimensional
tori.
\endclaim

\proof{Proof} Since $q=h^{1,0}>0$, the Albanese map $\alpha$ is non constant.
By simplicity, $X$ cannot possess any fibration with positive dimensional 
fibers, so we must have $\dim\alpha(X)=\dim X=3$, and as 
$q=h^{1,0}=h^0(X,\Omega^1_X)\le 3$ (Prop.~\ref{cohomology-bound} with 
$j=1$, $m=0$) the Albanese map $\alpha$ must be surjective. The function
$\det(d\alpha)$ cannot vanish, other\-wise we would get a non trivial divisor,
so $\alpha$ is \'etale. Therefore $X$ is a $3$-dimensional torus, 
as a finite \'etale cover of the $3$-dimensional torus $\Alb(X)$, and $\alpha$
must be an isomorphism.\QED
\endproof

In \cite{CHP14}, the following stronger result is established as a consequence
of the existence of good minimal models for K\"ahler threefolds:

\claim{\label{theorem-kummer}{6.8}.\ Theorem} Let $X$ be smooth compact K\"ahler threefold. If $X$ is simple, there exists a bimeromorphic morphism $X \rightarrow T/G$ where $T$ is a torus and $G$ a finite group acting on $T$. 
\endclaim

\vskip7mm

\noindent
{\twelvebf References}
\medskip
{\eightpoint\let\em=\it
\parskip=3pt plus 1pt minus 1pt

\bibitem[Bay82]{Bay82}
David~Allen Bayer.
\newblock {\em The division algorithm and the {H}ilbert scheme}.
\newblock Thesis (Ph.D.)--Harvard University, ProQuest LLC, Ann Arbor, MI, 1982.
\endbib

\bibitem[BEGZ10]{BEGZ}
S{\'e}bastien Boucksom, Philippe Eyssidieux, Vincent Guedj, and Ahmed Zeriahi.
\newblock Monge-{A}mp\`ere equations in big cohomology classes.
\newblock {\em Acta Math.}, 205(2):199--262, 2010.
\endbib

\bibitem[{Ber}13]{Bern13}
B.~{Berndtsson}.
\newblock {The openness conjecture for plurisubharmonic functions}.
\newblock {\em ArXiv e-prints, math.CV/ 1305.5781}, May 2013.
\endbib

\bibitem[BFJ08]{BFJ08}
S{\'e}bastien Boucksom, Charles Favre, and Mattias Jonsson.
\newblock Valuations and plurisubharmonic singularities.
\newblock {\em Publ. Res. Inst. Math. Sci.}, 44(2):449--494, 2008.
\endbib

\bibitem[BM87]{BM87}
E.~Bierstone and P.~D. Milman.
\newblock Relations among analytic functions. {I}.
\newblock {\em Ann. Inst. Fourier (Grenoble)}, 37(1):187--239, 1987.
\endbib

\bibitem[BM89]{BM89}
Edward Bierstone and Pierre~D. Milman.
\newblock Uniformization of analytic spaces.
\newblock {\em J. Amer. Math. Soc.}, 2(4):801--836, 1989.
\endbib

\bibitem[Bou02]{Bou02}
S{\'e}bastien Boucksom.
\newblock C{\^o}nes positifs des vari\'et\'es complexes compactes.
\newblock PhD Thesis, Grenoble, 2002.
\endbib

\bibitem[Bou04]{Bou04}
S{\'e}bastien Boucksom.
\newblock Divisorial {Z}ariski decompositions on compact complex manifolds.
\newblock {\em Ann. Sci. \'Ecole Norm. Sup. (4)}, 37(1):45--76, 2004.
\endbib

\bibitem[Bru10]{Bru10}
Marco Brunella.
\newblock Uniformisation of foliations by curves.
\newblock In {\em Holomorphic dynamical systems}, volume 1998 of {\em Lecture
  Notes in Math.}, pages 105--163. Springer, Berlin, 2010, also available from
  arXiv e-prints, math.CV/0802.4432\ignore{2010}.
\endbib

\bibitem[BT76]{BT76}
Eric Bedford and B.~A. Taylor.
\newblock The {D}irichlet problem for a complex {M}onge-{A}mp\`ere equation.
\newblock {\em Invent. Math.}, 37(1):1--44, 1976.
\endbib

\bibitem[Cam80]{Cam80}
Fr{\'e}d{\'e}ric Campana.
\newblock Alg\'ebricit\'e et compacit\'e dans l'espace des cycles d'un espace
  analytique complexe.
\newblock {\em Math. Ann.}, 251(1):7--18, 1980.
\endbib

\bibitem[Cam85]{Cam85}
F.~Campana.
\newblock R\'eduction d'{A}lbanese d'un morphisme propre et faiblement
  k\"ahl\'erien. {II}. {G}roupes d'automorphismes relatifs.
\newblock {\em Compositio Math.}, 54(3):399--416, 1985.
\endbib

\bibitem[Cao13]{JC13}
Junyan Cao.
\newblock Th{\'e}or{\`e}mes d'annulation et th{\'e}or{\`e}mes de structure sur
  les vari{\'e}t{\'e}s k{\"a}hleriennes compactes.
\newblock PhD Thesis, Grenoble, 2013.
\endbib

\bibitem[Cao14]{JC14}
Junyan Cao.
\newblock Numerical dimension and a {K}awamata-{V}iehweg-{N}adel type vanishing
  theorem on compact {K}\"ahler manifolds.
\newblock {\em Compositio Math.}, \/(to appear, cf.\ {\em arXiv e-prints:
  math.AG/1210.5692}), 2014.
\endbib

\bibitem[Cao15]{JC15}
Junyan Cao.
\newblock Additivity of the approximation functional of currents induced by
  {B}ergman kernels.
\newblock {\em C. R. Math. Acad. Sci. Paris}, \/(to appear, cf.\ doi:
  10.1016/j.crma.2014.11.004, {\em arXiv e-prints: math.AG/1410.8288}), 2015.
\endbib

\bibitem[CDV13]{CDV13}
{F}r\'ed\'eric {C}ampana, Jean-Pierre {D}emailly, and Misha Verbitsky.
\newblock {Compact {K}\"ahler 3-manifolds without non-trivial subvarieties}.
\newblock {\em ArXiv e-prints, math.CV/1305.5781}, April 2013.
\endbib

\bibitem[CHP14]{CHP14}
Fr{\'e}d{\'e}ric Campana, Andreas H{\"o}ring, and Thomas Peternell.
\newblock Abundance for {K}\"ahler threefolds.
\newblock {\em ArXiv e-prints, math.AG/1403.3175}, 2014.
\endbib

\bibitem[DEL00]{DEL00}
Jean-Pierre Demailly, Lawrence Ein, and Robert Lazarsfeld.
\newblock A subadditivity property of multiplier ideals.
\newblock {\em Michigan Math. J.}, 48:137--156, 2000,
Dedicated to William Fulton on the occasion of his 60th birthday.
\endbib

\bibitem[Dem82]{Dem82}
Jean-Pierre Demailly.
\newblock Estimations {$L^{2}$} pour l'op\'erateur {$\bar \partial $} d'un
  fibr\'e vectoriel holomorphe semi-positif au-dessus d'une vari\'et\'e
  k\"ahl\'erienne compl\`ete.
\newblock {\em Ann. Sci. \'Ecole Norm. Sup. (4)}, 15(3): 457--511, 1982.
\endbib

\bibitem[Dem85]{Dem85b}
Jean-Pierre Demailly.
\newblock Champs magn\'etiques et in\'egalit\'es de {M}orse pour la
  {$d''$}-cohomologie.
\newblock {\em Ann. Inst. Fourier (Grenoble)}, 35(4):189--229, 1985.
\endbib

\bibitem[Dem92]{Dem92}
Jean-Pierre Demailly.
\newblock Regularization of closed positive currents and intersection theory.
\newblock {\em J. Algebraic Geom.}, 1(3):361--409, 1992.
\endbib

\bibitem[Dem93]{Dem93b}
Jean-Pierre Demailly.
\newblock A numerical criterion for very ample line bundles.
\newblock {\em J. Differential Geom.}, 37(2):323--374, 1993.
\endbib

\bibitem[Dem01]{Dem00}
Jean-Pierre Demailly.
\newblock Multiplier ideal sheaves and analytic methods in algebraic geometry.
\newblock In {\em School on {V}anishing {T}heorems and {E}ffective {R}esults in
  {A}lgebraic {G}eometry ({T}rieste, 2000)}, volume~6 of {\em ICTP Lect.
  Notes}, pages 1--148. Abdus Salam Int. Cent. Theoret. Phys., Trieste, 2001.
\endbib

\bibitem[Dem12]{Dembook}
Jean-Pierre Demailly.
\newblock Complex analytic and differential geometry, {\rm online book
  at:}\hfil\break
  {\tt ht{}tp:/{}/www-fourier.ujf-grenoble.fr/\~{}demailly/documents.html}.
\newblock Institut Fourier, 2012.
\endbib

\bibitem[dFEM03]{FEM03}
Tommaso de~Fernex, Lawrence Ein, and Mircea Musta{\c{t}}{\u{a}}.
\newblock Bounds for log canonical thresholds with applications to birational
  rigidity.
\newblock {\em Math. Res. Lett.}, 10(2-3):219--236, 2003.
\endbib

\bibitem[dFEM10]{FEM10}
Tommaso de~Fernex, Lawrence Ein, and Mircea Musta{\c{t}}{\u{a}}.
\newblock Shokurov's {ACC} conjecture for log canonical thresholds on smooth
  varieties.
\newblock {\em Duke Math. J.}, 152(1):93--114, 2010.
\endbib

\bibitem[DK01]{DK01}
Jean-Pierre Demailly and J{\'a}nos Koll{\'a}r.
\newblock Semi-continuity of complex singularity exponents and
  {K}\"ahler-{E}instein metrics on {F}ano orbifolds.
\newblock {\em Ann. Sci. \'Ecole Norm. Sup. (4)}, 34(4):525--556, 2001.
\endbib

\bibitem[DP03]{DP02}
Jean-Pierre Demailly and Thomas Peternell.
\newblock A {K}awamata-{V}iehweg vanishing theorem on compact {K}\"ahler
  manifolds.
\newblock In {\em Surveys in differential geometry, {V}ol.\ {VIII} ({B}oston,
  {MA}, 2002)}, Surv. Differ. Geom., VIII, pages 139--169. Int. Press,
  Somerville, MA, 2003.
\endbib

\bibitem[DP14]{DH14}
Jean-Pierre Demailly and Ho{\`a}ng~Hi{\cfudot{e}}p Ph{\d{a}}m.
\newblock A sharp lower bound for the log canonical threshold.
\newblock {\em Acta Math.}, 212(1):1--9, 2014.
\endbib

\bibitem[DPS94]{DPS94}
Jean-Pierre Demailly, Thomas Peternell, and Michael Schneider.
\newblock Compact complex manifolds with numerically effective tangent bundles.
\newblock {\em J. Algebraic Geom.}, 3(2):295--345, 1994.
\endbib

\bibitem[DPS01]{DPS01}
Jean-Pierre Demailly, Thomas Peternell, and Michael Schneider.
\newblock Pseudo-effective line bundles on compact {K}\"ahler manifolds.
\newblock {\em Internat. J. Math.}, 12(6):689--741, 2001.
\endbib

\bibitem[Eis95]{Eis95}
David Eisenbud.
\newblock {\em Commutative algebra with a view toward algebraic geometry},
  volume 150 of {\em Graduate Texts in Mathematics}.
\newblock Springer-Verlag, New York, 1995.
\endbib

\bibitem[Eno93]{Eno93}
Ichiro Enoki.
\newblock Strong-{L}efschetz-type theorem for semi-positive line bundles over
  compact {K}\"ahler manifolds.
\newblock In {\em Geometry and global analysis ({S}endai, 1993)}, pages
  211--212. Tohoku Univ., Sendai, 1993.
\endbib

\bibitem[FJ05]{FJ05}
Charles Favre and Mattias Jonsson.
\newblock Valuations and multiplier ideals.
\newblock {\em J. Amer. Math. Soc.}, 18(3):655--684 (electronic), 2005.
\endbib

\bibitem[Gal79]{Gal79}
Andr{\'e} Galligo.
\newblock Th\'eor\`eme de division et stabilit\'e en g\'eom\'etrie analytique
  locale.
\newblock {\em Ann. Inst. Fourier (Grenoble)}, 29(2):vii, 107--184, 1979.
\endbib

\bibitem[GZ13]{GZ13}
Qi'an Guan and Xiangyu Zhou.
\newblock Strong openness conjecture for plurisubharmonic functions.
\newblock {\em ArXiv e-prints, math.CV/1311.3781}, 2013.
\endbib

\bibitem[GZ14a]{GZ14a}
Qi'an Guan and Xiangyu Zhou.
\newblock {\kern0pt}strong openness conjecture and related problems for
  plurisubharmonic functions.
\newblock {\em ArXiv e-prints, math.CV/1401.7158}, 2014.
\endbib

\bibitem[GZ14b]{GZ14b}
Qi'an Guan and Xiangyu Zhou.
\newblock Effectiveness of demailly's strong openness conjecture and related
  problems.
\newblock {\em ArXiv e-prints, math.CV/1403.7247}, 2014.
\endbib

\bibitem[Hir64]{Hir64}
Heisuke Hironaka.
\newblock Resolution of singularities of an algebraic variety over a field of
  characteristic zero. {I}, {II}.
\newblock {\em Ann. of Math. (2) 79 (1964), 109--203; ibid. (2)}, 79:205--326,
  1964.
\endbib

\bibitem[H{\"o}r66]{Hor66}
Lars H{\"o}rmander.
\newblock {\em An introduction to complex analysis in several variables},
  volume~7 of {\em North-Holland Mathematical Library}.
\newblock North-Holland Publishing Co., Amsterdam, third edition, 1990, (first
  edition 1966).
\endbib

\bibitem[HPR11]{HPR11}
Andreas H{\"o}ring, Thomas Peternell, and Ivo Radloff.
\newblock Uniformisation in dimension four: towards a conjecture of iitaka.
\newblock {\em arXiv e-prints, math.AG/1103.5392, to appear in Math.
  Zeitschrift}, 2011.
\endbib

\bibitem[JM12]{JM12}
Mattias Jonsson and Mircea Musta{\c{t}}{\u{a}}.
\newblock Valuations and asymptotic invariants for sequences of ideals.
\newblock {\em Ann. Inst. Fourier (Grenoble)}, 62(6):2145--2209 (2013), 2012.
\endbib

\bibitem[JM14]{JM14}
Mattias Jonsson and Mircea Musta{\c{t}}{\u{a}}.
\newblock An algebraic approach to the openness conjecture of {D}emailly and
  {K}oll\'ar.
\newblock {\em J. Inst. Math. Jussieu}, 13(1):119--144, 2014.
\endbib

\bibitem[Kim14]{Kim13}
Dano Kim.
\newblock A remark on the approximation of plurisubharmonic functions.
\newblock {\em C. R. Math. Acad. Sci. Paris}, \/(to appear), 2014.
\endbib

\bibitem[Laz04]{Laz04}
Robert Lazarsfeld.
\newblock {\em Positivity in algebraic geometry. {I,II}}, volume~48 of {\em
  Ergebnisse der Mathematik und ihrer Grenzgebiete.}
\newblock Springer-Verlag, Berlin, 2004.
\endbib

\bibitem[Mo{\u\i}66]{Moi66}
B.~G. Mo{\u\i}{\v{s}}ezon.
\newblock On {$n$}-dimensional compact complex manifolds having {$n$}
  algebraically independent meromorphic functions. {I, II, III}.
\newblock {\em Izv. Akad. Nauk SSSR Ser. Mat.}, 30:133--174, 345--386,
  621--656, 1966.
\endbib

\bibitem[Mou95]{Mou95}
Christophe Mourougane.
\newblock Versions k\"ahl\'eriennes du th\'eor\`eme d'annulation de
  {B}ogomolov-{S}ommese.
\newblock {\em C. R. Acad. Sci. Paris S\'er. I Math.}, 321(11):1459--1462,
  1995.
\endbib

\bibitem[Nad90]{Nad90}
Alan~Michael Nadel.
\newblock Multiplier ideal sheaves and {K}\"ahler-{E}instein metrics of
  positive scalar curvature.
\newblock {\em Ann. of Math. (2)}, 132(3):549--596, 1990.
\endbib

\bibitem[OT87]{OT87}
Takeo Ohsawa and Kensh{\=o} Takegoshi.
\newblock On the extension of {$L^2$} holomorphic functions.
\newblock {\em Math. Z.}, 195(2):197--204, 1987.
\endbib

\bibitem[Pet86]{Pet86}
Thomas Peternell.
\newblock Algebraicity criteria for compact complex manifolds.
\newblock {\em Math. Ann.}, 275(4): 653--672, 1986.
\endbib

\bibitem[Pet98]{Pet98a}
Thomas Peternell.
\newblock Moishezon manifolds and rigidity theorems.
\newblock {\em Bayreuth. Math. Schr.}, \/(54): 1--108, 1998.
\endbib

\bibitem[Ph{\d{a}}14]{Hiep14}
Ho{\`a}ng~Hi{\cfudot{e}}p Ph{\d{a}}m.
\newblock The weighted log canonical threshold.
\newblock {\em C. R. Math. Acad. Sci. Paris}, 352(4):283--288, 2014.
\endbib

\bibitem[PS00]{PhS00}
D.~H. Phong and Jacob Sturm.
\newblock On a conjecture of {D}emailly and {K}oll\'ar.
\newblock {\em Asian J. Math.}, 4(1):221--226, 2000.
\endbib

\bibitem[Ric68]{Ric68}
Rolf Richberg.
\newblock Stetige streng pseudokonvexe {F}unktionen.
\newblock {\em Math. Ann.}, 175:257--286, 1968.
\endbib

\bibitem[Siu74]{Siu74}
Yum~Tong Siu.
\newblock Analyticity of sets associated to {L}elong numbers and the extension
  of closed positive currents.
\newblock {\em Invent. Math.}, 27:53--156, 1974.
\endbib

\bibitem[Siu84]{Siu84}
Yum~Tong Siu.
\newblock A vanishing theorem for semipositive line bundles over non-{K}\"ahler
  manifolds.
\newblock {\em J. Differential Geom.}, 19(2):431--452, 1984.
\endbib

\bibitem[Siu85]{Siu85}
Yum~Tong Siu.
\newblock Some recent results in complex manifold theory related to vanishing
  theorems for the semipositive case.
\newblock In {\em Workshop {B}onn 1984 ({B}onn, 1984)}, volume 1111 of {\em
  Lecture Notes in Math.}, pages 169--192. Springer, Berlin, 1985.
\endbib

\bibitem[Sko72a]{Sko72a}
Henri Skoda.
\newblock Application des techniques {$L^{2}$} \`a la th\'eorie des id\'eaux
  d'une alg\`ebre de fonctions holomorphes avec poids.
\newblock {\em Ann. Sci. \'Ecole Norm. Sup. (4)}, 5:545--579, 1972.
\endbib

\bibitem[Sko72b]{Sko72b}
Henri Skoda.
\newblock Sous-ensembles analytiques d'ordre fini ou infini dans {${\Bbb
  C}^{n}$}.
\newblock {\em Bull. Soc. Math. France}, 100:353--408, 1972.
\endbib

\bibitem[Tak97]{Tak97}
Kensho Takegoshi.
\newblock On cohomology groups of nef line bundles tensorized with multiplier
  ideal sheaves on compact {K}\"ahler manifolds.
\newblock {\em Osaka J. Math.}, 34(4):783--802, 1997.
\endbib

\bibitem[{Tsu}07]{Tsu07}
H.~{Tsuji}.
\newblock {Extension of log pluricanonical forms from subvarieties}.
\newblock {\em ArXiv e-prints}, math.AG/ 0709.2710, September 2007.
\endbib

\bibitem[Yau78]{Yau78}
Shing~Tung Yau.
\newblock On the {R}icci curvature of a compact {K}\"ahler manifold and the
  complex {M}onge-{A}mp\`ere equation. {I}.
\newblock {\em Comm. Pure Appl. Math.}, 31(3):339--411, 1978.
\endbib
\medskip
}

\medskip\medskip
\noindent
Jean-Pierre Demailly\\
Universit\'e de Grenoble I, Institut Fourier, UMR 5582 du CNRS\\
BP$\,$74, 100 rue des Maths, 38402 Saint-Martin d'H\`eres, France\\
{\it e-mail}: {\tt jean-pierre.demailly@ujf-grenoble.fr}

\vskip30pt
\noindent
(December 31, 2014; printed on \today)

\ifuseauxfile \closeout\outauxfile \fi
\bye